\documentclass[12pt]{article}
\usepackage{graphicx}
\usepackage{float}
\usepackage{amssymb,amsmath,amsfonts,relsize,mathrsfs}
\usepackage[usenames]{color}
\usepackage{cite}
\usepackage{esvect}
\usepackage{txfonts}
\usepackage{cmdtrack}

\makeatletter
\renewcommand{\fnum@figure}{Fig. \thefigure}
\makeatother
\newtheorem{theorem}{Theorem}[section]
\newtheorem{definition}[theorem]{Definition}
\newtheorem{corollary}[theorem]{Corollary}
\newtheorem{proposition}[theorem]{Proposition}
\newtheorem{remark}[theorem]{Remark}
\newtheorem{lemma}[theorem]{Lemma}

\newtheorem{problem}[theorem]{Problem}

\newtheorem{claim}[theorem]{Claim}
\newtheorem{condition}[theorem]{Condition}
\newtheorem{pra}[theorem]{The Projection Algorithm}
\newtheorem{alg}[theorem]{Algorithm}
\newcommand {\Ac}      {{\mathcal A}}
\newcommand {\Bc}      {{\mathcal B}}

\newcommand {\Fc}      {{\mathcal F}}

\newcommand {\Hc}      {{\mathcal H}}
\newcommand {\Ic}      {{\mathcal I}}

\newcommand {\Kc}      {{\mathcal K}}
\newcommand {\Lc}      {{\mathcal L}}
\newcommand {\Mc}      {{\mathcal M}}

\newcommand {\Pc}      {{\mathcal P}}

\newcommand {\Rc}      {{\mathcal R}}

\newcommand {\Tc}      {{\mathcal T}}

\newcommand {\Vc}      {{\mathcal V}}
\newcommand {\Wc}      {{\mathcal W}}

\newcommand {\Acr}      {{\mathscr A}}
\newcommand {\Bcr}      {{\mathscr B}}

\newcommand {\Fcr}      {{\mathscr F}}
\newcommand {\Gcr}      {{\mathscr G}}
\newcommand {\Hcr}      {{\mathscr H}}
\newcommand {\R}       {{\bf R}}
\newcommand {\RPL}     {{\bf R_+}}

\newcommand {\RN}      {\R^n}
\newcommand {\RD}      {\R^D}

\newcommand {\RT}      {\R^2}
\newcommand {\RTP}     {\R^2_+}

\newcommand {\lh}      {{\bf h}}

\newcommand {\tF}      {\widetilde{F}}

\newcommand {\tK}      {\widetilde{K}}

\newcommand {\tQ}      {\widetilde{Q}}

\newcommand {\tT}      {\widetilde{T}}

\newcommand {\tMc}     {\widetilde{{\mathcal M}}}

\newcommand {\tMf}     {\widetilde{\mathfrak M}}
\newcommand {\ta}      {\tilde{a}}

\newcommand {\tf}      {\tilde{f}}

\newcommand {\tu}      {\tilde{u}}

\newcommand {\tx}      {\tilde{x}}

\newcommand {\tz}      {\tilde{z}}

\newcommand {\tlm}     {\tilde{\lambda}}
\newcommand {\trh}     {\tilde{\rho}}

\newcommand {\tta}     {\tilde{\tau}}
\newcommand {\tell}    {\tilde{\ell}}

\newcommand {\hM}      {\widehat{M}}
\newcommand {\hx}      {\hat{x}}

\newcommand {\hCF}     {\widehat{\Cf}}
\newcommand {\ellh}    {\hat{\ell}}
\newcommand {\lmh}     {\hat{\lambda}}
\newcommand {\brx}     {\bar{x}}

\newcommand {\Cf}      {{\mathfrak C}}
\newcommand {\Df}      {{\mathfrak D}}

\newcommand {\Mf}      {{\mathfrak M}}

\newcommand {\Rf}      {{\mathfrak R}}

\newcommand {\Tf}      {{\mathfrak T}}

\newcommand {\ve}      {\varepsilon}

\newcommand {\MR}      {(\Mc,\rho)}
\newcommand {\MS}      {\Mf}
\newcommand {\lmv}     {\vec{\lambda}}
\newcommand {\vf}      {\varphi}
\newcommand {\emp}     {\emptyset}

\newcommand {\CRT}     {\Conv(\RT)}
\newcommand {\RCT}     {\Rf(\RT)}
\newcommand {\SO}      {{\bf S}_1}
\newcommand {\RL}      {\Rc_{F}}

\newcommand {\HR}      {\Hc}

\newcommand {\AF}      {\vf_F}

\newcommand {\HY}      {\mathbb{H}}
\newcommand {\HYG}     {\mathbb{H}_G}
\newcommand {\WTC}     {\widetilde{\Cf}}

\newcommand {\Tfw}     {\widetilde{\Tf}}

\newcommand {\cbg}    {\,\mathlarger{\mathlarger{\cap}}\,}

\newcommand {\BXR}     {I_0}

\newcommand {\LTI}     {\ell^2_\infty}

\newcommand {\LTT}     {\ell^2_2}

\newcommand {\KRT}     {\Kc(\RT)}
\newcommand {\CNV}     {\Conv(\RT)}
\newcommand {\DXY}     {\Delta_n(x,y)}
\newcommand {\DXP}     {\Delta_n(x,x')}
\newcommand {\DYP}     {\Delta_n(y,y')}

\newcommand {\HPL}     {\Hc\Pc(\RT)}
\newcommand {\SXP}     {\sin{\AF(x,x')}}
\newcommand {\SYP}     {\sin{\AF(y,y')}}

\newcommand {\ip}[1]   {\langle{#1}\rangle}

\newcommand {\FM}      {|F|_{\Mf}}
\newcommand {\n}       {{\bf n}}

\newcommand {\LF}      {(\lambda;F)}

\newcommand {\PL}      {{\mathscr P}_L(\RT)}

\newcommand {\al}      {\alpha}
\newcommand {\VST}     {\vspace*{1mm}}

\newcommand {\Lip}     {\operatorname{Lip}}

\newcommand {\dhf}     {\operatorname{d_H}}

\newcommand {\Prj}     {\operatorname{\text{\sc Pr}}}
\newcommand {\Prm}     {\operatorname{\text{\bf Pr}}}
\newcommand {\diam}    {\operatorname{diam}}
\newcommand {\dist}    {\operatorname{dist}}

\newcommand {\Conv}    {\operatorname{Conv}}
\newcommand {\sign}    {\operatorname{sign}}
\newcommand {\cent}    {\operatorname{center}}
\newcommand {\cl}      {\operatorname{\bf cl}}
\newcommand {\smsk}    {\smallskip}
\newcommand {\msk}     {\medskip}
\newcommand {\bsk}     {\bigskip}
\newcommand {\bx}      {\hspace{10mm}$\blacksquare$}
\newcommand {\rbx}     {\hspace{10mm}$\blacktriangleleft$}

\newcommand {\nn}      {\nonumber}
\newcommand {\rf}[1]    {(\ref{#1})}      
\newcommand {\reff}[1] {\ref{#1}}         
\newcommand{\lbl}[1]      {\label{#1}}       
\newcommand{\be}          {\begin{eqnarray}}
\newcommand{\bel}[1]      {\begin{eqnarray} \label{#1}}
\newcommand{\ee}           {\end{eqnarray}}
\newcommand {\SECT}[2] {\section*{\centerline{\normalsize
{\bf #1}}} \setcounter{section}{#2}
\setcounter{theorem}{0}\setcounter{equation}{0}}
\begin{document}
\parindent 1em
\parskip 0mm
\centerline{{\bf Existence Criteria for Lipschitz Selections of Set-Valued Mappings in $\RT$}}
\vspace*{2mm}
\centerline{By~ {\sc Pavel Shvartsman}}\vspace*{2 mm}
\centerline {\it Department of Mathematics, Technion - Israel Institute of Technology,}\vspace*{1 mm}
\centerline{\it 32000 Haifa, Israel}\vspace*{1 mm}
\centerline{\it e-mail: pshv@technion.ac.il}
\renewcommand{\thefootnote}{ }
\footnotetext[1]{{\it\hspace{-6mm}Math Subject
Classification:} 46E35
\par This research was supported by the ISRAEL SCIENCE FOUNDATION (grant No. 520/22).}
\begin{abstract} Let $F$ be a set-valued mapping which to each point $x$ of a metric space $(\Mc,\rho)$ assigns a convex closed set $F(x)\subset\RT$. We present several constructive criteria for the existence of a Lipschitz selection of $F$, i.e., a Lipschitz mapping $f:\Mc\to\RT$ such that $f(x)\in F(x)$ for every $x\in\Mc$. The geometric methods we develop to prove these criteria provide efficient algorithms for constructing nearly optimal Lipschitz selections and computing the order of magnitude of their Lipschitz seminorms.
\end{abstract}
\vspace*{-10mm}
\renewcommand{\contentsname}{ }
\tableofcontents
\addtocontents{toc}{{\centerline{\sc{Contents}}}
\vspace*{2mm}\par}

\SECT{1. Introduction.}{1}
\addtocontents{toc}{\hspace*{3.2mm} 1. Introduction.\hfill \thepage\par\VST}
\bsk


\indent\par Let $\Mf=(\Mc,\rho)$ be a {\it pseudometric space}, i.e., suppose that the ``distance function'' $\rho:\Mc\times\Mc\to [0,+\infty]$ satisfies
$$
\rho(x,x)=0,~ \rho(x,y)=\rho(y,x),~~~\text{and}~~~\rho(x,y)\le \rho(x,z)+\rho(z,y)
$$
for all $x,y,z\in\Mc$. Note that $\rho(x,y)=0$ may hold with $x\ne y$, and $\rho(x,y)$ may be $+\infty$.
\par By $\Lip(\Mc)$ we denote the space of all Lipschitz mappings from $\Mc$ into $\RT$ equipped with the Lipschitz seminorm
$$
\|f\|_{\Lip(\Mc)}=\inf\{\,\lambda\ge 0:\|f(x)-f(y)\|
\le\lambda\,\rho(x,y)~~\text{for all}~~x,y\in\Mc\,\}.
$$
\par Hereafter $\|\cdot\|$ denotes the uniform norm in $\RT$, i.e.,
$$
\|x\|=\max\{|x_1|,|x_2|\}~~~~\text{for}~~~ x=(x_1,x_2)\in\RT.
$$
\par  Let $F$ be a set-valued mapping which to each element $x\in\Mc$ assigns a non-empty convex closed set $F(x)\subset\RT$. A {\it selection} of $F$ is a map $f:\Mc\to \RT$ such that $f(x)\in F(x)$ for all $x\in\Mc$.
A selection $f$ is said to be Lipschitz if $f\in\Lip(\Mc)$.
\par We introduce the quantity $\FM$ by letting
\bel{FM}
\FM=\inf\{\,\|f\|_{\Lip(\Mc)}: f~~\text{is a Lipschitz selection of}~~F\}
\ee
whenever $F$ has a Lipschitz selection, and we set $\FM=+\infty$ otherwise.
\par Let $\Tf$ be a family of convex closed subsets of $\RT$. In this paper, we study two main problems related to calculation of the quantity $\FM$ and constructing of a nearly optimal Lipschitz selection of $F$ respectively.
\par Here is the first of these problems.
\begin{problem}\lbl{PR1} {\em  {\it Find an explicit constructive formula for the order of magnitude of the quantity $\FM$ where $F:\Mc\to\Tf$ is an arbitrary set-valued mapping.}
\par More specifically, we would like to find an {\it explicit} formula for computing the value of $\FM$ (up to an absolute positive constant) which exploits only the pseudometric $\rho$ and certain geometric characteristics of the sets $F(x)\in\Tf$, $x\in\Mc$.}
\end{problem}
\par By ``order of magnitude'' we mean the following: Two numbers $A,B\ge 0$ are said to have ``the same order of magnitude'' provided that $cA\le B\le CA$, with absolute positive constants $c$ and $C$. To ``compute the order of magnitude of $A$'' is to compute a number $B$ such that $A$ and $B$ have the same order of magnitude.
\par Let us formulate the second main problem.
\begin{problem}\lbl{PR2} {\em {\it Find a constructive algorithm which given a set-valued mapping $F:\Mc\to\Tf$ assigns a nearly optimal Lipschitz selection $f$ of $F$.}
\par Thus, we are looking for an efficient constructive algorithm which proceeds a selection $f$ of $F$ such that $\|f\|_{\Lip(\Mc)}\le \gamma\,\FM$ where $\gamma\ge 1$ is an absolute constant.
\par We expect that this algorithm uses only the pseudometric $\rho$ and the geometrical parameters which determine the mapping $F$.}
\end{problem}
\smsk
\par In this paper, the family $\Tf$ is one of the following families of convex sets:
\smsk
\par (i) $\Tf=\CRT$ where
\bel{CRT-DF}
\CRT=\{C\subset\RT: C~~\text{is non-empty convex and closed}\};
\ee
\par (ii) $\Tf=\KRT$, where
\bel{KRT-DF}
\KRT=\{C\subset\RT: C~~\text{is non-empty convex and bounded}\};
\ee
\par (iii) $\Tf=\HPL$ where
\bel{HPL-DF}
\HPL=\{H\subset\RT: H~~\text{is a closed half-plane}\}.
\ee
\par In Theorem \reff{FOR-1} and Theorem \reff{FOR-2}  below, we present two solutions to Problem \reff{PR1} by exhibiting two different explicit formulae for the order of magnitude of the quantity $\FM$ where $F$ is an arbitrary set-valued mapping from $\Mc$ into the family $\Tf=\KRT$. These formulae are expressed in terms of the diameters of the four-point subsets of $\Mc$ and the angles between the supporting half-planes of the sets $F(x)$, $x\in\Mc$.
\smsk
\par Problems \reff{PR1} and \reff{PR2} are special cases of the general Lipschitz selection problem studying the existence and properties of Lipschitz selections of set-valued mappings from (pseudo)metric spaces into various families of convex subsets of Banach spaces. The Lipschitz selection problem may be regarded as a search for a Lipschitz mapping that agrees approximately with data.
\par There is an extensive literature devoted to different aspects of the Lipschitz and the related smooth selection problems. Among the multitude of results known so far we mention those in the papers and monographs \cite{Ar-1991,AF-1990,BL-2000,FIL-2017,FP-2019,FS-2018,
JLO-2022-1,JLO-2022-2,PR-1992,PY-1989,PY-1995,S-2001,S-2002,
S-2004,S-2008}. We refer the reader to all of these works and references therein, for numerous results and techniques concerning this topic.
\smsk
\par The Lipschitz selection problem is of great interest in recent years, mainly due to its close connections with the classical {\it Whitney Extension Problem} \cite{Wh-1934}: {\it Given a positive integer $m$ and a function $f$ deﬁned on a closed subset of $\RN$, how can one tell whether $f$ extends to a $C^m$-function on all of $\RN$?}
\smsk
\par Over the years (since 1934) this problem has attracted a lot of attention, and there is an extensive literature devoted to this problem and its analogues for various spaces of smooth functions. For a detailed account of the history of extension and restriction problems for $m$-smooth functions, and various references related to this topic, we refer the reader to
\cite{BS-1994,BS-2001,F-2005,F-2009,FK-2009,
FI-2020,JL-2021,FJL-2023}.
\smsk
\par As an example, let us illustrate the connection between the Lipschitz selection problem and the Whitney problem for the space $C^2(\RN)$.
In \cite{S-1987,BS-2001,S-2002} we show that the Whitney problem for the restrictions of $C^2$-functions to {\it finite} subsets of $\RN$ can be reduced to a certain Lipschitz selection problem for {\it affine-set valued} mappings.
\par A solution to this special case of the Lipschitz selection problem given in \cite{S-2001,S-2002,S-2004} led us to an in\-te\-resting property of the restrictions of  $C^2$-functions called by C. Fefferman \cite{F-2009} (for the general case of $C^m$-spaces) as the {\it Finiteness Principle}. This principle enables us to reduce the Whitney problem for $C^2(\RN)$-restrictions to an {\it arbitrary finite subsets} of $\RN$ to a similar problem but for $C^2(\RN)$-restrictions to {\it finite sets consisting of at most $k^{\#}=3\cdot 2^{n-1}$ points}. See \cite{S-1987}. In \cite{F-2005}, C. Fefferman showed that this version of the Finiteness Principle holds for the space $C^m(\RN)$ for arbitrary $m,n\ge 1$ with a certain constant $k^{\#}=k^{\#}(m,n)$ depending only on $m$ and $n$.
\smsk
\par Furthermore, in \cite{S-2002} we solved Problem \reff{PR1} for the special case of the {\it line-set valued} mappings in $\RT$, and showed how constructive geometrical criteria for Lipschitz selections of such mappings are transformed into purely analytical descriptions of the restrictions of $C^2$-functions to finite subsets of the plane.
\smsk
\par There is also a {\it Finiteness Principle for Lipschitz selections} proven in the recent joint work with C. Fefferman \cite{FS-2018}. In two dimensional case, this principle states the following (see \cite{S-2002}):
\par {\it Let $F:\Mc\to\KRT$ be a set-valued mapping. If for every $\Mc'\subset\Mc$ consisting of at most \underline{four} points, the restriction $F|_{\Mc'}$ of $F$ to $\Mc'$ has a Lipschitz selection $f_{\Mc'}$ with $\|f_{\Mc'}\|_{\Lip(\Mc')}\le 1$, then $F$ has a Lipschitz selection with Lipschitz seminorm at most $\gamma$. Here, $\gamma>0$ is an absolute constant.}
\smsk
\par This statement is equivalent to the following inequality:
\bel{FP-FM}
\FM\le \gamma\,\sup\{|F|_{\Mc'}|_{\Mf'}: \Mf'=(\Mc',\rho), \Mc'\subset\Mc, \#\Mc'\le N\}~~~~\text{with}~~~N=4.
\ee
(Clearly, the converse inequality trivially holds with $\gamma=1$.) It is shown in \cite{FS-2018} that the Finiteness Principle holds for set-valued mappings taking values in the family $\Kc(\RN)$ of all compact convex subsets of $\RN$ with $N=2^n$ and a constant $\gamma=\gamma(n)$ depending only on $n$.
\smsk
\par In particular, the Finiteness Principle enables us to reduce Problem \reff{PR1} to the same one but for a {\it finite} pseudometric spaces $\Mf=(\Mc,\rho)$ consisting of at most {\it four} elements. Nevertheless, as we will see below, even for a four-element pseudometric space, the solution to Problem \reff{PR1} (especially with good lower and upper bounds for the quantity $\FM$) remains a rather difficult problem.
\smsk

\smsk
\par Our main results, Theorems \reff{FOR-1} and \reff{FOR-2}, are corollaries of Theorems \reff{CR-1L} and \reff{CF-CR} which provide two different solutions to Problem \reff{PR1} for the case of set-valued mappings taking values in the family $\Tf=\HPL$ of all closed half-planes in $\RT$.
\par  Let us prepare the ingredients that are needed to formulate these results. Let $\SO$ be the unit circle in $\RT$, and let
$$
\n:\Mc\to\SO~~~~\text{and}~~~~\al:\Mc\to\R
$$
be two mappings defined on $\Mc$. These mappings determine a set-valued mapping $F:\Mc\to \HPL$ defined by
\bel{F-NAL-I}
F(x)=\{a\in\RT:\ip{a,\n(x)}+\al(x)\le 0\},~~~~x\in\Mc.
\ee
Here, given $a=(a_1,a_2)$, $\n(x)=(n_1(x),n_2(x))\in\RT$, by $$\ip{a,\n(x)}=a_1 n_1(x)+a_2 n_2(x)$$ we denote the standard inner product in $\RT$. Thus, for each $x\in\Mc$, the set $F(x)$ is a half-plane in $\RT$ whose boundary is the straight line
$$
\ell_F(x)=\{a\in\RT:\ip{a,\n(x)}+\al(x)=0\}.
$$
The unit vector $\n(x)$ is directed outside of the half-plane $F(x)$ and is orthogonal to the line  $\ell_F(x)$.
\begin{definition}\lbl{G-PS} {\em Let $F:\Mc\to\HPL$ be a set-valued mapping defined by \rf{F-NAL-I}. We say that the half-planes $\{F(x):x\in\Mc\}$ are in {\it general position} if there exist elements $x_1,...,x_m\in\Mc$ such that}
$$
\text{\it the interior of the convex hull of the vectors} ~~\n(x_1),...,\n(x_m)~~\text{\it contains}~~0.~~~
\text{\rbx}
$$
\end{definition}
\par Given $x,y\in\Mc$, we let $\DXY$ denote the determinant
\bel{DET-XY}
\DXY=\left|
\begin{array}{ll}
n_1(x)& n_1(y)\vspace{2mm}\\
n_2(x)& n_2(y)
\end{array}
\right|=n_1(x)\,n_2(y)-n_2(x)\,n_1(y).
\ee
\par Next, let $\n(x),\n(y)\in \SO$ be two non-collinear vectors (we write $\n(x)\nparallel \n(y)$). Let
$$
w(x,y:F)=(w_1(x,y:F),w_2(x,y:F))=\ell_F(x)\cap\ell_F(y)
$$
be the point of intersection of the straight lines $\ell_F(x)$ and $\ell_F(y)$. Clearly,
$$
w_1(x,y:F)=\frac{\left|
\begin{array}{ll}
\al(y)& \al(x)\\
n_2(y)& n_2(x)
\end{array}
\right|}{\DXY}
~~~~~\text{and}~~~~~
w_2(x,y:F)=\frac{\left|
\begin{array}{ll}
\al(x)& \al(y)\\
n_1(x)& n_1(y)
\end{array}
\right|}{\DXY}.
$$
\par Finally, given $x,x',y,y'\in \Mc$ such that $\n(x)\nparallel \n(x')$ and $\n(y)\nparallel \n(y')$, we introduce the following quantities:
\bel{D1-DF}
D_1[x,x':y,y']=
\frac{\rho(x,x')}{|\DXP|}\min\{|n_2(x)|,|n_2(x')|\}+
\frac{\rho(y,y')}{|\DYP|}\min\{|n_2(y)|,|n_2(y')|\}+
\rho(x,y)
\ee
and
\bel{D2-DF}
D_2[x,x':y,y']=
\frac{\rho(x,x')}{|\DXP|}\min\{|n_1(x)|,|n_1(x')|\}+
\frac{\rho(y,y')}{|\DYP|}\min\{|n_1(y)|,|n_1(y')|\}+
\rho(x,y).
\ee

\begin{theorem}\lbl{CR-1L} Let $F:\Mc\to\HPL$ be a set-valued mapping defined by \rf{F-NAL-I}. Assume that either $\Mc$ is finite or the half-planes $\{F(x):x\in\Mc\}$ are in general position.
(See Definition \reff{G-PS}.)
\par Then $F$ has a Lipschitz selection if and only if there exists a constant $\lambda\ge 0$ such that the following conditions hold:
\msk
\par ($\bigstar 1$) $\al(x)+\al(y)\le \lambda\,\rho(x,y)$ for every $x,y\in\Mc$ such that $\n(y)=-\n(x)$;
\smsk
\par ($\bigstar 2$) Let $x,x',y,y'\in\Mc$ be arbitrary elements such that $\n(x)\nparallel \n(x')$ and $\n(y)\nparallel \n(y')$. Then the following conditions are satisfied:
\smsk
\par (i). If
\bel{N2X}
n_2(x)\,n_2(x')\le 0,~n_1(x)+n_1(x')\le 0~~~~~\text{and}~~~~~n_2(y)\,n_2(y')\le 0, ~n_1(y)+n_1(y')\ge 0,
\ee
then
$$
w_1(x,x':F)-w_1(y,y':F)\le \lambda\,D_1[x,x':y,y'];
$$
\par (ii). If
\bel{N1X}
n_1(x)\,n_1(x')\le 0,~ n_2(x)+n_2(x')\le 0,~~~~~\text{and}~~~~~n_1(y)\,n_1(y')\le 0,~ n_2(y)+n_2(y')\ge 0.
\ee
then
$$
w_2(x,x':F)-w_2(y,y':F)\le \lambda\,D_2[x,x':y,y'].
$$
\par Furthermore,
$$
\tfrac{1}{\sqrt{2}}\,\inf\lambda
\le \,|F|_{\Mf} \le 5\inf\lambda.
$$
See Fig. 1.
\end{theorem}

\begin{figure}[H]
\hspace{17mm}
\includegraphics[scale=0.33]{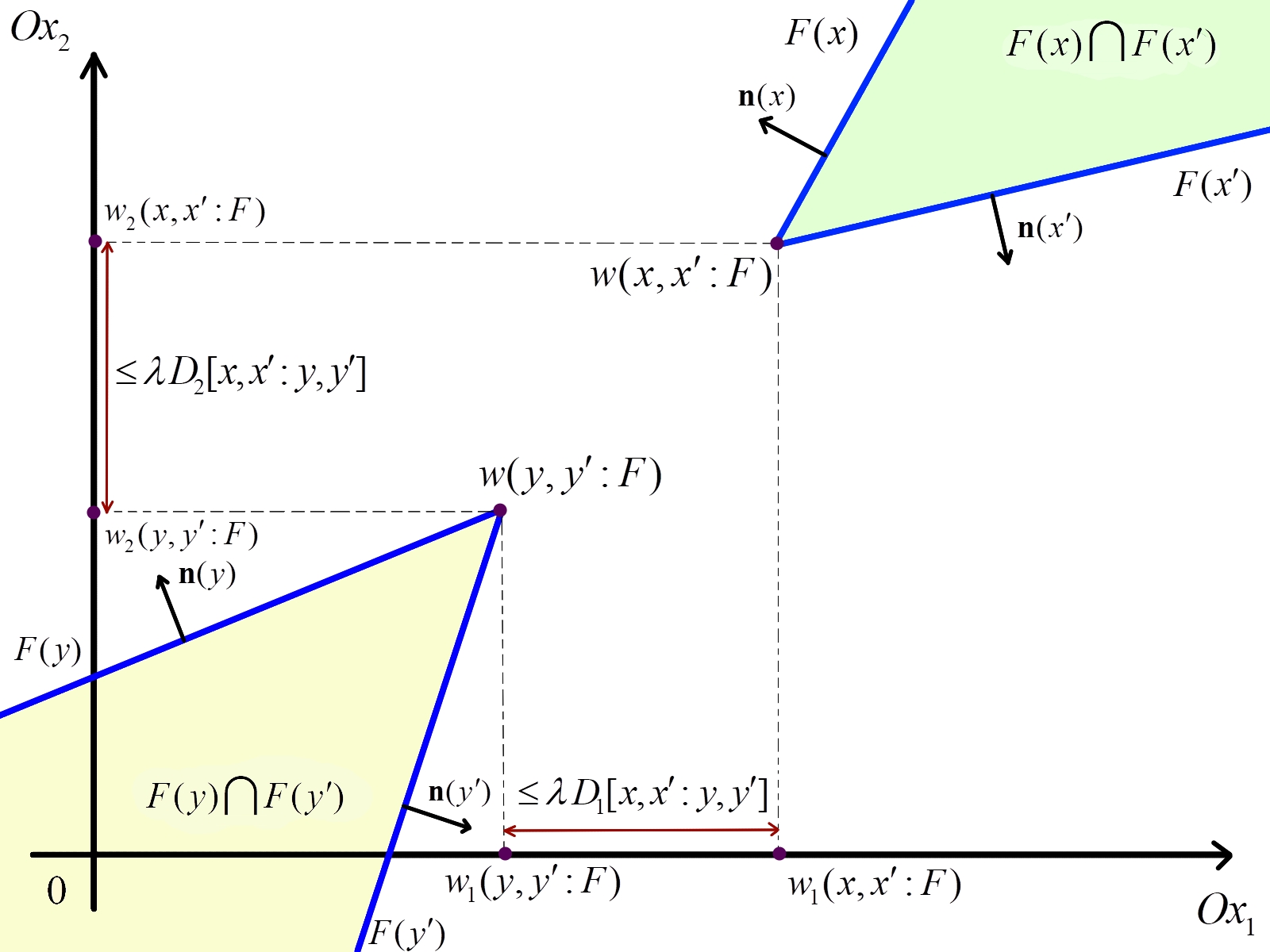}
\caption{The Lipschitz selection criterion for half-planes.}
\end{figure}
\msk
\begin{remark} {\em Let us explain the geometrical meaning of inequalities \rf{N2X} and \rf{N1X}.
\par We let $\Prj_i$, $i=1,2$, denote the operator of the orthogonal projection onto the axis $Ox_i$. It can be shown that, if these inequalities hold then $\Prj_i[F(x)\cbg F(x')]$ and $\Prj_i[F(y)\cbg F(y')]$ are two closed {\it unbounded} intervals in $Ox_i$ such that {\it the first of them is bounded from below}, while {\it the second is bounded from above} (with respect to the standard ordering on the coordinate axes). In other words, the inequalities in \rf{N2X} are equivalent to the following equalities:
$$
\Prj_1[F(x)\cbg F(x')]=\{t\,(w_1(x,x':F),0):t\ge 1\},~
\Prj_1[F(y)\cbg F(y')]=\{t\,(w_1(y,y':F),0):t\le 1\}.
$$
In turn, the inequalities in \rf{N1X} mean the following:
$$
\Prj_2[F(x)\cbg F(x')]=\{t\,(0,w_2(x,x':F)):t\ge 1\},~
\Prj_2[F(y)\cbg F(y')]=\{t\,(0,w_2(y,y':F)):t\le 1\}.
$$
\par Thus, if \rf{N2X} and \rf{N1X} hold, condition ($\bigstar 2$) of Theorem \reff{CR-1L} can be reformulated as follows:
\smsk
\par ($\bigstar 2'$) Let $x,x',y,y'\in\Mc$ be arbitrary elements such that $\n(x)\nparallel \n(x')$ and $\n(y)\nparallel \n(y')$. Then for every $i=1,2$ the following inequality
\bel{CR1-PR}
\dist(\Prj_i[F(x)\cbg F(x')],\Prj_i[F(y)\cbg F(y')])
\le \lambda\,D_i[x,x':y,y']
\ee
holds. See Fig. 2.

\begin{figure}[H]
\hspace{22mm}
\includegraphics[scale=0.3]{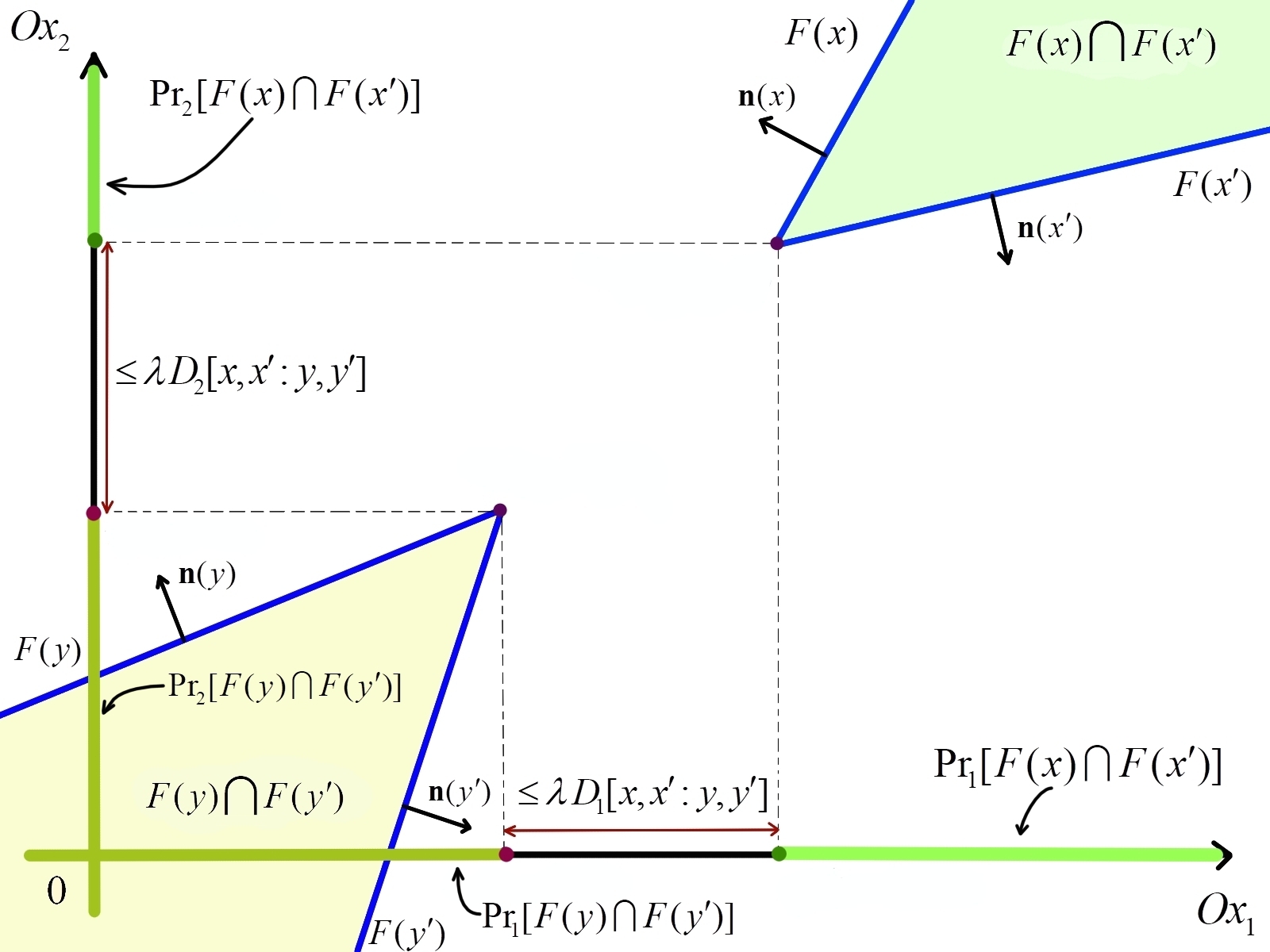}
\caption{A reformulation of the Lipschitz selection criterion for half-planes.}
\end{figure}
\par In general, we can omit the requirement that inequalities \rf{N2X} and \rf{N1X} are satisfied. In other words, in the formulation of Theorem \reff{CR-1L} we can replace condition ($\bigstar 2$) with ($\bigstar 2'$).\rbx}
\end{remark}
\par We turn to the second criterion for the seminorm of a nearly optimal Lipschitz selection. This criterion can be regarded as a certain modification of the criterion \rf{CR1-PR}. This modification is motivated by the following observation: {\it Theorem \reff{CR-1L} includes the quantities $D_1[\cdot,\cdot:\cdot,\cdot]$ and $D_2[\cdot,\cdot:\cdot,\cdot]$ which depend on the Cartesian coordinates of the vectors $\n(x)$, $x\in\Mc$.} See \rf{D1-DF} and \rf{D2-DF}.
\par Theorem \reff{CF-CR} below provides another explicit criterion for Lipschitz selections of a set-valued mapping $F:\Mc\to\HPL$. This criterion is formulated in terms of geometric objects that depend only on $F$ and do not depend on the choice of the coordinate system in $\RT$. We refer to this criterion as a ``coordinate-free'' Lipschitz selection criterion.
\smsk
\par Let us introduce additional definitions and notation necessary for its formulation. Given $x,y\in\Mc$, we let $\AF(x,y)\in[0,\pi/2]$ denote the angle between the boundaries of $F(x)$ and $F(y)$, i.e., between the straight lines $\ell_F(x)$ and $\ell_F(y)$. Clearly, $\sin{\AF(x,y)}=|\DXY|$, see \rf{DET-XY}.
\par Given $\Mc'\subset\Mc$, by $\diam_\rho(\Mc')$ we denote the diameter of $\Mc'$ in the pseudometric space $\MR$. We set $0/0=0$, $a/0=+\infty$ for every $a>0$, and $\dist(\hspace{0.2mm}\emp,A)=0$ for $A\subset\RT$.
Finally, we put
$$
\Df[x,x';y,y']=
\frac{\rho(x,x')}{\SXP}+\frac{\rho(y,y')}{\SYP}+
\diam_\rho\{x,x',y,y'\}~~~\text{provided}~~~
x,x',y,y'\in\Mc.
$$
\begin{theorem}\lbl{CF-CR} Let $\Mf=\MR$ be a pseudometric space, and let $F:\Mc\to\HPL$ be a set-valued mapping. Assume that either $\Mc$ is finite or the half-planes $\{F(x):x\in\Mc\}$ are in general position.
\par The mapping $F$ has a Lipschitz selection $f:\Mc\to\RT$ if and only if there exists a constant
$\lambda\ge 0$ such that for every four elements $x,x',y,y'\in\Mc$ the following inequality
\bel{MC-2}
\dist(F(x)\cap F(x'),F(y)\cap F(y'))\le
\lambda\,\Df[x,x';y,y']
\ee
holds. Furthermore,
$\tfrac{1}{\sqrt{2}}\inf\lambda
\le |F|_{\Mf}\le \gamma\inf\lambda$
where $\gamma>0$ is an absolute constant. See Fig. 3.
\end{theorem}

\begin{figure}[h!]
\hspace{10mm}
\includegraphics[scale=0.9]{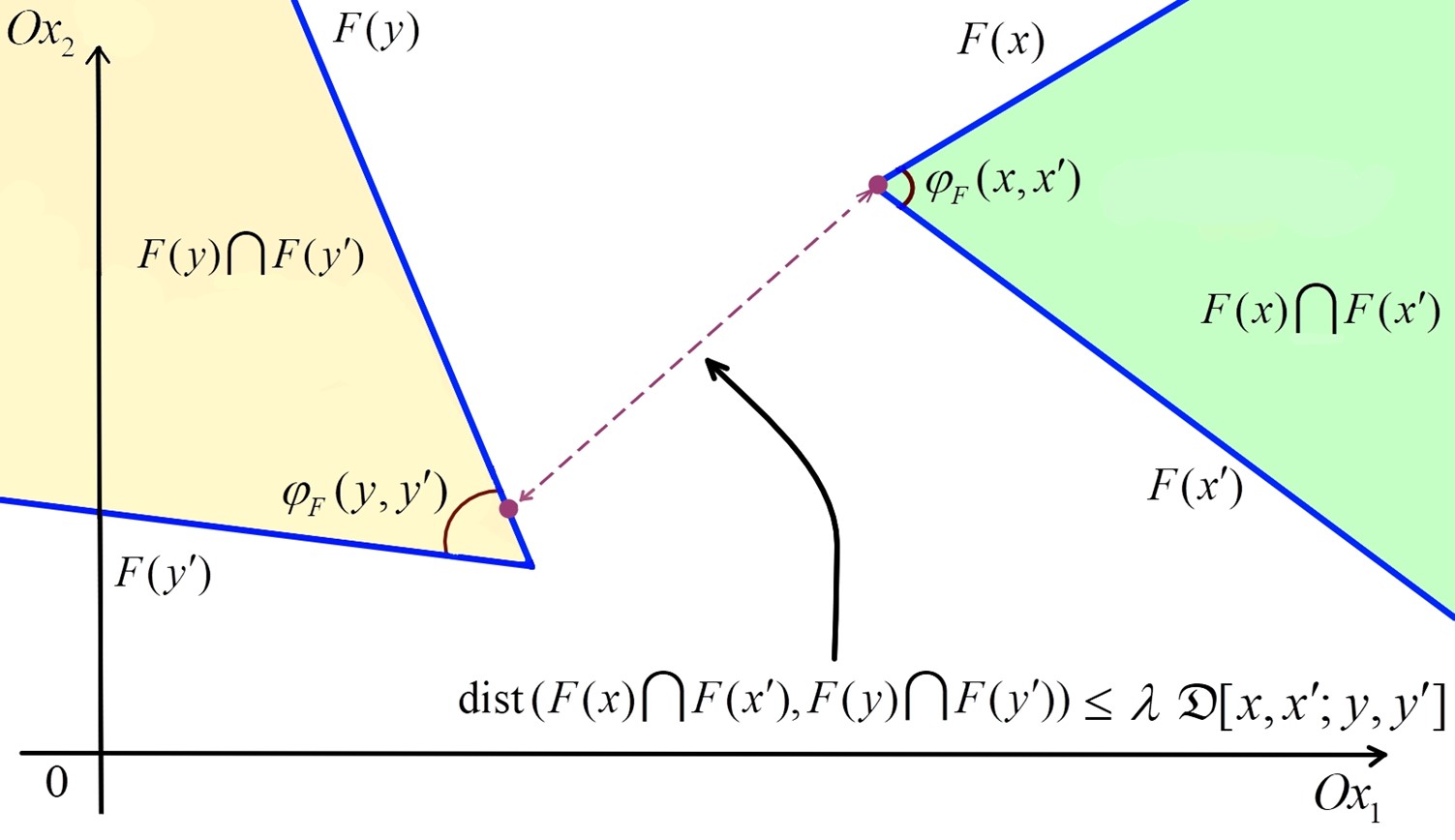}
\caption{The coordinate-free criterion for Lipschitz selections.}
\end{figure}
\par Theorem \reff{CR-1L} and Theorem \reff{CF-CR} lead us to two different criteria for Lipschitz selections which provide two different solutions to Problem \reff{PR1} for the family $\Tc=\KRT$ of all convex compact subsets of $\RT$.
\par Let us fix some notation that we need to formulate these results. Given $\n\in\SO$ and $\alpha\in\R$, we let $H(\n,\alpha)$ denote a closed half-plane defined by 
$$
H(\n,\alpha)=\{a\in\RT: \ip{\n,a}+\alpha\le 0\}.
$$
\par Let $\Mf=\MR$ be a pseudometric space, and let $G:\Mc\to\KRT$ be a set-valued mapping. For each  $x\in\Mc$, let us fix a family of half-planes $\HYG(x)\subset \HPL$ such that
\bel{H-NA}
G(x)=\cap\left\{H: H\in\HYG(x)\right\}.
\ee
\par Of course, the family $\HYG(x)$ can be defined in various ways: for instance, thanks to separation theorem, one can set $\HYG(x)=\{H\in\HPL: H\supset G(x)\}$.
\par A smaller family $\HYG(x)$ satisfying \rf{H-NA} one can define as follows: to each $x\in\Mc$ and each $\n\in\SO$ we assign a half-plane $V(x;\n)$ defined by
$$
V(x;\n)=\{a\in\RT:\ip{\n,a}\le h(\n:G(x))\}.
$$
Here
$h(\n:A)=\sup\{\ip{\n,a}:a\in A\}$
is {\it the support function} of a convex set $A\subset\RT$.
\par Then we set $$\HYG(x)=\{V(x;\n):\n\in\SO\},~~~~x\in\Mc.$$ Thus, in this case, the family $\HYG(x)$ is the family of all {\it supporting half-planes} of the set $G(x)$.
\begin{theorem}\lbl{FOR-1} Let $\Mf=\MR$ be a pseudometric space, and let $G:\Mc\to\KRT$ be a set-valued mapping.
This mapping has a Lipschitz selection if and only if
there exists a constant $\lambda\ge 0$ such that the following two conditions are satisfied:
\smsk
\par (i) $\dist(G(x),G(y))\le \lambda\,\rho(x,y)$ for every $x,y\in\Mc$;
\smsk
\par (ii) Let $\Mc'=\{x,x',y,y'\}\subset\Mc$ be an arbitrary set, and let $\n:\Mc'\to \SO$ and $\alpha:\Mc'\to\R$ be arbitrary mappings such that $$H(\n(z),\alpha(z))\in\HYG(z)~~~\text{for every}~~~z\in\Mc'.$$
Let $F(z)= H(\n(z),\alpha(z))$ on $\Mc'$.
\par Then condition ($\bigstar 2$) of Theorem \reff{CR-1L} holds for $x,x',y,y'$ and $F$. (Recall that this condition is equivalent to condition ($\bigstar 2'$), see \rf{CR1-PR}.)
\smsk
\par Furthermore,
$$
\tfrac{1}{\sqrt{2}}\,\inf\lambda\le \,|G|_{\Mf} \le 5\inf\lambda.
$$
See Fig. 4.
\end{theorem}
\begin{figure}[H]
\hspace{17mm}
\includegraphics[scale=0.33]{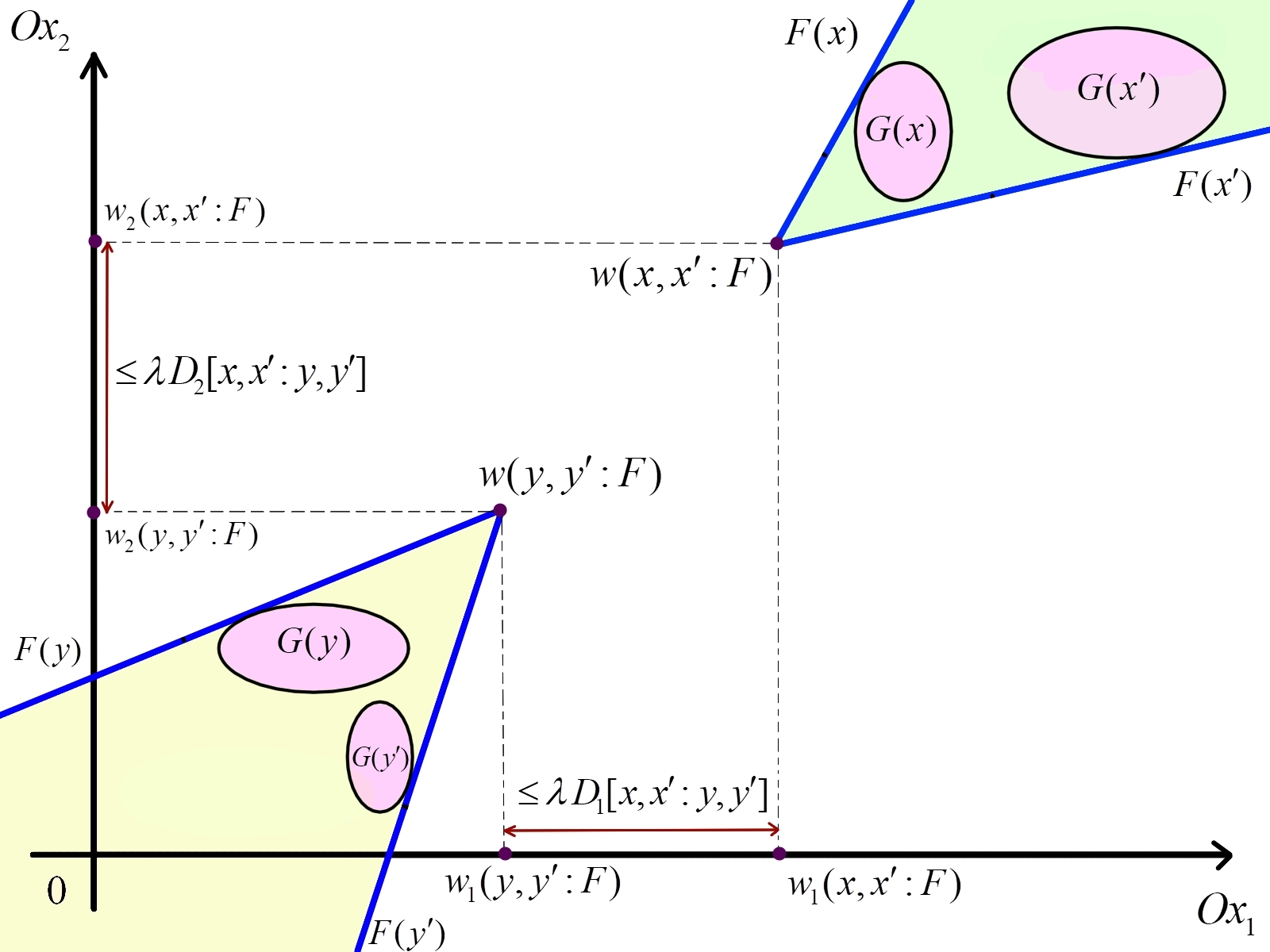}
\caption{The first Lipschitz selection criterion for the family $\KRT$.}
\end{figure}
\par Let us formulate the second criterion for Lipschitz selections.
\begin{theorem}\lbl{FOR-2} Let $\Mf=\MR$ be a pseudometric space, and let $G:\Mc\to\KRT$ be a set-valued mapping. The mapping $G$ has a Lipschitz selection if and only if there exists a constant $\lambda\ge 0$ such that for every four elements $x,x',y,y'\in\Mc$ and every four half-planes
$$
F(x)\in\HY(x),~~F(x')\in\HY(x'),
~~F(y)\in\HY(y),
~~F(y')\in\HY(y'),
$$
inequality \rf{MC-2}  holds.
\par Furthermore,
$$
\tfrac{1}{\sqrt{2}}\,\inf\lambda
\le \,|G|_{\Mf}\le \gamma\,\inf\lambda
$$
where $\gamma>0$ is an absolute constant.
\smsk
\par See Fig. 5.
\end{theorem}
\msk

\begin{figure}[h!]
\hspace{10mm}
\includegraphics[scale=0.9]{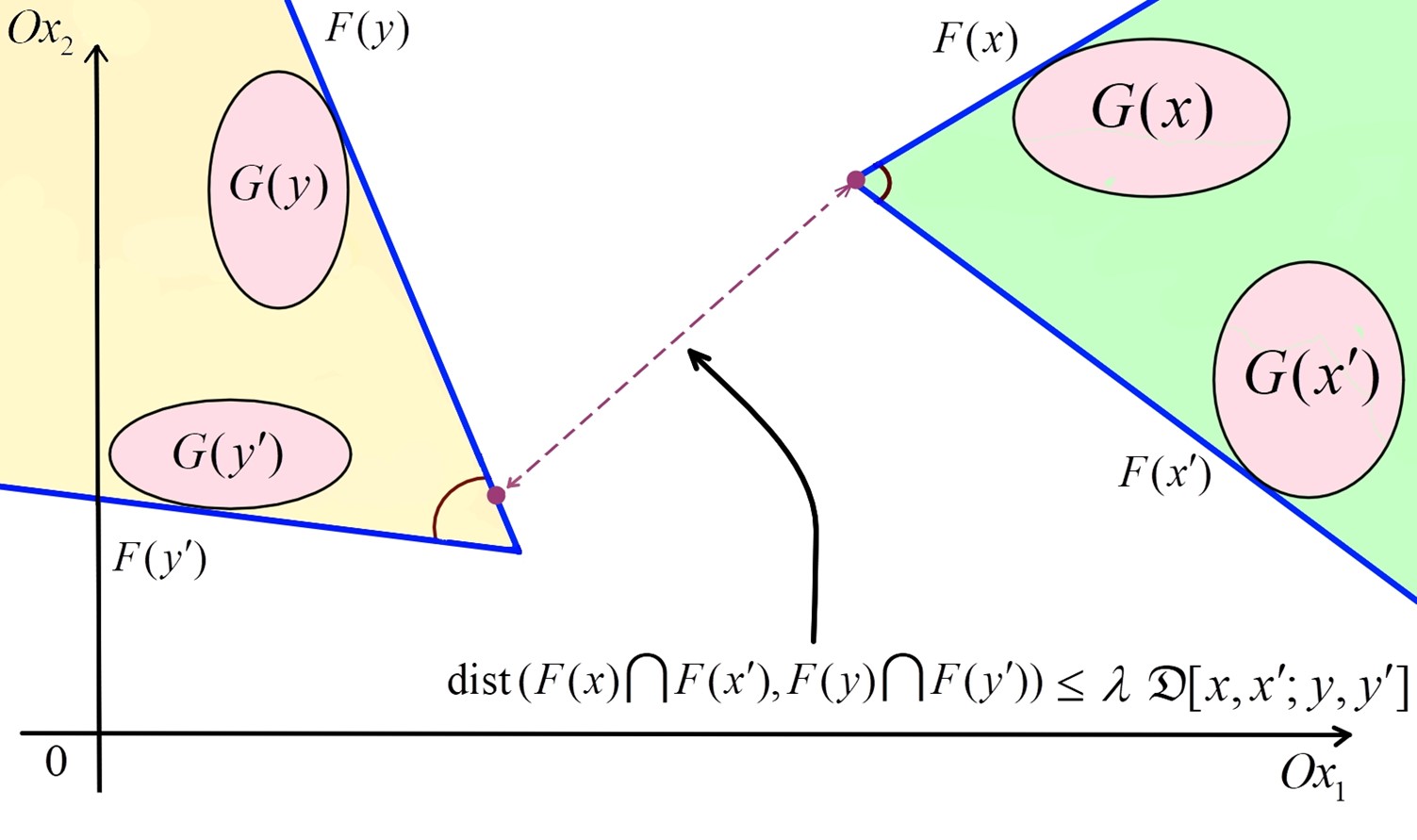}
\caption{The second Lipschitz selection criterion for the family $\KRT$.}
\end{figure}
\bsk
\par Let us now briefly describe the structure of the present paper and the main ideas of the proofs of the results stated above.
\par Given a convex set $S\subset\RT$, we let $\HR[S]$ denote  {\it the smallest rectangle (possibly unbounded) with sides parallel to the coordinate axes, containing} $S$. Thus, 
$\HR[S]=\Prj_1[S]\times \Prj_2[S]$.
(Recall that  $\Prj_i$ is the operator of the orthogonal projection onto the axis $Ox_i$, $i=1,2$.)
\par We refer to $\HR[S]$ as {\it the rectangular hull} of $S$. (See Section 2.2 for more details.)

\smsk
\par Let $Q_0=[-1,1]\times[-1,1]$. Given a set-valued mapping $F:\Mc\to\CRT$, a constant $\lambda\ge 0$ and elements $x,x'\in\Mc$, we let $\RL[x,x':\lambda]$ denote the rectangle defined by
\bel{RL}
\RL[x,x':\lambda]=
\HR[F(x)\cap\{F(x')+\lambda\,\rho(x,x') Q_0\}].
\ee
See Fig. 6.
\begin{figure}[h!]
\hspace{32mm}
\includegraphics[scale=0.7]{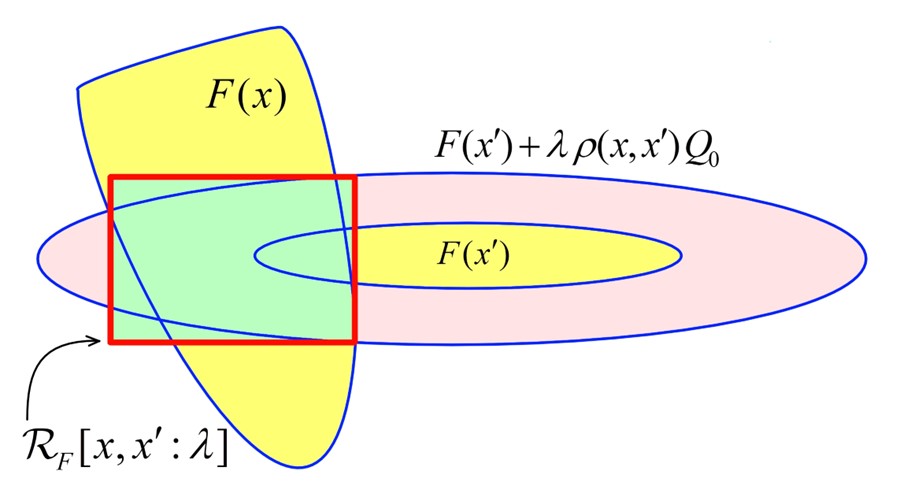}
\caption{The rectangle $\RL[x,x':\lambda]$ for $x,x'\in\Mc$ and $\lambda\ge 0$.}
\end{figure}
\par In this paper we prove a number of Lipschitz selection criteria for set-valued mappings $F:\Mc\to\CRT$ which work under certain natural conditions on the pseudometric space $\Mf=\MR$ and $F$. Here are these conditions.
\begin{condition}\lbl{CND-T} Either $\Mc$ is finite or
there exist a constant $\alpha\ge 0$ and elements $\brx_1,...,\brx_m\in\Mc$ such that the set
$\cap\{F(\brx_k)+\alpha Q_0:k=1,...,m\}$ is non-empty and bounded.
\end{condition}
\par The following result is one of the main ingredients of our approach.
\begin{theorem}\lbl{CR-LS1} Let $\Mf=\MR$ be a pseudometric space, and let $F:\Mc\to\CRT$ be a set-valued mapping. Suppose that $\Mf$ and $F$ satisfy Condition \reff{CND-T}.
\par The mapping $F$ has a Lipschitz selection if and only if there exists a constant $\lambda\ge 0$ such that
\bel{R-MLM}
\RL[x,x':\lambda]\cap \{\RL[y,y':\lambda]+\lambda\,\rho(x,y)\,Q_0\}\ne\emp
~~~~~\text{for every}~~x,x',y,y'\in\Mc.
\ee
\par Furthermore, the following inequalities hold:
\bel{OP-NF}
\inf\lambda\le \,|F|_{\Mf}\le 5\inf\lambda.
\ee
\end{theorem}
\par Clearly, \rf{R-MLM} is equivalent to the existence of points $u\in \RL[x,x':\lambda]$ and $v\in\RL[y,y':\lambda]$ (depending on $F$, $\rho$ and $x,x',y,y'$) such that $\|u-v\|\le\lambda\,\rho(x,y)$. See Fig. 7.
\begin{figure}[h!]
\hspace{18mm}
\includegraphics[scale=0.84]{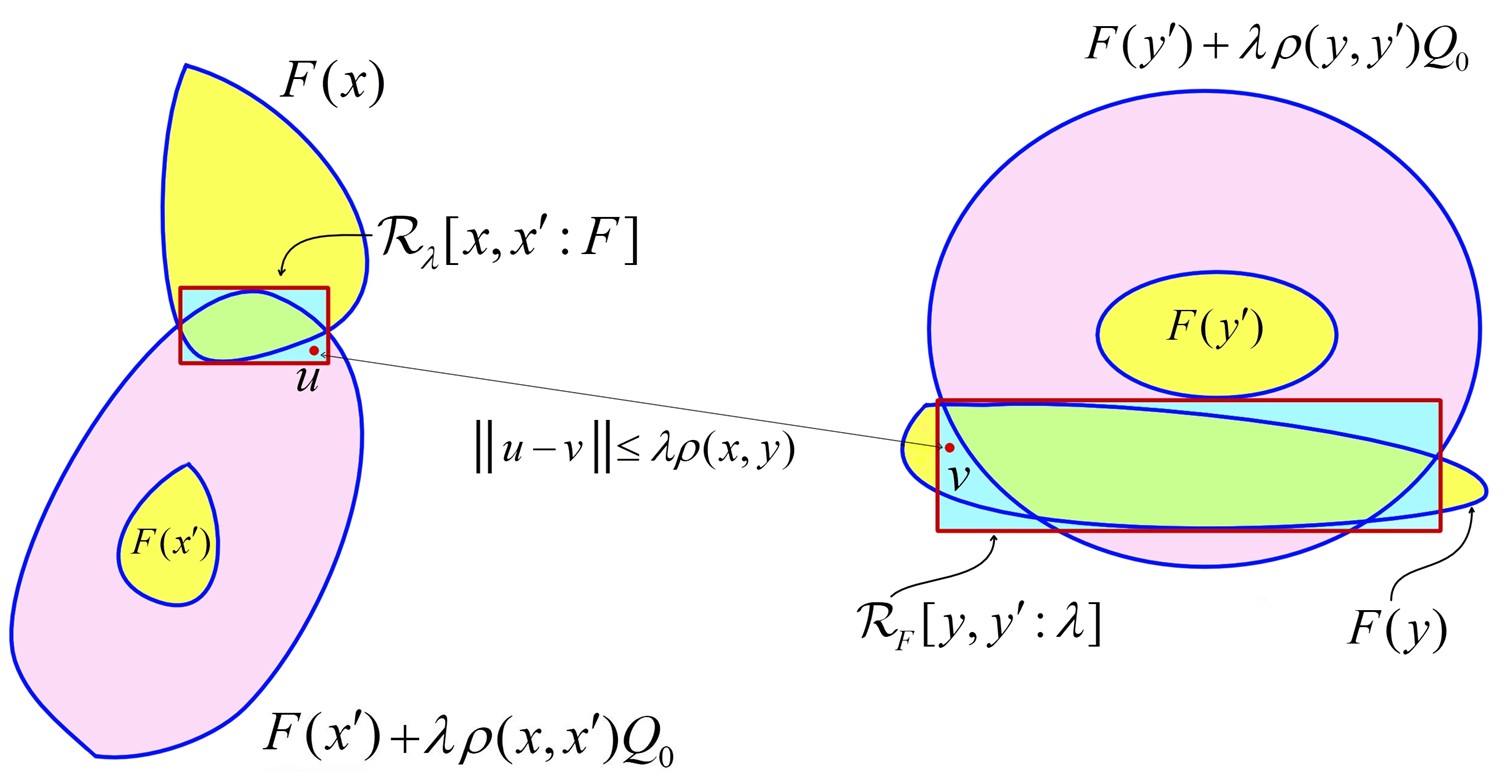}
\caption{The Lipschitz selection criterion in $\RT$.}
\end{figure}
\par We give a proof of Theorem \reff{CR-LS1} in Section 5. This proof relies on the criterion for Lipschitz selections given in Theorem \reff{W-CR}. In turn, the proof of Theorem \reff{W-CR} is based on the further development and generalization of the ideas and methods of work \cite{S-2002} devoted to  set-valued mappings taking values in the family $\KRT$, see \rf{KRT-DF}. Theorem \reff{W-CR} is the most technically difficult part of the present paper. We prove this theorem in Sections 3 and 4.

\begin{remark}\lbl{RMK-A} {\em (i) Let us note that, given  $\lambda\ge 0$ and $x,x',y,y'\in\Mc$, the intersection
$$
\RL[x,x':\lambda]\cap \{\RL[y,y':\lambda]+\lambda\,\rho(x,y)\,Q_0\}\ne\emp
$$
if and only if there exist points
$$
A(i)=(a_1(i),a_2(i))\in F(x), A'(i)\in F(x'), B(i)=(b_1(i),b_2(i))\in F(y)~~~\text{and}~~~
B'(i)\in F(y')
$$
such that
$$
\|A(i)-A'(i)\|\le \lambda\,\rho(x,x'),~~~
\|B(i)-B'(i)\|\le \lambda\,\rho(y,y')~~\text{and}~~|a_i(i)-b_i(i)|\le \lambda\,\rho(x,y).
$$
\par See Fig. 8.
\begin{figure}[h!]
\hspace{8mm}
\includegraphics[scale=0.7]{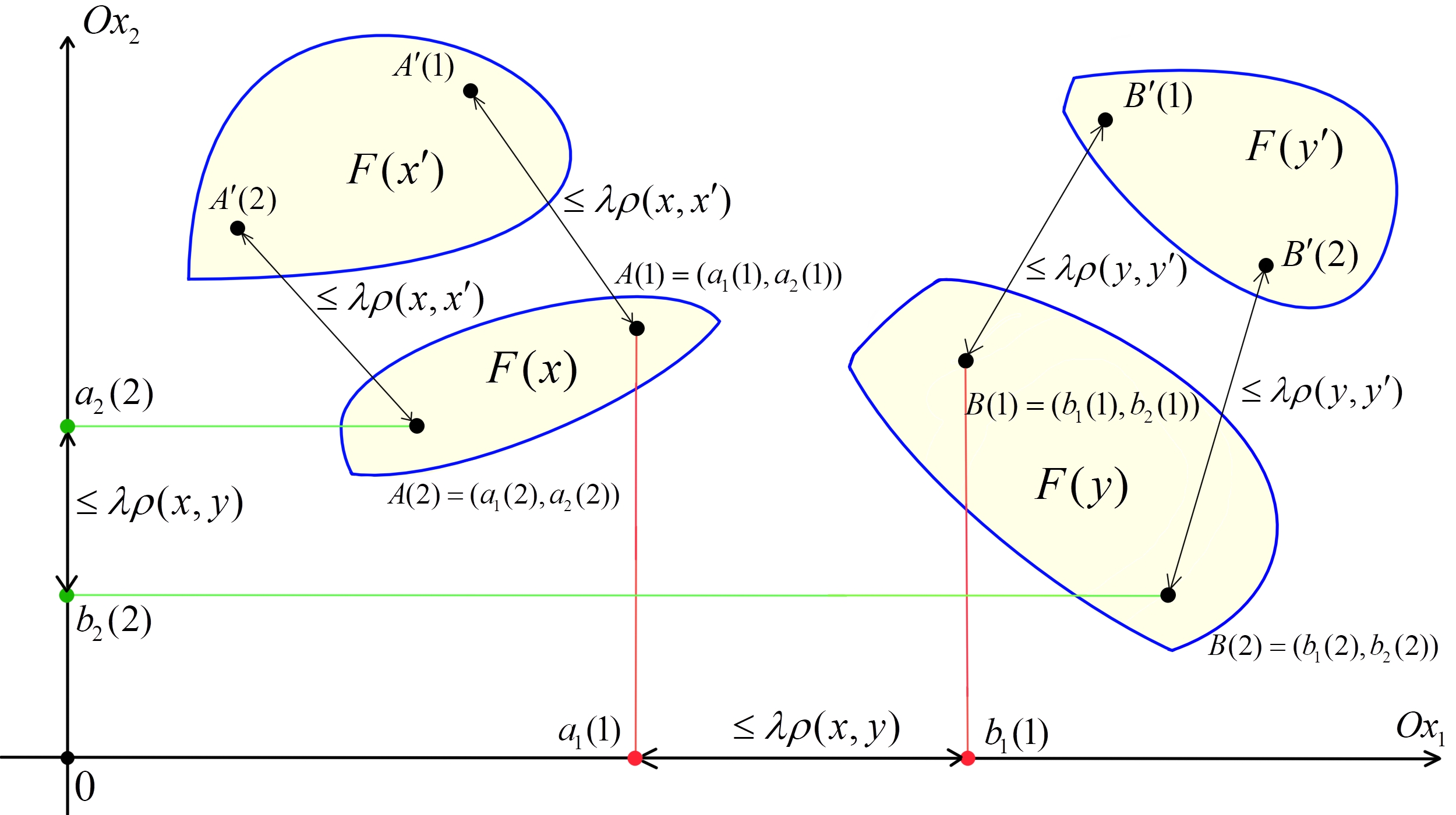}
\caption{The Lipschitz selection criterion in $\RT$.}
\end{figure}
\smsk
\par (ii) Suppose that for every $x,x'\in\Mc$, the rectangle $\RL[x,x':\lambda]$ is a {\it closed set}. (In particular, this property holds provided $F:\Mc\to\HPL$, or $\Mc$ and $F$ satisfy Condition \reff{CND-T}. However, in general, this property of $F$ does not hold. See Remarks \reff{RM-RCB} and \reff{RM-TAU} below.)
\par In this case, \rf{R-MLM} is equivalent to the following two conditions: for every $x,y\in\Mc$,
$$
\dist(F(x),F(y))\le\lambda\,\rho(x,y),
$$
and
$$
\dist\left(\,\RL[x,x':\lambda],\RL[y,y':\lambda]\,\right)
\le\lambda\,\rho(x,y)~~~~\text{for every}~~x,x',y,y'\in\Mc.~\text{\rbx}
$$
}
\end{remark}
\par The criterion of Theorem \reff{CR-LS1} provides important geometrical tool for tackling Problem \reff{PR1}. It tells us that the quantity $\FM$ (i.e., the Lipschitz seminorm of the optimal Lipschitz selection of $F$, see \rf{FM}) is equivalent to the smallest value of the constant $\lambda\ge 0$ such that condition \rf{R-MLM} is satisfied {\it for any choice} of elements $x,y,x',y'\in\Mc$.
\smsk
\par We recall that Problem \reff{PR1} is the problem of finding an {\it explicit} formula for $\FM$. Theorem \reff{CR-LS1} takes an important step in this direction by reduction the original problem to the purely geometrical problem of calculation the optimal $\lambda$ satisfying condition \rf{R-MLM}.
\par We note that this is also a rather difficult optimization problem, and, in general, given elements $x,y,x',y'\in\Mc$ we are unable to give an explicit formula for the order of magnitude of the smallest value of $\lambda$ for which \rf{R-MLM} is satisfied. (Even if each $F(x)$ is a {\it disk} on the plane with center $c(x)$ and radius $r(x)$, we do not know how to express $\inf \lambda$ in terms of the functions $c(\cdot)$ and $r(\cdot)$.)
\smsk
\par Fortunately, we are able to solve the above problem and express the order of magnitude of the optimal $\lambda$ from \rf{R-MLM} in explicit geometrical terms provided $F$ is a mapping from $\Mc$ into the family $\HPL$ of all closed half-planes. This leads us to conditions ($\bigstar 1$) and ($\bigstar 2$) of Theorem \reff{CR-1L} providing the proof of this result. 
\smsk
\par In turn, the proof of Theorem \reff{CF-CR} relies on the finiteness principle for Lipschitz selections \rf{FP-FM} and the Lipschitz selection criterion given in Theorem \reff{CR-1L}. Finally, we prove Theorems \reff{FOR-1} and \reff{FOR-2} by applying Theorems \reff{CR-1L} and \reff{CF-CR} to a new pseudometric space
$(\tMc,\trh)$ where
$$
\tMc=\{(x,H):x\in\Mc, H\in\HYG(x)\}
$$
and
$$
\trh((x,H),(x',H'))=\rho(x,x')~~~\text{for}~~~x,x'\in\Mc
~~~\text{and}~~~H\in\HYG(x),H'\in\HYG(x').
$$
\par We refer the reader to paper \cite{S-2020-L} for the detailed proofs of these theorems.
\smsk
\par In Section 5, we present three corollaries of Theorem \reff{W-CR}. These are Theorems \reff{FP-RT}, \reff{CR-W5} and \reff{CR-LS-RI}. Each of these results provides a criterion for Lipschitz selections in $\RT$. In particular, Theorem \reff{FP-RT} is the Finiteness Principle itself, but with the constant $\gamma=3$ in inequality \rf{FP-FM} which the smallest value of this constant known so far.
\smsk
\par Theorem \reff{CR-W5} is a variant of Theorem
\reff{CR-LS1} in which the rectangles $\RL[\cdot,\cdot:\lambda]$ defined by \rf{RL} are replaced with the rectangles $\Wc_F[\cdot,\cdot,\cdot:\lambda]$ defined by \rf{WC-DF}. As a result, we obtain inequality \rf{OP-NF} with a smaller constant in the right hand side ($3$ instead of $5$); on the other hand, the rectangles $\Wc_F$ are more complicated objects each depending on three elements of the set $\Mc$, while each rectangle $\RL$ depends only on two elements of $\Mc$.
\smsk
\par In turn, Theorem \reff{CR-LS-RI} enables us to reformulate the Lipschitz selection criterion of Theorem \reff{CR-LS1} in terms of intersections of certain rectangles in $\RT$. In particular, the combination of this criterion with the linear-time algorithms for linear programming developed by N. Megiddo \cite{M-1983} and  M. E. Dyer \cite{Dy-1984} provides an efficient algorithm for computing the quantity $|F|_{\Mf}$ for set-valued mappings $F:\Mc\to\HPL$ defined on finite pseudometric space $\Mf=\MR$. See Section 6.
\smsk
\par In Section 6 and Section 7 we present a solution to Problem \reff{PR2}. More specifically, in these sections we introduce and study two different constructive algorithms for Lipschitz selections which solve this problem for families $\Tf$ of convex closed sets satisfying rather mild geometrical conditions. We refer to these algorithms as the {\it ``Projection Algorithm''} and the {\it ``Iterative Algorithm''} respectively. Let us note that the Projection Algorithm is based on the approach to the Lipschitz selection problem developed in the proof of Theorem \reff{W-CR}, while the Iterative Algorithm relies on the ideas of the proof of the ``Stabilization Principle'' introduced in our recent paper \cite{S-2022}.
\smsk
\begin{remark} {\em We would like to point out that in this paper we deal only with {\it theoretical aspects of Projection and Iterative Algorithms}. We give a detailed description and justification of these algorithms, as well as some preliminary estimates of the their computational efficiency.
\par But we do not cover all the computational aspects of their work. This issue will be discussed in the next paper \cite{S-2023}, where we will show how these algorithms can be implemented on an (idealized) computer with the standard von Neumann architecture.
\par We will prove that these algorithms are efficient, i.e., they require a minimal use of computer resources. In \cite{S-2023} we will give efficient estimates of the ``work'' of these algorithms (i.e., the number of machine operations needed to carry them out), as well as estimates of the amount of computer memory required to process them.\rbx}
\end{remark}

\smsk
\par {\bf Acknowledgements.} I am very grateful to Charles Fefferman for stimulating discussions and valuable advice.

\SECT{2. Notation and preliminaries.}{2}
\addtocontents{toc}{2. Notation and preliminaries.\hfill \thepage\par\VST}

\indent\par {\bf 2.1 Background notation.}
\addtocontents{toc}{~~~~2.1 Background notation.\hfill \thepage\par\VST}
\msk
\indent
\par Let $A$ and $B$ be non-empty subsets of $\RT$. We let
$$
A+B=\{a+b: a\in A, b\in B\}
$$
denote the Minkowski sum of these sets. Given $\lambda\ge 0$, by $\lambda A$ we denote the set $\lambda A =\{\lambda a: a\in A\}$.
\par We write
$$
\dist(A,B)=\inf\{\|a-b\|:a\in A,~b\in B\}
$$
to denote the distance between $A$ and $B$. For $x\in\RT$ we also set $\dist(x,A)=\dist(\{x\},A)$. We put $\dist(\emp,A)=0$ provided $A$ is an arbitrary (possibly empty) subset of $\RT$.
\par We let $A^{\cl}$ denote the {\it closure} of the set $A$.
\par For a Banach space $X$ with the unit ball $B_X$, and  two non-empty subsets $A$ and $B$ of $X$, we let $\dhf(A,B)$ denote the Hausdorff distance between $A$ and $B$ in $X$:
\bel{HD-DF}
\dhf(A,B)=\inf\{r>0: A+rB_X\supset B,~B+rB_X\supset A\}.
\ee
(Of course, $\dhf(A,B)$ also depends on $X$. However, we use $\dhf$ only in those places in the paper where $X$ is clear from the context. Therefore, we omit $X$ in the Hausdorff distance notation.)
\par Let $\RPL=\{a\in\R: a\ge 0\}$, and let $\RTP=\{a=(a_1,a_2)\in\RT: a_1,a_2\ge 0\}$. Given $a,b\in\RT$, $a\ne b$, by $[a,b]$ we denote the closed interval with the ends in $a$ and $b$:
$$
[a,b]=\{x\in\RT:x=(1-t)\,a+t\,b, 0\le t\le 1\}.
$$
\par Given a set $A\subset\R$ we put
$\min A=\{\min x:x\in A\}$ and $\max A=\{\max x:x\in A\}$ provided $A$ is a closed subset of $\R$ bounded from above or below respectively. We write $[x]_+$ for the positive part of the real $x$, i.e., $[x]_+=\max\{x,0\}$. We also use the natural convention that
\bel{INF-S}
\frac{0}{0}=0,~~\frac{a}{0}=+\infty~~\text{for}~~a>0,
~~a-b=0~~~\text{if}~~a=b=\pm\infty,~~\text{and}~~
(\pm\infty)-(\mp\infty)=\pm\infty.
\ee
\smsk
\par If $S$ is a finite set, by $\#S$ we denote the number of elements of $S$.
\par By
$$
Ox_1=\{x=(t,0):t\in\R\}~~~\text{and}~~~  Ox_2=\{x=(0,t):t\in\R\}
$$
we denote the coordinate axes in $\RT$. We recall that by $\Prj_i$, $i=1,2$, we denote the operator of the orthogonal projection onto the axis $Ox_i$.
Thus, given $x=(x_1,x_2)\in\RT$, we have $\Prj_1[x]=(x_1,0)$ and $\Prj_2[x]=(0,x_2)$. Given sets $A_i\subset Ox_i$, $i=1,2$, we let $A_1\times A_2$ denote a subset of $\RT$ defined by
\bel{A1TA2}
A_1\times A_2=\{a=(a_1,a_2)\in\RT: (a_1,0)\in A_1, (0,a_2)\in A_2\}.
\ee
In other words,
\bel{A12-PR}
A=A_1\times A_2~~~~\text{if and only if}~~~~
\Prj_1(A)=A_1~~~\text{and}~~~\Prj_2(A)=A_2.
\ee
\par Given $a\in\RT$ and $r>0$, we let $Q(a,r)$ denote the square with center $a$ and length of side $2r$:
\bel{SQ-DF}
Q(a,r)=\{y\in\RT:\|y-a\|\le r\}.
\ee
In particular, $Q_0=[-1,1]^2=Q(0,1)$ is the unit ball of the Banach space $\LTI=(\RT,\|\cdot\|)$.
\smsk
\par Let $S$ be a non-empty {\it convex} closed subset of $\RT$.  By $\Prm(\cdot,S)$ we denote the operator of metric projection onto $S$ in $\LTI$-norm. To each $a\in\RT$ this operator assigns the set of all points in $S$ that are nearest to $a$ on $S$ in the uniform norm. Thus,
\bel{MPR}
\Prm(a,S)=S\cap Q(a,\dist(a,S)).
\ee
Clearly, the set $\Prm(a,S)$ is either a singleton or a line segment in $\RT$ parallel to one of the coordinate axes.
\par If $S\subset\RT$ is convex bounded and centrally symmetric, by $\cent(S)$ we denote the center of $S$.
\smsk
\par We let $\LTT=(\RT,\|\cdot\|_{\LTT})$ denote $\RT$ equipped with the standard Euclidean norm $\|a\|_{\LTT}=(a_1^2+a_2^2)^{\frac{1}{2}}$, $a=(a_1,a_2)$. Let
$$
B_0=\{a=(a_1,a_2)\in\RT: a_1^2+a_2^2\le 1\}
 ~~~~\text{and}~~~~
\SO=\{a=(a_1,a_2)\in\RT:
a_1^2+a_2^2=1\}
$$
be the closed unit disk and the unit circle in $\RT$ respectively. Given non-zero vectors $u,v\in\RT$ we write $u\parallel v$ if $u$ and $v$ are collinear, and we write $u\nparallel v$ whenever these vectors are non-collinear.
We say that the vectors $u,v\in\RT$ are {\it co-directed} if $u,v\ne 0$, $u$ and $v$ are collinear and have the same direction, i.e., $v=\alpha u$ for some $\alpha>0$.
\par By $\theta(u,v)\in[0,2\pi)$ we denote
$$
\text{the angle of rotation from}~~~u/\|u\|_{\LTT}~~~\text{to}~~~  v/\|v\|_{\LTT}~~~\text{in the counterclockwise direction}.
$$
(Thus, $\theta(v,u)=2\pi-\theta(u,v)$.) We refer to $\theta(u,v)$ as the angle between the vectors $u$ and $v$.
\smsk
\par Let $\ell_1$ and $\ell_2$ be two non-parallel straight lines in $\RT$; in this case, we write $\ell_1\nparallel\ell_2$. Let $V=\ell_1\cap \ell_2$. These two lines form two angles $\vf_1,\vf_2\in[0,\pi)$, $\vf_1+\vf_2=\pi$, with the vertex at the point $V$. Let
$$
\vf(\ell_1,\ell_2)=\min\{\vf_1,\vf_2\};~~~~~
\text{clearly,}~~~~ \vf(\ell_1,\ell_2)\in[0,\pi/2].
$$
\par We refer to $\vf(\ell_1,\ell_2)$ as {\it ``the angle between the straight lines $\ell_1$ and $\ell_2$''.} If $\ell_1\parallel \ell_2$ (i.e., $\ell_1$ and $\ell_2$ are parallel), we set $\vf(\ell_1,\ell_2)=0$.

\indent\par {\bf 2.2 Rectangles and rectangular hulls.}
\addtocontents{toc}{~~~~2.2 Rectangles and rectangular hulls.\hfill \thepage\par\VST}
\msk
\indent
\par Let $\Ic(Ox_i)$, $i=1,2$, be the family of all non-empty convex subsets of the coordinate axis $Ox_i$. In other words, $\Ic(Ox_i)$ is the family of all non-empty intervals (bounded or unbounded) lying on the $Ox_i$ axis. We set
\bel{RCT-D}
\RCT=\{\Pi=I_1\times I_2: I_1\in\Ic(Ox_1),I_2\in\Ic(Ox_2)\}.
\ee
We refer to every member of the family $\RCT$ as a {\it ``rectangle''}. Furthermore, throughout the paper, the word ``rectangle'' will mean an element of $\RCT$, i.e., a rectangle (possibly unbounded and not necessarily closed) with ``sides'' parallel to the coordinate axes.
\par Clearly, thanks to definition \rf{RCT-D},
\bel{PR-RW}
\Pi=\Prj_1[\Pi]\times \Prj_2[\Pi]~~~~\text{for every rectangle}~~~\Pi\in\RCT.
\ee
\par Because $\Prj_i$, $i=1,2$, is a continuous operator, for every {\it convex} set $S\subset\RT$ its orthogonal projection $\Prj_i[S]$ onto the axis $Ox_i$ is a {\it convex} subset of $Ox_i$, i.e., $\Prj_i[S]\in\Ic(Ox_i)$. Thus, thanks to this property, \rf{A1TA2} and \rf{A12-PR}, we have the following:
\bel{RC-PD}
\text{a convex set}~~\Pi\in\RCT~~\text{if and only if}~~   \Pi=\Prj_1[\Pi]\times \Prj_2[\Pi].
\ee
\par Let us also note the following intersection property of rectangles: non-empty intersection of any collection $\Rc$ of rectangles is a rectangle as well. Furthermore,
in this case,
\bel{RC-IP}
\bigcap_{\Pi\in\Rc} \Pi=\left\{\bigcap_{\Pi\in\Rc}\Prj_1[\Pi]\right\}
\times \left\{\bigcap_{\Pi\in\Rc}\Prj_2[\Pi]\right\}.
\ee
\par We also note that any interval in $\RT$ lying on a line parallel to a coordinate axis is a ``rectangle''. In particular, every interval on the axis $Ox_1$ or $Ox_2$ belongs to the family $\RCT$.
\par Finally, given a bounded rectangle $\Pi\in\RCT$ (which is not necessarily centrally symmetric set) we set $\cent(\Pi)=\cent(\Pi^{\cl})$.
\par Given a set $S\in\CRT$, see \rf{CRT-DF}, we let $\HR[S]$ denote the {\it ``rectangular hull``} of $S$, i.e., the smallest (with respect to inclusion) rectangle containing $S$. Thus,
\bel{HRS}
\HR[S]=\cap\{\Pi: \Pi\in\RCT, \Pi\supset S\}.
\ee
\par Combining this definition with \rf{RC-IP}, we conclude that
\bel{P1XP2}
\HR[S]=\Prj_1[S]\times\Prj_2[S].
\ee
\par Thus, given $S\in\CRT$, its rectangular hull $\HR[S]$ is the only rectangle $\Pi$ for which
\bel{HS-U}
\Prj_1[\Pi]=\Prj_1[S]~~~~\text{and}~~~~~
\Prj_2[\Pi]=\Prj_2[S].
\ee
\par Let us also note the following elementary property of rectangles: for every closed convex set $S\subset\RT$ and every $r\ge 0$ we have
\bel{N-HS}
\HR[S+r Q_0]=\HR[S]+r Q_0.
\ee
\begin{remark}\lbl{RM-RCB} {\em If $S$ is a convex {\it compact} subset of $\RT$ then its orthogonal projections onto the coordinate axes are closed bounded intervals, so that, thanks to \rf{P1XP2}, $\HR[S]=\Prj_1[S]\times\Prj_2[S]$ is {\it compact} as well.
\par Clearly, if the set $S\in\CRT$ is unbounded, then its rectangular hull $\HR[S]$ is also unbounded. But {\it we cannot guarantee that in this case $\HR[S]$ is closed} (in spite of $S$ itself is closed). More specifically, $\HR[S]$ is not closed if $\partial S$, the boundary of $S$, has {\it either a horizontal or vertical asymptote} (i.e., a straight line parallel to one of the coordinate axis having the $0$-distance to $S$ and empty intersection with $S$). For instance, let $S=\{x=(x_1,x_2)\in\RT:x_1,x_2>0,\, x_1\cdot x_2\ge 1\}$ be the epigraph of the function $f(t)=1/t$, $t>0$. Then $\HR[S]=\{x=(x_1,x_2)\in\RT:x_1,x_2>0$\} is the positive (open) octant.\rbx}
\end{remark}

\par In this section we present two important auxiliary results. The first of them is a variant of the classical Helly's intersection theorem for rectangles.
\begin{lemma}\lbl{H-R} Let $\Kc\subset\RCT$ be a collection of rectangles in $\RT$. Suppose that either (i) $\Kc$ is finite or (ii) every rectangle from $\Kc$ is closed, and there exists a finite subfamily of $\Kc$ with a non-empty bounded intersection.
\par Under these conditions, the following is true: If the intersection of every two rectangles from $\Kc$ is non-empty, then there exists a point in $\RT$ common to all of the family $\Kc$.
\end{lemma}
\par {\it Proof.} Representation \rf{PR-RW} reduces the problem to the one dimensional case. In this case the statement of the lemma is a variant of Helly's theorem in $\R$. See, e.g. \cite{DGK-1963}.\bx
\smsk
\par The second auxiliary result is a Helly-type theorem formulated in terms of the orthogonal projections onto the coordinate axes.
\begin{proposition}\lbl{INT-RE} Let $\Cf$ be a family of convex closed subsets in $\RT$. Suppose that either (i) $\Cf$ is finite or (ii) there exists a finite subfamily $\WTC\subset\Cf$ such that the intersection $\cap\{C:C\in\WTC\}$ is non-empty and bounded. If
\bel{P1-C}
\Prj_1[C_1\cap C_1']\,\cbg \Prj_1[C_2\cap C_2'] \ne\emp.
\ee
for every $C_1,C_1',C_2,C_2'\in\Cf$,\, then
\bel{IP-C}
\cap\{C:C\in\Cf\}\ne\emp.
\ee
\par Furthermore, in this case
\bel{PH-C}
\HR\left[\cap\{C:C\in\Cf\}\right]=
\cap\{\HR[C\cap C']: C,C'\in\Cf\}\,.
\ee
\end{proposition}
\par {\it Proof.}  Condition \rf{P1-C} tells us that for every $C,C'\in\Cf$ the set $C\cap C'$ is a non-empty. Because $C\cap C'$ is a convex subset of $\RT$, its projection onto $Ox_1$, the set $\Prj_1[C\cap C']\subset Ox_1$, is convex as well, i.e., this set is an interval in $Ox_1$.
\par First, let us prove the proposition provided condition (i) of the proposition's hypothesis holds, i.e., $\Cf$ is a {\it finite} family. Let $\Wc=\{\Prj_1[C\cap C']:C,C'\in\Cf\}$. Then $\Wc$ is a finite family of intervals, and, thanks to \rf{P1-C}, every two members of this family have a common point. Helly's theorem tells us that in this case there exists a point in $Ox_1$ common to all of the family $\Wc$. See Lemma \reff{H-R}. Thus,
\bel{V-L}
V=\bigcap_{C,C'\in\,\Cf}\Prj_1[C\cap C']\ne\emp\,.
\ee
\par Fix a point $v\in V$. Then, thanks to \rf{V-L},
\bel{SV-LB}
v\in \Prj_1[C\cap C']~~~~\text{for every}~~~~C,C'\in\Cf.
\ee
\par Let
\bel{L-DFB}
L=\{w\in\RT:\Prj_1[w]=v\}
\ee
be the straight line through $v$ orthogonal to the axis $Ox_1$.
\par Given $C\in\Cf$, we set $K(C)=C\cap L$. Thanks to \rf{SV-LB}, $v\in\Prj_1[C\cap C]=\Prj_1[C]$ so that there exists $u_C\in C$ such that $\Prj_1[u_C]=v$. From this and \rf{L-DFB}, we have $u_C\in C\cap L$ proving that $K(C)\ne\emp$ for every $C\in\Cf$. Clearly, each $K(C)$ is a closed interval lying on the straight line $L$. Let us show that there exists a point in $L$ common to all these intervals.
\smsk
\par Property \rf{SV-LB} tells us that for every $C,C'\in\Cf$ there exists a point $\tu\in C\cap C'$ such that $\Prj_1[\tu]=v$. Hence, thanks to \rf{L-DFB},
$\tu\in L$ so that
$$
\tu\in L\cap C\cap C'=(L\cap C)\cap (L\cap C')=K(C)\cap K(C').
$$
\par This proves that any two members of the family $\Kc=\{K(C): C\in\Cf\}$ have a common point. Furthermore,  $\Kc$ is a {\it finite} (because $\Cf$ is finite) family of intervals lying in $L$. Helly's theorem tells us that in this case $\cap\{K(C):C\in\Cf\}\ne\emp$. Thus,
$$
\cap\{K(C):C\in\Cf\}=\cap\{L\cap C:C\in\Cf\}
=L\cap\left(\cap\{C:C\in\Cf\}\right)\ne\emp
$$
proving \rf{IP-C} for a finite family $\Cf$.
\smsk
\par Let us see that \rf{IP-C} holds provided condition (ii) of the proposition's hypothesis holds. As we have proved above, $\cap\{C:C\in\Cf'\}\ne\emp$ for every {\it finite} subfamily $\Cf'$ of the family $\Cf$. In particular, {\it every three members of $\Cf$ have a common point}. This property and condition (ii) of the hypothesis tell us that $\Cf$ satisfies the hypothesis of {\it the two dimensional Helly's theorem} \cite{DGK-1963}. Thanks to this theorem, \rf{IP-C} holds  completing the proof of this property.
\smsk
\par Let us prove \rf{PH-C}. Clearly, the left hand side of the equality \rf{PH-C} is contained in its right hand side.
Let us prove the converse statement.
\par Fix a point
\bel{U-FS}
u=(u_1,u_2)\in\cap \{\HR[C\cap C']: C,C'\in\Cf\}
\ee
and prove that $u\in\HR\left[\cap\{C:C\in\Cf\}\right]$.
Thanks to \rf{P1XP2}, property \rf{U-FS} is equivalent to the following one:
\bel{PCO-12}
u_1\in\Prj_1[\Tc]~~~\text{and}~~~u_2\in\Prj_2[\Tc]~~~~
\text{where}~~~~\Tc=\cap\{C:C\in\Cf\}.
\ee
\par Prove that $u_1\in\Prj_1[\Tc]$. We let $L_1$ denote the straight line in $\RT$ through the point $u=(u_1,u_2)$ orthogonal to the axis $Ox_1$. Thus,
\bel{L1-OR}
L_1=\{w\in\RT:\Prj_1[w]=u_1\}.
\ee
\par Let us show that $L_1\cap\Tc\ne\emp$. Indeed, thanks to \rf{U-FS}, $u\in\HR[C\cap C']$ provided  $C,C'\in\Cf$ so that, thanks to \rf{P1XP2},
$$
u_1\in\Prj_1[C\cap C']~~~~\text{for every}~~~~ C,C'\in\Cf.
$$
Combining this property with definition \rf{L1-OR}, we conclude that $L_1\cap C\cap C'\ne\emp$ for all $C,C'\in\Cf$. Thus, every two members of the family $\Kc_1=\{L_1\cap C:C\in\Cf\}$ have a common point.
\par If the family $\Cf$ is finite, i.e., condition (i) of the proposition's hypothesis holds, then, thanks to Helly's theorem,
\bel{TU-C}
\text{there exists a point}~~~
\tu\in \cap\{L_1\cap C:C\in\Cf\}=L_1\cap\{C:C\in\Cf\}=L_1\cap\Tc.
\ee
\par Thus, $\tu\in L_1$ so that, thanks \rf{L1-OR}, $u_1=\Prj[\tu]$. Furthermore, $\tu\in \Tc=\cap\{C:C\in\Cf\}$ so that $u_1\in \Prj_1[\Tc]$.
This proves the first statement in \rf{PCO-12} provided condition (i) holds.
\smsk
\par Let us prove this statement provided the family $\Cf$ satisfies condition (ii). As we have seen above, it suffices to show that in this case the statement \rf{TU-C} holds.
\par We recall that in case (ii), there exists a {\it finite} subfamily $\WTC\subset\Cf$ such that  $\cap\{C:C\in\WTC\}$ is non-empty and bounded. Thus, condition (i) holds for  $\WTC$ so that
$\cap\{L_1\cap C:C\in\WTC\}\ne\emp$. Furthermore,
$$
\cap\{L_1\cap C:C\in\WTC\}\subset \cap\{C:C\in\WTC\}.
$$
But $\cap\{C:C\in\WTC\}$ is bounded so that the set $\cap\{L_1\cap C:C\in\WTC\}$ is bounded as well. Thus, the family
$$
\widetilde{\Kc_1}=\{L_1\cap C:C\in\WTC\}\subset\Kc_1
$$ 
has a {\it non-empty and bounded} intersection. This proves that the family $\Kc_1$ of closed intervals lying on the straight line $L_1$ satisfies all conditions of Helly's theorem for rectangles, see Lemma \reff{H-R}. Thanks to this lemma,
$\cap\{L_1\cap C:C\in\Cf\}\ne\emp$ proving the statement
\rf{TU-C} in the case under consideration.
\par Thus, we have proved the first statement of  \rf{PCO-12}, i.e., the property $u_1\in\Prj_1[\Tc]$. In the same way we prove that $u_2\in\Prj_2[\Tc]$ completing the proof of the proposition.\bx

\indent\par {\bf 2.3 Rectangles: intersections, neighborhoods and selections.}
\addtocontents{toc}{~~~~2.3 Rectangles: intersections, neighborhoods and selections.\hfill \thepage\par\VST}
\msk
\indent

\par In this section we present several criteria and several constructive formulae for the optimal Lipschitz selections of set-valued mappings taking values in the
family $\RCT$ of all rectangles in $\RT$ with sides parallel to the coordinate axes. See \rf{RCT-D}.
\par Let $I_0=[-1,1]$. Given $a\in\R$ and $r\ge 0$, we set $rI_0=[-r,r]$. We also recall that given a bounded interval $I\in\Ic(\R)$, by $\cent(I)$ we denote the center of $I$.
\begin{lemma}\lbl{H-RB} Let $\Kc\subset\RCT$ be a family of rectangles in $\RT$ with non-empty intersection. Then for every $r\ge 0$ the following equality
\bel{H-SM}
\left(\,\bigcap_{K\in\,\Kc} K\right) +r Q_0
=\bigcap_{K\in\,\Kc}\,\left\{\,K+rQ_0\,\right\}
\ee
holds.
\end{lemma}
\par {\it Proof.} Obviously, the right hand side of \rf{H-SM} contains its left hand side. Let us prove that
\bel{H-SM-LR}
\left(\,\bigcap_{K\in\,\Kc} K\right) +r Q_0
\supset\bigcap_{K\in\,\Kc}\,\left\{\,K+rQ_0\,\right\}.
\ee
\par This inclusion is based on the following simple claim: Let $\Ic$ be a family of convex subsets of $\R$ (intervals) with non-empty intersection, and let $K=[a,b]\subset\R$, be a closed bounded interval. Suppose that $K\cap I\ne\emp$ for every $I\in\Ic$. Then there exists a point common to $K$ and all of the members of the family $\Ic$.
\par Indeed, let $c\in\cap\{I:I\in\Ic\}$. If $c\in K$ then the claim holds. Suppose that $c\notin K$, say $c>b$. Prove that $b\in I$ for every $I\in \Ic$. In fact, because $K\cap I\ne\emp$, there exists a point $p_I\in K\cap I$. Then $p_K\le b$ so that $b\in[p_K,c]\subset I$ proving the claim.
\par This claim implies the following one dimensional variant of inclusion \rf{H-SM-LR}: Let $\Ic$ be a family of intervals in $\R$ with non-empty intersection. Then
\bel{H-SM-R1}
\left(\,\bigcap_{I\in\,\Ic} I\right) +r I_0
\supset\bigcap_{I\in\,\Ic}\,\left\{\,I+rI_0\,\right\}
~~~~\text{where}~~~~I_0=[-1,1].
\ee
\par Indeed, if $u\in\cap\{I+rI_0:I\in\Ic\}$ then $[u-r,u+r]\cap I\ne\emp$ for every $I\in\Ic$. Therefore, thanks to the above claim, $[u-r,u+r]\cap(\cap\{I:I\in\Ic\})\ne\emp$ proving that $u$ belongs to the left hand side of \rf{H-SM-R1}.
\smsk
\par Now, let us prove \rf{H-SM-LR} using \rf{H-SM-R1} and properties \rf{RC-PD} and \rf{RC-IP} of rectangles. For every $i=1,2$, we have
$$
\Prj_i\left[\left(\,\bigcap_{K\in\,\Kc} K\right)
+r Q_0\right]=
\left(\,\bigcap_{K\in\,\Kc} \Prj_i[K]\right)
+r\Prj_i[Q_0]=U_i.
$$
Furthermore,
$$
\Prj_i\left[
\bigcap_{K\in\,\Kc}\,\left\{\,K+rQ_0\,\right\}
\right]=
\bigcap_{K\in\,\Kc}\left\{\Prj_i[K]+r\Prj_i[Q_0]\right\}
=V_i.
$$
\par Thanks to inclusion \rf{H-SM-R1}, $U_i\supset V_i$, $i=1,2$, proving that the orthogonal projections onto the coordinate axes of the left hand side of \rf{H-SM-LR} contain the corresponding projections of its right hand side. Because the left and right hand sides of \rf{H-SM-LR} are {\it rectangles} inclusion \rf{H-SM-LR} holds.
\par The proof of the lemma is complete.\bx

\begin{lemma}\lbl{R-12} Let $\Rc_1,\Rc_2\subset\Rf(\RT)$ be two families of rectangles in $\RT$. If each family has a non-empty intersection, then
$$
\dist\left(\,\bigcap_{\Pi\in\Rc_1}\Pi,
\bigcap_{\Pi\in\Rc_2}\Pi\right)=
\sup_{\Pi_1\in\,\Rc_1,\Pi_2\in\,\Rc_2}
\dist(\,\Pi_1,\Pi_2)\,.
$$
\end{lemma}
\par {\it Proof.} The definition of the uniform norm in $\RT$ and representation \rf{PR-RW} easily imply the following formula for the distance between rectangles $\Pi_1,\Pi_2\in\RCT$:
$$
\dist\left(\Pi_1,\Pi_2\right)=
\max\{\dist\left(\,\Prj_1[\Pi_1],\Prj_1[\Pi_2]\right),
\dist\left(\,\Prj_2[\Pi_1],\Prj_2[\Pi_2]\right)\}.
$$
\par It is also clear that for every family $\Rc$ of rectangles in $\RT$ with non-empty intersection, we have
$$
\Prj_i\left[\,\bigcap_{\Pi\in\Rc}\Pi\right]=
\bigcap_{\Pi\in\Rc}\Prj_i[\Pi],~~~~i=1,2.
$$
\par These observations reduces the problem for the one dimensional case which we prove below.
\par Let $\Ic_1,\Ic_2\subset\Ic(\R)$ be two families of intervals in $\R$. We assume that each family has a non-empty intersection. Our aim is to show that the following equality
\bel{DS-R1}
\dist\left(\,\bigcap_{I\in\Ic_1}I,
\bigcap_{I\in\Ic_2}I\right)=
\sup_{I_1\in\,\Ic_1,I_2\in\,\Ic_2}
\dist(I_1,I_2)
\ee
holds.
\par Clearly, $\dist(I_1,I_2)=\dist(I_1^{\cl},I_2^{\cl})$ for every $I_1\in\Ic_1$, $I_2\in\Ic_2$. (Recall that the sign ``$\cl$'' denotes the closure of a set.) It is also can be readily seen that for any family $\Ac$ of intervals with non-empty intersection, we have $(\cap\{I:I\in\Ac\})^{\cl}=\cap\{I^{\cl}:I\in\Ac\}$. These remarks show that without loss of generality we may assume that all intervals from the families $\Ic_1$ and $\Ic_2$ are {\it closed}.
\smsk
\par Clearly, the right hand side of \rf{DS-R1} is majorized by its left hand side. Let us prove the converse statement. We know that
\bel{HH}
H_1=\cap\{\,I:I\in\Ic_1\}\ne\emp~~~\text{and}~~~
H_2=\cap\{\,I:I\in\Ic_2\}\ne\emp.
\ee
\par Our aim is to prove that
\bel{D-R1}
\dist(H_1,H_2)\le \sup\{\,\dist(I_1,I_2):I_1\in\Ic_1,I_2\in\Ic_2\,\}.
\ee
\par Let $r=\sup\{\,\dist(I_1,I_2):I_1\in\Ic_1,I_2\in\Ic_2\,\}$. Because all sets from $\Ic_1$ and $\Ic_2$ are closed, we have
\bel{RST}
(I_1+rI_0)\cap I_2\ne\emp~~~~
\text{for every}~~~~I_1\in\Ic_1, I_2\in\Ic_2.
\ee
\par We recall that $H_1=\cap\{\,I:I\in\Ic_1\}$ is a non-empty set, so that, thanks to Lemma \reff{H-RB},
\bel{H1-DF}
H_1+rI_0=\cap\{I+rI_0:I\in\Ic_1\}.
\ee
\par Let us put
$\Kc=\{I+rI_0:I\in\Ic_1\}\cup\{I:I\in\Ic_2\}$ and prove that $\cap\{\,I: I\in\Kc\}\ne\emp$. (If this property holds, then, thanks to \rf{H1-DF}, $(H_1+rI_0)\cap H_2\ne\emp$ proving the required inequality \rf{D-R1}.)
\smsk
\par Our proof will rely on Helly's theorem. First, we note that, thanks to \rf{HH} and \rf{RST}, $I_1\cap I_2\ne\emp$ for every $I_1,I_2\in\Kc$.
\par We also note that, thanks to \rf{HH}, the family $\Kc$ has a non-empty intersection provided every member of $\Kc$ is an unbounded interval of the form $I=[a,+\infty)$, $a\in\R$. The same is true provided each interval $I\in\Kc$ is of the form $I=(-\infty,a]$, $a\in\R$. Therefore, we may assume that either $\Kc$ contains a {\it bounded} interval or there exist intervals $I_1,I_2\in\Kc$ of the forms $I_1=[a_1,+\infty)$ and $I_2=(-\infty,a_2]$ respectively. Clearly, in this case $I_1\cap I_2$ is a {\it bounded} set.
\par Thus, in this case, all conditions of Lemma \reff{H-R} are satisfied. Thanks to this lemma, there exists a point common to all members of $\Kc$, and the proof of Lemma \reff{R-12} is complete.\bx
\smsk
\par We will also need the following simple claim.
\begin{claim}\lbl{TWO} (i) Let $A$ and $B$ be two intervals in $\R$. Then
\bel{AB-END}
\dhf(A,B)=\max\{|\inf A-\inf B|,|\sup A-\sup B|\}.
\ee
(See also our convention \rf{INF-S} for the cases of $\inf A,\inf B=-\infty$ and $\sup A,\sup B=+\infty$.)
\smsk
\par (ii) Let $\Ac,\Bc\in\RCT$ be two bounded rectangles in $\RT$. Then
\bel{CAB-H}
\|\cent(\Ac)-\cent(\Bc)\|\le\dhf(\Ac,\Bc).
\ee
\end{claim}
\par {\it Proof.} (i) Clearly, $\inf I=\inf I^{\cl}$ and
$\sup I=\sup I^{\cl}$ for any interval $I\in\Ic(\R)$. Let us also note that $\dhf(A,B)=\dhf(A^{\cl},B^{\cl})$. These properties show that without loss of generality we may assume that $A$ and $B$ are {\it closed} intervals.
\par Let
$$
r=\dhf(A,B)~~~\text{and}~~~
\delta=\max\{|\inf A-\inf B|,|\sup A-\sup B|\}.
$$
Then, thanks to \rf{HD-DF}, $A+rI_0\supset B$ proving that $\sup A+r\ge \sup B$ and $\inf A-r\le \inf B$. By interchanging the roles of $A$ and $B$ we obtain also $\sup B+r\ge \sup A$ and $\inf B-r\le \inf A$ proving that $\delta\le r$.
\par Let us prove that $r\le \delta$. Suppose that both $A$ and $B$ are bounded, i.e., $A=[\inf A,\sup A]$ and
$B=[\inf B,\sup B]$. Let $\alpha\in[0,1]$, and let
$$
a_\alpha=\alpha\inf A+(1-\alpha)\sup A~~~\text{and}~~~
b_\alpha=\alpha\inf B+(1-\alpha)\sup B.
$$
Then $a_\alpha\in A$, $b_\alpha\in B$, and $|a_\alpha-b_\alpha|\le \delta$ proving that $\dist(a,B)\le r$ for every $a\in A$ and $\dist(b,A)\le r$ for every $b\in B$. Hence, $A+rI_0\supset B$ and $B+rI_0\supset A$, so that, thanks to \rf{HD-DF}, $r\le \delta$. In a similar way we prove this inequality whenever one of the intervals is unbounded. We leave the details to the interested reader as an easy exercise.
\par (ii) By orthogonal projecting to the coordinate axes, we can reduce the problem to the one dimensional case. In this case, given bounded intervals $A,B\in\Ic(\R)$, we have
$$
\cent(A)=(\inf A+\sup A)/2~~~\text{and}~~~
\cent(B)=(\inf B+\sup B)/2.
$$
This and inequality \rf{AB-END} imply the required inequality $|\cent(A)-\cent(B)|\le\dhf(A,B)$ proving the claim.\bx
\smsk
\par Let $\MR$ be a pseudometric space and let $\Tc:\Mc\to \RCT$ be a set-valued mapping. Given $\eta\ge 0$, we define a set-valued mapping on $\Mc$ by
\bel{TAU-1}
\Tc^{[1]}[x:\eta]=\bigcap_{z\in\Mc}\,
\left[\Tc(z)+\eta\,\rho(x,z)\,Q_0\right],~~~x\in\Mc.
\ee
\begin{proposition}\lbl{X2-C} (a) Suppose that either $\Mc$ is finite or each rectangle $\Tc(x)$, $x\in\Mc$, is closed. If
\bel{TAU-NE}
\Tc^{[1]}[x:\eta]\ne\emp~~~\text{for every}~~~x\in\Mc,
\ee
then
\bel{TAM-1}
\dhf\left(\Tc^{[1]}[x:\eta],\Tc^{[1]}[y:\eta]\right)
\le \eta\,\rho(x,y)~~~\text{for all}~~~x,y\in\Mc.
\ee
\par (b) Let us assume that either $\Mc$ is finite or all rectangles $\Tc(x)$, $x\in\Mc$, are closed and at least one of them is bounded. If
\bel{TD-DR}
\Tc(x)\cap\{\Tc(y)+\eta\,\rho(x,y)Q_0\}\ne\emp~~~\text{for all}~~~ x,y\in\Mc,
\ee
then properties \rf{TAU-NE} and \rf{TAM-1} hold.
\par Furthermore, if the set $\Tc^{[1]}[x:\eta]$ is bounded for every $x\in\Mc$, then the mapping
$\tau(x)=\cent\left(\Tc^{[1]}[x:\eta]\right)$, $x\in\Mc$,
is a Lipschitz selection of $\Tc$ with $\|\tau\|_{\Lip(\Mc)}\le\eta$.
\end{proposition}
\par {\it Proof.} (a) We know that $\Tc^{[1]}[x:\eta]\ne\emp$ so that, thanks to \rf{TAU-NE}, Lemma \reff{H-RB} and definition \rf{TAU-1}, we have
\be
\Tc^{[1]}[x:\eta]+\eta\,\rho(x,y)Q_0&=&
\left\{\bigcap_{z\in \Mc}
\left[\Tc(z)+\eta\,\rho(x,z)Q_0\right]\right\}
+\eta\,\rho(x,y)Q_0
\nn\\
&=&\bigcap_{z\in \Mc}
\left[\Tc(z)+(\eta\,\rho(x,z)+\eta\,\rho(x,y))\,Q_0\right].
\nn
\ee
Hence, thanks to the triangle inequality, we have
$$
\Tc^{[1]}[x:\eta]+\eta\,\rho(x,y)\,Q_0\supset
\bigcap_{z\in \Mc}\,
\left[\Tc(z)+\eta\,\rho(y,z)\,Q_0\right]=\Tc^{[1]}[y:\eta].
$$
\par By interchanging the roles of $x$ and $y$ we obtain also $\Tc^{[1]}[y:\eta]+\eta\,\rho(x,y)\,Q_0\supset \Tc^{[1]}[x:\eta]$. These two inclusions imply the required inequality \rf{TAM-1}.
\smsk
\par (b) Let us fix $x\in\Mc$ and prove that
\bel{FFP}
\{\Tc(z)+\eta\,\rho(x,z)\,Q_0\}
\cap\{\Tc(z')+\eta\,\rho(x,z')\,Q_0\}
\ne\emp~~~~~\text{for every}~~~z,z'\in\Mc.
\ee
\par Thanks to \rf{TD-DR}, there exist points $g(z)\in \Tc(z)$ and $g(z')\in \Tc(z')$ such that $\|g(z)-g(z')\|\le \eta\,\rho(z,z')$. From this and the triangle inequality, we have
$\|g(z)-g(z')\|\le \eta\,\rho(z,x)+\eta\,\rho(x,z')$. This
implies the existence of a point $w\in\RT$ such that
$\|g(z)-w\|\le \eta\,\rho(z,x)$ and $\|g(z')-w\|\le \eta\,\rho(z',x)$.
\par But $g(z)\in \Tc(z)$ and $g(z')\in \Tc(z')$ so that $w$ belongs to the left hand side of \rf{FFP} proving this property.
\msk
\par Let $\Kc=\{\Tc(z)+\eta\,\rho(x,z)\,Q_0:z\in\Mc\}$.
If $\Mc$ is infinite, then, thanks to the hypothesis of part (b), each member of $\Kc$ is closed. Furthermore, there exists $\tx\in\Mc$ such that the set $\Tc(\tx)$ is bounded so that the set
$\tK=\Tc(\tx)+\eta\,\rho(x,\tx)\,Q_0$ is bounded as well and belongs to the family $\Kc$. Combining these properties of $\Kc$ with \rf{FFP}, we conclude that this family satisfies the hypothesis of Lemma \reff{H-R}. Thanks to this lemma, the set $\Tc^{[1]}[x:\eta]=\cap\{K:K\in\Kc\}$ is non-empty which proves property \rf{TAU-NE}. In turn, inequality \rf{TAM-1} follows from part (a) of the proposition.
\msk
\par Finally, inequalities \rf{CAB-H} and \rf{TAM-1} imply the required inequality $\|\tau\|_{\Lip(\Mc)}\le\eta$ completing the proof of the proposition.\bx
\msk
\par Let us give several explicit formulae for Lipschitz selections of set-valued mappings in the one dimensional case. Let $F:\Mc\to \Ic(\R)$ be a set-valued mapping. We set
$$
a_F(x)=\inf F(x)~~~~\text{and}~~~~b_F(x)=\sup F(x).
$$
Thus, $a_F$ and $b_F$ are two functions on $\Mc$ such that
$$
a_F:\Mc\to\R\cup\{-\infty\},~~~
b_F:\Mc\to\R\cup\{+\infty\}~~~~\text{and}~~~~a_F(x)\le b_F(x)~~~\text{for all}~~~x\in\Mc.
$$
\par Clearly,
\bel{DR-XY}
\dist(F(x),F(y))=\max\{[a_F(x)-b_F(y)]_+,[a_F(y)-b_F(x)]_+\}.
\ee
(See our convention \rf{INF-S} for the case of $a_F(x)=-\infty$, $b_F(x)=+\infty$.)
\par Given $\eta\ge 0$, we introduce the following functions on $\Mc$:
\bel{FP-D1}
a^{[1]}_F[x:\eta]=\sup_{y\in\Mc}
\,\left\{a_F(y)-\eta\,\rho(x,y)\right\},~~~~~
b^{[1]}_F[x:\eta]=
\inf_{y\in\Mc}\,\left\{b_F(y)+\eta\,\rho(x,y)\right\}
\ee
and
\bel{FS-D1}
c_F[x:\eta]=\left(a^{[1]}_F[x:\eta]
+b^{[1]}_F[x:\eta]\right)/2.
\ee
\par We also define a set-valued mapping on $\Mc$ by
\bel{F1-PR1}
F^{[1]}[x:\eta]=
\bigcap_{z\in\Mc}\,\left[F(z)+
\eta\,\rho(x,z)\,\BXR\right],
~~~x\in\Mc.
\ee
\par Comparing this definition with definitions \rf{FP-D1} and \rf{FS-D1}, we conclude that for every $x\in\Mc$,
\bel{F1-ENDS}
a^{[1]}_F[x:\eta]=\inf F^{[1]}[x:\eta],~~~~~
b^{[1]}_F[x:\eta]=\sup F^{[1]}[x:\eta],
\ee
and
\bel{F1-CNT}
c_F[x:\eta]=\cent\left(F^{[1]}[x:\eta]\right)
\ee
provided $F^{[1]}[x:\eta]$ is {\it bounded}.
\smsk
\begin{remark}\lbl{RM-ABF} {\em We note that the function
\bel{FPL}
f^+=b^{[1]}_F[\cdot:\eta]~~~\text{maps}~~~\Mc~~~
\text{into}~~~\R
\ee
if and only if $f^+\nequiv+\infty$, i.e., $f^+(x^+)<\infty$ for some $x^+\in\Mc$. Analogously, the mapping
\bel{FMN}
f^-=a^{[1]}_F[\cdot:\eta]~~~\text{maps}~~~\Mc~~~
\text{into}~~~\R
\ee
if and only if $f^-\nequiv-\infty$. Finally, the mapping
$c^{[1]}_F[\cdot:\eta]=(f^++f^-)/2$, see \rf{FS-D1} and \rf{F1-CNT}, is well defined if and only if
both $f^+\nequiv+\infty$ and $f^-\nequiv-\infty$. \rbx}
\end{remark}
\smsk

\begin{proposition}\lbl{FP-R1} (i) (The Finiteness Principle for Lipschitz selections in $\R$.) Let
$F:\Mc\to \Ic(\R)$ be a set-valued mapping. Let us assume that either $\Mc$ is finite or all intervals $F(x)$, $x\in\Mc$ are closed and at least one
of them is bounded.
\par Let $\eta\ge 0$. Suppose that for every $x,y\in\Mc$ the restriction $F|_{\{x,y\}}$ of $F$ to $\{x,y\}$ has a Lipschitz selection $f_{\{x,y\}}$ with $\|f_{\{x,y\}}\|_{\Lip(\{x,y\},\R)}\le \eta$. Then $F$ has a Lipschitz selection $f:\Mc\to\R$ with Lipschitz seminorm $\|f\|_{\Lip(\Mc,\R)}\le \eta$.
\par Furthermore, one can set $f=c^{[1]}_F[\cdot:\eta]$ provided there exist $x^+,x^-\in\Mc$ such that $\inf F(x^-)>-\infty$ and $\sup F(x^+)<\infty$. If all intervals $F(x)$, $x\in\Mc$, are closed one can set $f=b^{[1]}_F[\cdot:\eta]$ if $F(x^+)$ is bounded from above for some $x^+\in\Mc$, or $f=a^{[1]}_F[\cdot:\eta]$ if $F(x^-)$ is bounded from below for some $x^-\in\Mc$.
\smsk
\par (ii) Suppose that all intervals $F(x)$, $x\in\Mc$ are closed and at least one of them is bounded. Let
$$
\lambda_F=\sup_{x,y\in\Mc}
[\inf F(x)-\sup F(y)]_+/\rho(x,y)
=\sup_{x,y\in\Mc}\dist(F(x),F(y))/\rho(x,y).
$$
(When calculating $\lambda_F$, we use convention \rf{INF-S}. We also note that the second equality in this definition is due to \rf{DR-XY}.)
\par Then $F$ has a Lipschitz selection if and only if $\lambda_F<\infty$. Moreover, in this case, $\lambda_F=\FM$, see \rf{FM}, and each of the mappings $f=a^{[1]}_F[\cdot:\eta]$, $f=b^{[1]}_F[\cdot:\eta]$ and $f=c^{[1]}_F[\cdot:\eta]$, provides an optimal Lipschitz selection of $F$, i.e.,  $\|f\|_{\Lip(\Mc,\R)}=\FM$.
\end{proposition}
\par {\it Proof.} (i) Thanks to the hypothesis of part (i) of the proposition and part (a) and (b) of Proposition \reff{X2-C}, we have
$$
F^{[1]}[x:\eta]\ne\emp~~~\text{for every}~~~x\in\Mc,
$$
and
$$
\dhf\left(F^{[1]}[x:\eta],F^{[1]}[y:\eta]\right)
\le \eta\,\rho(x,y)~~~\text{for all}~~~x,y\in\Mc.
$$
From this, part (i) of Claim \reff{TWO}, and definitions \rf{F1-ENDS} it follows that the inequality
\bel{AB-HDF}
\max\left\{\,|a^{[1]}_F[x:\eta]-a^{[1]}_F[y:\eta]|,\,
|b^{[1]}_F[x:\eta]-b^{[1]}_F[y:\eta]|\,\right\}
\le \eta\,\rho(x,y)
\ee
holds for all $x,y\in\Mc$.
\par Clearly, if $F\equiv \R$ then the constant mapping $f\equiv \{0\}$ on $\Mc$ is a Lipschitz selection of $F$ (with $\|f\|_{\Lip(\Mc,\R)}=0$). Otherwise, either $f^+=b^{[1]}_F[\cdot:\eta]\nequiv+\infty$ or $f^-=a^{[1]}_F[\cdot:\eta]\nequiv-\infty$. Therefore, thanks to Remark \reff{RM-ABF}, either $f^+:\Mc\to\R$
or $f^-:\Mc\to\R$. See \rf{FPL} and \rf{FMN}.
\par Let us note that if each set $F(x)$ is {\it closed}, and at least on of the sets $F(x)$ is bounded, then the set $F^{[1]}[x:\eta]$ is closed and bounded as well. In this case, the points
\bel{AB-INF}
a^{[1]}_F[x:\eta],b^{[1]}_F[x:\eta]\in F^{[1]}[x:\eta]\subset F(x).
\ee
See definition \rf{F1-ENDS}. Therefore, in this case we can set either $f=f^+$ or $f=f^-$. Then, thanks to \rf{AB-HDF}, in both cases the mapping $f:\Mc\to\R$ will be a Lipschitz selection of $F$ with $\|f\|_{\Lip(\Mc,\R)}\le\eta$. Also from this it follows that the mapping $f=c_F[\cdot:\eta]=(f^+ + f^-)/2$, see \rf{F1-CNT}, has the same properties.
\par Now, let $\Mc$ be {\it finite}. In this case we can not guarantee that the set $F^{[1]}[x:\eta]$ is closed and property \rf{AB-INF} holds. However, if   $f^+=b^{[1]}_F[\cdot:\eta]\nequiv+\infty$ but  $f^-=a^{[1]}_F[\cdot:\eta]\equiv-\infty$, then each interval $F^{[1]}[x:\eta]$ is {\it unbounded from below}. Because $\Mc$ is finite, all these intervals have a common point, say $A$. Then the constant mapping $f\equiv \{A\}$ is a Lipschitz selection of $F$ (with $\|f\|_{\Lip(\Mc,\R)}=0$). Analogously, if $f^-\nequiv-\infty$ but $f^+\equiv+\infty$, there is a constant mapping which provides a Lipschitz selection of $F$.
\par Let us suppose that both $f^-=a^{[1]}_F[\cdot:\eta]\nequiv-\infty$ and $f^-=a^{[1]}_F[\cdot:\eta]\nequiv-\infty$. In this case, thanks to Remark \reff{RM-ABF}, the mapping
$c^{[1]}_F[\cdot:\eta]=(f^++f^-)/2$, see \rf{FS-D1} and \rf{F1-CNT}, is well defined, i.e., each set $F^{[1]}[x:\eta]$ is non-empty and bounded (but not necessarily closed!). Clearly,
$$
f(x)=c^{[1]}_F[x:\eta]\in F^{[1]}[x:\eta]\subset F(x)
~~~~\text{for every}~~~~x\in\Mc,
$$
proving that $f$ is a selection of $F$. Thanks to \rf{AB-HDF}, its Lipschitz seminorm  $\|f\|_{\Lip(\Mc,\R)}\le\eta$, and the proof of part (i) of the proposition is complete.
\msk
\par (ii) This criterion for Lipschitz selections is immediate from part (i) of the proposition. We leave the details of the proof to the interested reader.\bx
\msk
\par Part (i) of Proposition \reff{FP-R1} implies the following Finiteness Principle for rectangles in $\RT$.
\begin{proposition}\lbl{FP-RC} Let $\eta\ge 0$. Let $\MR$ be a pseudometric space, and let $\Tc:\Mc\to \RCT$ be a set-valued mapping. Let us assume that either $\Mc$ is finite or all rectangles $\Tc(x)$, $x\in\Mc$, are closed and at least one of them is bounded.
\par Suppose that for every $x,y\in\Mc$ the restriction $\Tc|_{\{x,y\}}$ of $\Tc$ to $\{x,y\}$ has a Lipschitz selection $\tau_{\{x,y\}}$ with $\|\tau_{\{x,y\}}\|_{\Lip(\{x,y\})}\le \eta$.
\par Then $\Tc$ has a Lipschitz selection $\tau:\Mc\to\RT$ with Lipschitz seminorm $\|\tau\|_{\Lip(\Mc)}\le \eta$.
\end{proposition}
\par {\it Proof.} By orthogonal projecting onto the coordinate axes, we reduce the problem to the one dimensional case, i.e., to the Finiteness Principle for Lipschitz selections in $\R$ proven in part (i) of Proposition \reff{FP-R1}.\bx
\smsk

\SECT{3. The key theorem: Lipschitz selections and rectangular hulls.}{3}
\addtocontents{toc}{3. The key theorem: Lipschitz selections and rectangular hulls.\hfill \thepage\par\VST}

\indent
\par Let $\MS=(\Mc,\rho)$ be a pseudometric space, and let
$F:\Mc\to \CNV$ be a set-valued mapping. Given $\lambda\ge 0$ and $x,x',x''\in\Mc$, we let   $\Wc_F[x,x',x'':\lambda]$ denote a (possibly empty) subset of $\RT$ defined by
\bel{WC-DF}
\Wc_F[x,x',x'':\lambda]=
\HR[\{F(x')+\lambda\,\rho(x',x)\,Q_0\}
\cap \{F(x'')+\lambda\,\rho(x'',x)\,Q_0\}].
\ee
\par We recall that by $\HR[\cdot]$ we denote the rectangular hull of a set (see \rf{HRS}), and by  $\RL[\cdot,\cdot:\lambda]$ the rectangle defined by formula \rf{RL}. Note that from \rf{RL} and \rf{WC-DF}, we have
\bel{RW-C}
\RL[x,x':\lambda]=\Wc_F[x,x,x':\lambda]~~~~\text{for every}~~~~x,x'\in\Mc.
\ee
\par The proofs of the necessity parts of Lipschitz selection criteria presented in this section rely on the following proposition.
\begin{proposition}\lbl{P-GR-RT} Let $F:\Mc\to \CNV$ be a set-valued mapping and let $\lambda\ge 0$. Suppose that $F$ has a Lipschitz selection $f:\Mc\to \RT$ with $\|f\|_{\Lip(\Mc)}\le\lambda$. Then for every $x,x',x'',y,y',y''\in\Mc$, we have
\smsk
\par (i)
$\RL[x,x':\lambda]\cap \{\RL[y,y':\lambda]+\lambda\,\rho(x,y)\,Q_0\}\ne\emp$;
\smsk
\par (ii) $\Wc_F[x,x',x'':\lambda]\cap \{\Wc_F[y,y',y'':\lambda]+\lambda\,\rho(x,y)\,Q_0\}\ne\emp$.
\end{proposition}
\par {\it Proof.} Clearly, part (i) of the proposition follows from part (ii) and \rf{RW-C}. Let us prove part (ii).
\par Because $f$ is a Lipschitz selection of $F$, for every $x,x',x''\in\Mc$ we have $f(x)\in F(x)$,
$f(x')\in F(x')$,  $f(x'')\in F(x'')$,
$$
\|f(x)-f(x')\|\le\lambda\,\rho(x,x')~~~~\text{and}~~~~
\|f(x)-f(x'')\|\le\lambda\,\rho(x,x'').
$$
\par Hence,
\be
f(x)&\in& \{F(x')+\lambda\,\rho(x,x')Q_0\}\cap
\{F(x'')+\lambda\,\rho(x,x'')Q_0\}\nn\\
&\subset&
\HR[\{F(x')+\lambda\,\rho(x,x')Q_0\}\cap
\{F(x'')+\lambda\,\rho(x,x'')Q_0\}]=
\Wc_F[x,x',x'':\lambda].
\nn
\ee
\par In the same fashion we show that $f(y)\in \Wc_F[y,y',y'':\lambda]$. These properties and inequality
$\|f(x)-f(y)\|\le \lambda\,\rho(x,y)$ imply (ii) proving the proposition.\bx
\msk
\par The main result of the present section is the following theorem.
\begin{theorem}\lbl{W-CR} Let $\Mf=\MR$ be a pseudometric space, and let $F:\Mc\to\CRT$ be a set-valued mapping satisfying Condition \reff{CND-T}.
\par Given non-negative constants $\tlm$ and $\lambda$, let us assume that for every $x,x',x'',y,y',y''\in\Mc$, we have
\bel{WNEW}
\Wc_F[x,x',x'':\tlm]\cap \{\Wc_F[y,y',y'':\tlm]+\lambda\,\rho(x,y)\,Q_0\}\ne\emp.
\ee
\par Then $F$ has a Lipschitz selection $f:\Mc\to\RT$ with $\|f\|_{\Lip(\Mc)}\le 2\lambda+\tlm$.
\end{theorem}
\par We refer to this result as {\it the key theorem.}
We start the proof of the key theorem in this section and complete it at the end of the next section.
\smsk
\par {\it Proof of the key theorem.} Suppose that for every $x,x',x'',y,y',y''\in\Mc$ condition \rf{WNEW} holds. Let us construct a Lipschitz selection $f:\Mc\to\RT$ of $F$ with Lipschitz seminorm $\|f\|_{\Lip(\Mc)}\le 2\lambda+\tlm$. We will do this in three steps.
\msk
\par \underline{\sc The First Step.} We introduce a set-valued mapping on $\Mc$ defined by the formula
\bel{WF1-3L}
F^{[1]}[x:\tlm]=
\bigcap_{y\in \Mc}\,\left[F(y)+\tlm\,\rho(x,y)\,Q_0\right],
~~~~x\in\Mc.
\ee
\begin{lemma}\lbl{WF1-NE} For each $x\in\Mc$ the set
$F^{[1]}[x:\tlm]$ is a non-empty closed convex subset of $\RT$. Moreover, for every $x\in\Mc$ the following representation holds:
\bel{WH-F1}
\HR[F^{[1]}[x:\tlm]]=\cap\{\Wc_F[x,y,y':\tlm]: y,y'\in\Mc\}.
\ee
\par Furthermore, if $\Mc$ is infinite, then the set
$F^{[1]}[x:\tlm]$ is bounded for all $x\in\Mc$.
\end{lemma}
\par {\it Proof.} Let us prove that
\bel{F1-IN}
F^{[1]}[x:\tlm]\ne\emp~~~~\text{for every}~~~~x\in\Mc.
\ee
\par Given $x\in\Mc$, we set
$$
\Cf_x=\{F(y)+\tlm\,\rho(x,y)\,Q_0:y\in\Mc\}.
$$
Then $F^{[1]}[x:\tlm]=\cap\{C:C\in\Cf_x\}$. See \rf{WF1-3L}.
\par Let us prove that for every $y_1,y_1',y_2,y_2'\in\Mc$ the sets
\bel{C-PT}
C_i=F(y_i)+\tlm\,\rho(x,y_i)Q_0~~~~\text{and}~~~~ C'_i=F(y'_i)+\tlm\,\rho(x,y'_i)Q_0,~~ i=1,2,
\ee
satisfy property \rf{P1-C}.
\par First, let us note that, thanks to \rf{WNEW}, the set
$$
\Wc_F[x,y,z:\tlm]\ne\emp~~~~\text{for all}~~~ x,y,z\in\Mc.
$$
In particular, from this and definition \rf{WC-DF}, it follows that
$$
\{F(y)+\tlm\,\rho(x,y)\,Q_0\}\cap
\{F(z)+\tlm\,\rho(x,z)\,Q_0\}\ne\emp~~~~
\text{for every}~~~~y,z\in\Mc
$$
proving that {\it any two elements of $\Cf_x$ have a common point}.
\par Property \rf{WNEW} tells us that
$$
\Wc_F[x,y_1,y'_1:\tlm]\cap\Wc_F[x,y_2,y'_2:\tlm]
\ne\emp~~~~\text{for every}~~~~y_1,y'_1,y_2,y'_2\in\Mc.
$$
Thanks to \rf{WC-DF} and \rf{C-PT}, $\Wc_F[x,y_1,y'_1:\tlm]=\Hc[C_1\cap C_1']$ and $\Wc_F[x,y_2,y'_2:\tlm]=\Hc[C_2\cap C_2']$ proving that
$$
\Hc[C_1\cap C_1']\,\cbg \Hc[C_2\cap C_2'] \ne\emp~~~~~\text{for every}~~~~~ C_1,C_1',C_2,C_2'\in\Cf.
$$
Hence,
$$
\Prj_1[\Hc[C_1\cap C_1']]\,\cbg \Prj_1[\Hc[C_2\cap C_2']] \ne\emp
$$
so that, thanks to \rf{HS-U},
$\Prj_1[C_1\cap C_1']\,\cbg \Prj_1[C_2\cap C_2'] \ne\emp$
proving \rf{P1-C}.
\par Now let us assume that the set $\Mc$ is {\it finite}.
In this case, the conditions (i) and (ii) of Proposition \reff{INT-RE} are satisfied. This proposition tells us that
$\cap\{C:C\in\Cf\}\ne\emp$ proving \rf{F1-IN} in the case under consideration.
\smsk
\par Now let $\Mc$ be infinite. In this case, thanks to the hypothesis of Theorem \reff{W-CR} and Condition \reff{CND-T}, there exist a constant $\alpha\ge 0$ and a finite set
$\hM=\{\hx_1,...,\hx_m\}\subset\Mc$
such that the set
$$
\cap\{F(y)+\alpha Q_0:y\in\hM\}~~~~\text{is non-empty and bounded}.
$$
\par Let $x\in\Mc$ and let
\bel{HCF}
\hCF_x=\{F(y)+\tlm\rho(x,y)Q_0:y\in\hM\}.
\ee
Let us prove that the set
\bel{N-B}
G_x=\cap\{C:C\in\hCF_x\}~~~~\text{is non-empty and bounded}.
\ee
\par As we have proved above, $G_x\ne\emp$ (because $\hM$ is finite). Let us see that $G_x$ is bounded.
\par Suppose that the set $G_x$ is unbounded.
We recall two properties of subsets from $\CRT$:
\smsk
\par ($\bigstar 1$) Let $K\in\CRT$. If $K$ is unbounded than it contains a ray;
\par ($\bigstar 2$) Let $K\in\CRT$, $h\in \SO$, and let $L_h=\{th:t\ge 0\}$ be the ray emanating from the origin in the direction of $h$. Let $x,y\in K$. Then $x + L_h\subset K$ if and only $y + L_h\subset K$.
\smsk
\par See, e.g.  \cite[Section 2.5, Lemma 1, 2 ]{Gr-2003}.
\smsk
\par Fix a point $\ta\in G_x$. The set $G_x\in \CRT$ and  unbounded, so that, thanks to ($\bigstar 1$), there exists a vector $h\in \SO$ such that the ray starting at $\ta$ in the direction of $h$ lies in $G$. Thus, $\ta+L_h\subset G_x$ where $L_h=\{th:t\ge 0\}$. Therefore, thanks to definition of $G_x$ (see \rf{N-B}) and \rf{HCF}, for every $y\in\hM$,
$$
\ta+L_h\subset F(y)+\tlm\rho(x,y)Q_0.
$$
\par Furthermore, combining this property with property ($\bigstar 2$), we obtain the following: {\it for every} $b\in F(y)+\tlm\rho(x,y)Q_0$ we have
\bel{B-GI}
b+L_h\subset F(y)+\tlm\rho(x,y)Q_0.
\ee
In particular, this property holds for every $b\in F(y)$.
\par Let us prove that property \rf{B-GI} implies a stronger one:
\bel{B-FG}
b+L_h\subset F(y)~~~~\text{for every}~~~~
b\in F(y).
\ee
\par Indeed, let $p\in b+L_h$, i.e., there exits $s\ge 0$ such that $p=b+sh$. We know that
$$
b+th\in F(y)+\tlm\rho(x,y)Q_0~~~~\text{for every}~~~~
t\ge s
$$
so that there exists a point $b_t\in F(y)$ such that $\|b_t-(b+th)\|\le \tlm\rho(x,y)$.
\par Furthermore, because $F(y)$ is convex, the line segment $[b,b_t]\subset F(y)$ so that the point
$$
u_t=b+\frac{s}{t}(b_t-b)\in F(y).
$$
Hence,
$$
\|u_t-p\|=\|\frac{s}{t}(b_t-b)-sh\|=
\|\frac{s}{t}(b_t-(b+th))\|\le \frac{s}{t}(\tlm\rho(x,y))
$$
proving that $u_t\to p$ as $t\to\infty$. But $u_t\in F(y)$ and the set $F(y)$ is closed, so that $p\in F(y)$ proving
\rf{B-FG}.
\par In particular, from \rf{B-FG} it follows that $b+L_h\subset F(y)+\alpha Q_0$ provided $b\in F(y)$. Therefore, thanks to property ($\bigstar 2$), for every $y\in\hM$,
\bel{BFP}
b+L_h\subset F(y)+\alpha Q_0~~~\text{for all}~~~b\in F(y)+\alpha Q_0.
\ee
\par We recall that there exists a point, say $b_0$, common to all of the sets $\{F(y)+\alpha Q_0:y\in\hM\}$. Therefore, thanks to \rf{BFP},
$$
b_0+L_h\subset F(y)+\alpha Q_0~~~\text{for every}~~~y\in\hM
$$
proving that
$b_0+L_h\subset \tF=\{F(y)+\alpha Q_0:y\in\hM\}$. Thus, the set $\tF$ is unbounded, a contradiction.
\smsk
\par This contradictions proves that the set $G_x$ defined by \rf{N-B}, is bounded.
\smsk
\par Thus, the family $\Cf_x$ satisfies the hypothesis of Proposition \reff{INT-RE}. This proposition tells us that for every $x\in\Mc$ the set $F^{[1]}[x:\tlm]=\cap\{C:C\in\Cf_x\}$ is non-empty. Note that $F^{[1]}[x:\tlm]\subset G_x$ so that $F^{[1]}[x:\tlm]$ is {\it bounded} for each $x\in\Mc$ provided $\Mc$ is infinite.
\par Finally, formula \rf{PH-C} and definition \rf{WC-DF} imply formula \rf{WH-F1}, and the proof of the lemma is complete.\bx
\msk
\par Given $x\in\Mc$, we let $\Tc_F(x)$ denote the rectangular hull of the set $F^{[1]}[x:\tlm]$. Thus, $\Tc_F:\Mc\to\RCT$ is a set-valued mapping defined by
\bel{WTH-F1}
\Tc_F(x)=\HR[F^{[1]}[x:\tlm]]=
\HR\left[\cap\left\{F(y)+\tlm\,\rho(x,y)\,Q_0:y\in \Mc\right\}\right], ~~~~~x\in\Mc.
\ee
Let us note that formula \rf{WH-F1} provides the following representation of the mapping $\Tc_F$:
\bel{WR-TC}
\Tc_F(x)=\cap\{\Wc_F[x,x',x'':\tlm]: x',x''\in\Mc\},~~~~x\in\Mc.
\ee
\begin{remark}\lbl{RM-TAU} {\em Let
$\Tf=\{\Tc_F(x):x\in\Mc\}$. Let us note the following properties of this family:
\smsk
\par (i) Each member of the family $\Tf$ is a {\it non-empty rectangle} in $\RT$. Indeed, thanks to Lemma \reff{WF1-NE}, $F^{[1]}[x:\tlm]\ne\emp$ for every $x\in\Mc$ so that, thanks to \rf{WTH-F1}, $\Tc_F(x)=\HR[F^{[1]}[x:\tlm]]\ne\emp$ as well;
\smsk
\par (ii) Either the family $\Tf$ is {\it finite} or every rectangle from $\Tf$ is a {\it compact set}. Indeed, if
$\Mc$ is infinite, then, thanks to Lemma \reff{WF1-NE},
the set $F^{[1]}[x:\tlm]$ is compact. Therefore, its orthogonal projections onto the coordinate axes are compact, so that, thanks to \rf{P1XP2} and \rf{WTH-F1}, the rectangle $\Tc_F(x)$ is compact as well.\rbx}
\end{remark}
\msk

\par \underline{\sc The Second Step.} At this step we prove the existence of a Lipschitz selection of the set-valued mapping $\Tc_F$.
\begin{proposition}\lbl{WLS-T} (i) The set-valued mapping $\Tc_F:\Mc\to\RCT$ has a Lipschitz selection $g:\Mc\to\RT$ with Lipschitz seminorm $\|g\|_{\Lip(\Mc)}\le \lambda$;
\smsk
\par (ii) We let $\Tc_F^{[1]}[\cdot:\lambda]$ denote a set-valued mapping on $\Mc$ defined by
\bel{T-ONE}
\Tc_F^{[1]}[x:\lambda]=
\bigcap_{z\in\Mc}\,\left[\Tc_F(z)+
\lambda\,\rho(x,z)\,Q_0\right],~~~x\in\Mc.
\ee
\par Then
\bel{TC-NEMP}
\Tc_F^{[1]}[x:\lambda]\ne\emp~~~~\text{for every}~~~x\in\Mc.
\ee
\par Furthermore, the following property
\bel{DRT-1}
\dhf(\Tc_F^{[1]}[x:\lambda],\Tc_F^{[1]}[y:\lambda])\le \lambda\,\rho(x,y)~~~~\text{for every}~~~~x,y\in\Mc;
\ee
holds. (Recall that $\dhf$ denotes the Hausdorff distance between sets.)
\smsk
\par (iii) If each rectangle $\Tc_F^{[1]}[x:\lambda]$, $x\in\Mc$, is bounded, then the mapping
$$
g_F(x)=\cent(\Tc_F^{[1]}[x:\lambda]),~~~~x\in\Mc,
$$
is a Lipschitz selection of~ $\Tc_F$ with  $\|g_F\|_{\Lip(\Mc)}\le\lambda$.
\end{proposition}
\par {\it Proof.} {\it (i)} Remark \reff{RM-TAU} tells us that the mapping $\Tc=\Tc_F$ satisfies the hypothesis of Proposition \reff{FP-RC}. Thanks to this proposition,
the required Lipschitz selection $g$ exists provided  for every $x,y\in\Mc$ the restriction $\Tc_F|_{\{x,y\}}$ of $\Tc_F$ to $\{x,y\}$ has a Lipschitz selection $g_{\{x,y\}}$ with Lipschitz  seminorm $\|g_{\{x,y\}}\|_{\Lip(\{x,y\})}\le \lambda$. Clearly, this requirement is equivalent to the following property:
\bel{WDTC}
\Tc_F(x)\cap\{\Tc_F(y)+\lambda\,\rho(x,y)Q_0\}\ne\emp
~~~~\text{for every}~~~~x,y\in\Mc.
\ee
\par Let
$$
\Tf_x=\{\Wc_F[x,x',x'':\tlm]: x',x''\in\Mc\} ~~~\text{and}~~~
\Tf_y=\{\Wc_F[y,y',y'':\tlm]: y',y''\in\Mc\}.
$$
Thanks to \rf{WR-TC},
\bel{TF-XY}
\Tc_F(x)=\cap\{W:W\in\Tf_x\}~~~~\text{and}~~~~
\Tc_F(y)=\cap\{W:W\in\Tf_y\}.
\ee
\par First, let us prove \rf{WDTC} provided $\Mc$ is {\it finite}. In this case, $\Tf_x$ and $\Tf_y$ are {\it finite} families of rectangles. Furthermore, thanks to part (i) of Remark \reff{RM-TAU}, each family has a non-empty intersection.
\par Let $r=\lambda\,\rho(x,y)$. Then, thanks to \rf{TF-XY} and Lemma \reff{H-RB},
$$
\Tc_F(y)+\lambda\,\rho(x,y)Q_0=\cap\{W:W\in\Tf_y\}+rQ_0=
\cap\{W+rQ_0:W\in\Tf_y\}
$$
so that
\bel{WQ-NM}
\Tc_F(x)\cap\{\Tc_F(y)+\lambda\,\rho(x,y)Q_0\}=
[\cap\{W:W\in\Tf_x\}]\cap [\cap\{W+rQ_0:W\in\Tf_y\}].
\ee
\par Let $\Tfw=\Tf_x\cup \Tf^+_y$ where $\Tf^+_y=\{W+rQ_0:W\in\Tf_y\}$. Thanks to \rf{WQ-NM}, property \rf{WDTC} holds provided the family of rectangles
$\Tfw$ has a common point. Because $\Tf_x$ and $\Tf_y$ are {\it finite} families, the family $\Tfw$ is finite as well.
Therefore, thanks to Helly's intersection theorem for rectangles, see Lemma \reff{H-R}, there exists a point common to all the family $\Tfw$ provided $W'\cap W''\ne\emp$ for every $W',W''\in\Tfw$.
\par Clearly, $W'\cap W''\ne\emp$ if $W',W''\in\Tf_x$ or
$W',W''\in\Tf^+_y$ because both $\Tf_x$ and $\Tf^+_y$ has a non-empty intersection. Let $W'=\Wc_F[x,x',x'':\tlm]$, $x',x''\in\Mc$, and let $W''=\Wc_F[y,y',y'':\tlm]+rQ_0$, $y',y''\in\Mc$, be two arbitrary members of $\Tf_x$ and  $\Tf^+_y$ respectively. Then, thanks to assumption \rf{WNEW} of Theorem \reff{W-CR}, $W'\cap W''\ne\emp$.
\par Thus, the hypothesis of Lemma \reff{H-R} holds for
$\Tfw$ so that this family has a common point. This proves that \rf{WDTC} holds provided $\Mc$ is finite.
\smsk
\par Now, let $\Mc$ be an infinite set. Remark \reff{RM-TAU} tells us that in this case the rectangles
$\Tc_F(x)$ and $\Tc_F(y)$ are {\it compact sets}. Therefore, property \rf{WDTC} is equivalent to the following inequality:
\bel{D-MA}
\dist(\Tc_F(x),\Tc_F(y))\le \lambda\,\rho(x,y).
\ee
\par We recall that, thanks to \rf{TF-XY}, $\Tf_x$ and $\Tf_y$ are two families of rectangles with non-empty intersections (equal to $\Tc_F(x)$ and $\Tc_F(y)$ respectively). Lemma \reff{R-12} tells us that in this case
\bel{FT-OP}
\dist(\Tc_F(x),\Tc_F(y))=
\dist\left(\,\bigcap_{W\in\Tf_x}W,
\bigcap_{W\in\Tf_y}W\right)=
\sup\{\,\dist(W',W''):
W'\in\Tf_x\,,W''\in\Tf_y\}.
\ee
\par But, thanks to assumption \rf{WNEW}, for every rectangle $W'=\Wc_F[x,x',x'':\tlm]\in \Tf_x$ and every rectangle $W''=\Wc_F[y,y',y'':\tlm]\in\Tf_y$ we have  $\dist(W',W'')\le \lambda\,\rho(x,y)$. This inequality and \rf{FT-OP} imply \rf{D-MA} proving the required property \rf{WDTC} and part {\it (i)} of the proposition.
\smsk
\par Let us prove parts {\it (ii)} and {\it (iii)}. We note that, thanks to property \rf{WDTC}, the mapping $\Tc=\Tc_F$ satisfies the conditions of the hypothesis of part (b) of Proposition \reff{X2-C} with $\eta=\lambda$. This proposition tells us that property \rf{TC-NEMP} and inequality \rf{DRT-1} hold proving part {\it (ii)}. Furthermore, part (b) of Proposition \reff{X2-C} proves part {\it (iii)}.
\par The proof of Proposition \reff{WLS-T} is complete.\bx
\bsk

\par \underline{\sc The Third Step.} At this step we construct a Lipschitz selection $f$ of the set-valued mapping $F$ with Lipschitz seminorm at most $2\lambda+\tlm$.
\msk
\par First, we recall that the set-valued mapping $F^{[1]}[\cdot:\tlm]$ and its rectangular hull, the set-valued mapping $\Tc_F=\HR[F^{[1]}[\cdot:\tlm]]$, are defined by formulae \rf{WF1-3L} and \rf{WTH-F1} respectively.
\par Part {\it (i)} of Proposition \reff{WLS-T} tells us that $\Tc_F$ has a Lipschitz selection with Lipschitz seminorm at most $\lambda$. As we have noted above, we cannot guarantee that the rectangle $\Tc_F(x)$ is a {\it closed set} for all $x\in\Mc$. Sometimes this leads to a certain complication of our constructive algorithm for a nearly optimal Lipschitz selection of $F$. To avoid these technical difficulties, we will work with the rectangles $\Tc_F(x)^{\cl}(x)$ (i.e., with {\it the closures} of $\Tc_F(x)$) rather than the rectangles $\Tc_F(x)$ themselves.
\smsk
\par Proposition \reff{WLS-T} tells us that
the set-valued mapping $\Tc_F:\Mc\to\RCT$ has a Lipschitz selection with Lipschitz seminorm at most $\lambda$. Clearly, $\Tc_F(x)(x)\subset\Tc_F(x)^{\cl}(x)$ so that
the set-valued mapping $\Tc_F^{\cl}$ also has a Lipschitz selection with Lipschitz seminorm at most $\lambda$. In other words, there exists a mapping $g:\Mc\to\RT$ such that
\bel{WG-HF}
g(x)\in \Tc_F^{\cl}(x)=\HR[F^{[1]}[x:\tlm]]^{\cl}~~~~~\text{for every}~~~~~x\in\Mc,
\ee
and
\bel{G-LIP}
\|g(x)-g(y)\|\le \lambda\,\rho(x,y)~~~~~\text{for all}~~~~~x,y\in\Mc.
\ee

\begin{proposition}\lbl{LSEL-F} Let $g:\Mc\to\RT$ be an arbitrary Lipschitz selection of the set-valued mapping $\Tc_F^{\cl}:\Mc\to\RCT$ with Lipschitz seminorm at most $\lambda$, i.e., a mapping satisfying conditions \rf{WG-HF} and \rf{G-LIP}.
\smsk
\par We define a mapping $f:\Mc\to\RT$ by letting
\bel{FX-PR}
f(x)=\Prm\left(g(x),F^{[1]}[x:\tlm]\right),~~~~~~x\in\Mc.
\ee
(Recall that $\Prm(\cdot,S)$ is the operator of metric projection onto a convex closed $S\subset\RT$. See \rf{MPR}.)
\smsk
\par Then the following properties hold:
\smsk
\par ($\bigstar 1$) The mapping $f$ is well defined, i.e., $f(x)$ is a singleton for every $x\in\Mc$. In this case
$$
f(x)=\Prm\left(g(x),F^{[1]}[x:\tlm]\right)\in F^{[1]}[x:\tlm]\subset F(x)~~~~\text{for every}~~~~x\in\Mc,
$$
so that $f$ is a {\it selection of $F$ on $\Mc$};
\smsk
\par ($\bigstar 2$) The mapping $f:\Mc\to\RT$ is Lipschitz with Lipschitz seminorm $\|f\|_{\Lip(\Mc)}\le 2\lambda+\tlm$.
\end{proposition}
\par The proof of this proposition is based on a number of auxiliary results. The first of these is the following lemma.
\begin{lemma}\lbl{WMP-S} Let $S\subset\RT$ be a non-empty convex closed set. Then for every point $a\in\HR[S]^{\cl}$ the metric projection $\Prm(a,S)$ is a singleton. Furthermore, $\Prm(a,S)$ coincides with a vertex of the square $Q(a,\dist(a,S))$.
\end{lemma}
\par {\it Proof.} Our proof is a slight modification of the proof of this lemma for the special case $a\in\HR[S]$ given in \cite[p. 301]{S-2002}. See also \cite[p. 67]{S-2020-L}.
\smsk
\par Clearly, if $a\in S$, nothing to prove. Suppose $a\notin S$ so that $r=\dist(a,S)>0$. Because $S$ is closed, $\Prj(a;S)\ne\emp$. Furthermore, $\Prj(a;S)=S\cap Q=S\cap \partial Q$ where  $Q=Q(a,r)$.
\par Because $\Prj(a;S)$ is {\it a non-empty convex set} lying on the boundary of $Q$, it belongs to a certain side of the square $Q$. In other words, there exist two distinct vertices of $Q$, say $A$ and $B$, such that $\Prj(a;S)\subset [A,B]$.
\par Let us prove that
\bel{AB-V}
\text{either}~~~~\Prj(a;S)=\{A\}~~~~\text{or}~~~~ \Prj(a;S)=\{B\}.
\ee
\par Indeed, otherwise there exists a point
$p\in (A,B)\cap\Prj(a;S)$. Let $\ell$ be the straight line passing through $A$ and $B$. Clearly, $\ell$ is parallel to a coordinate axis. Let $H_1,H_2$ be the closed half-planes determined by $\ell$. (Thus $\ell=H_1\cap H_2$.) Clearly, $Q$ is contained in one of these half-planes, say in $H_1$. Then $a\in H^{int}_1$ where $H^{int}_1$ denotes the interior of $H_1$ (because $\dist(a,\ell)=r>0$).
\smsk
\par Prove that in this case $S\subset H_2$, i.e., the straight line $\ell$ separates (not strictly) the square $Q$ and the set $S$. Indeed, suppose that there exists  a point $b\in S\cap H^{int}_1$. Then also $(p,b]\subset H^{int}_1$ because $p\in\partial H_1=\ell$. But $p\in(A,B)$ so that $(p,b]\cap Q^{int}\ne\emp$. On the other hand, because $S$ is convex and $p\in\partial S$, the interval $(p,b]\subset S$ proving that $S\cap Q^{int}\ne\emp$. But $S\cap Q \subset \partial Q$, a contradiction.
\par Thus, $S\subset H_2$ and $Q\subset H_1$. But $a\in H^{int}_1$ so that $a\notin H_2$. But the half-plane $H_2\in\RCT$, because its boundary, the straight line $\ell$, is parallel to one of the coordinate axis. In other words, $H_2$ is an (unbounded) closed rectangle. Therefore $\HR[S]\subset H_2$, see definition \rf{HRS}. Because $H_2$ is closed, we have $\HR[S]^{\cl}\subset H_2$ so that, thanks to the lemma's hypothesis, $a\in\HR[S]^{\cl}\subset H_2$, a contradiction.
\par This contradiction implies \rf{AB-V} completing the proof of the lemma.\bx
\smsk
\par Clearly, this lemma implies the statement
($\bigstar 1$) of Proposition \reff{LSEL-F}.
\smsk
\par Let us prove the statement ($\bigstar 2$) which is equivalent to the inequality
\bel{LIP-FL}
\|f(x)-f(y)\|\le (2\lambda+\tlm)\,\rho(x,y),~~~~~x,y\in\Mc.
\ee
\par The proof of this inequality relies on a number of auxiliary results which we present in the next section.

\SECT{4. Proof of the key theorem: the final step.}{4}
\addtocontents{toc}{4. Proof of the key theorem: the final step.\hfill \thepage\par\VST}

\indent

\begin{lemma}\lbl{AB-PR} Let $A,B\subset \R^{2}$ be non-empty convex closed sets such that $A\subset B$, and let $a\in \HR[A]^{\cl}$.
\par Then $\Prm(a,A)$ and $\Prm(a,B)$ are singletons having the following properties:
\msk
\par (i) $\Prm(a,B)\in [\Prm(a,A),a]$;
\smsk
\par  (ii) the following equality
$$
\|\Prm(a,A)-\Prm(a,B)\|=\dist(a,A)-\dist(a,B)
$$
holds.
\end{lemma}
\par {\it Proof.} For the proof of the special case of the lemma for $a\in\HR[S]$ we refer the reader to \cite[p. 301]{S-2002}. See also \cite[p. 67]{S-2020-L}.
\smsk
\par First, we note that if $a\in B$, the statement of the lemma is immediate from Lemma \reff{WMP-S}.
\par Suppose that $a\notin B$. In this case, Lemma \reff{WMP-S} tells us that $\Prj(a;A)$ is one of the vertices of the square $Q(a,r)$ with $r=\dist(a,A)>0$. Because $A\subset B$, the point $a\in \HR[B]^{\cl}$ so that, thanks to Lemma \reff{WMP-S}, $\Prj(a;B)$ is a vertex of the square $Q(a,\alpha)$ where $\alpha=\dist(a,B)>0$.
\par Using a suitable shift and dilation, without loss of generality, we can assume that
$$
a=(0,0),~~r=\dist(a,A)=1,~~~\text{and}~~~ \Prj(a;A)=(1,1).
$$
Clearly, in this case $0<\alpha\le 1$. Furthermore, in these settings the statement of the lemma is equivalent to the property
\bel{P-AL}
\Prj(a;B)=(\alpha,\alpha).
\ee
\par Suppose that this property does not hold, i.e.,
$\Prj(a;B)\in\{(\alpha,-\alpha),(-\alpha,\alpha),
(-\alpha,-\alpha)\}$.
\par In order to get a contradiction, we construct a straight line $\ell_{A}$ which passes through $(1,1)$ and separates (not strictly) the square $Q(a,r)=[-1,1]^{2}$ and $A$. This line determines two closed half-planes, $S^{+}_{A}$ and $S^{-}_{A}$, with the common boundary (i.e., the line $\ell_A$) such that $\RT=S^{+}_{A}\cup S^{-}_{A}$. One of them, say $S^{+}_{A}$, contains $A$, so that $S^{-}_{A}\supset Q(a,r)$. We know that $S^{+}_{A}$ contains $(1,1)$ and does not contain intrinsic points of the square $[-1,1]^{2}$, so that $Q(a,r)\cap\ell_A=(1,1)$. Therefore, the half-plane $S^{+}_{A}$ admits the following representation:
\bel{SA-R2}
S^{+}_{A}=\{x=(x_1,x_2)\in \R^{2}:
(x_1-1)\,h_1+(x_2-1)\,h_2\ge 0\}
\ee
where $h_{1},h_{2}>0$ are certain numbers.
\smsk
\par Let us assume that $\Prj(a;B)=(-\alpha,\alpha)$ and show that this assumption leads to a contradiction. We let $\ell_B$ denote a straight line which passes through the point $(-\alpha,\alpha)$ and separates the square $Q(a,\dist(a,B))=[-\alpha,\alpha]^2$ and the set $B$. Let $S^{+}_{B}$ be that of the two half-planes determined by $\ell_B$ which contains $B$. Then another half-plane, $S^{-}_{B}$, contains $Q(a,\dist(a,B))$, and $S^{+}_{B}\cap S^{-}_{B}=\ell_B$.
\smsk
\par We know that $S^{+}_{B}$ contains the point $(-\alpha,\alpha)$ on its boundary and does not contain intrinsic points of the square $[-\alpha,\alpha]^2$. Therefore, this half-plane can be represented in the form
\bel{SB-N}
S^{+}_{B}=\{(x_1,x_2)\in \R^{2}:
-(x_1+\alpha)\,s_{1}+(x_2-\alpha)\,s_2\ge 0\}
\ee
with certain $s_1,s_2>0$.
\par Thus, $A\subset S^{+}_{A}$ and $A\subset B\subset S^{+}_{B}$, so that $A\subset S^{+}_{A}\cap S^{+}_{B}$ proving that for every $x=(x_1,x_2)\in A$ we have
\bel{S-ABA}
(x_1-1)\,h_1+(x_2-1)\,h_2\ge 0~~~\text{and}~~~
-(x_1+\alpha)\,s_{1}+(x_2-\alpha)\,s_2\ge 0.
\ee
See \rf{SA-R2} and \rf{SB-N}. Note also that since  $S^{+}_{A}\cap S^{+}_{B}\supset A\ne\emp$, we have $h_2+s_2>0$.
\par Let us prove that inequalities \rf{S-ABA} imply the following inclusion:
\bel{X2-A}
A\subset \Hc_\alpha=\{x=(x_1,x_2)\in\RT:x_2\ge\alpha\}.
\ee
\par Indeed, it is easy to see that from \rf{S-ABA} we have
$$
x_2-\alpha\ge \frac{s_1((1+\alpha)h_1+(1-\alpha)h_2))}
{s_1h_2+s_2h_1}\ge 0,~~~~~~x=(x_1,x_2)\in A,
$$
proving \rf{X2-A}.
\par Let us note that the half-plane $\Hc_\alpha$ is a rectangle, i.e., $\Hc_\alpha\in\RCT$. Therefore, the rectangle hull $\Hc[A]\subset\Hc_\alpha$. Furthermore, because $\Hc_\alpha$ is closed, we have $\Hc[A]^{\cl}\subset\Hc_\alpha$ so that, thanks to the lemma's assumption $a=(0,0)\in \Hc_\alpha$. But $\alpha>0$
so that $a=(0,0)\notin \Hc_\alpha$, a contradiction.
\par In a similar way we get a contradiction provided
$\Prj(a;B)=(\alpha,-\alpha)$ or $\Prj(a;B)=(-\alpha,-\alpha)$ proving the required property \rf{P-AL} and the lemma.\bx

\begin{lemma} (i) Let $u\in\Mc$, and let $a\in\HR[F^{[1]}[u:\tlm]]^{\cl}$. Then
\bel{DS-F1}
\dist(a,F^{[1]}[u:\tlm])=\sup_{z\in\Mc}\,\dist(a,F(z)+
\tlm\,\rho(u,z)Q_0)
=\sup_{z\in\Mc}\,[\dist(a,F(z))-\tlm\,\rho(u,z)]_+\,;
\ee
\par (ii) Let $u,v\in\Mc$, and let $a\in\HR[F^{[1]}[u:\tlm]]^{\cl}$ and $b\in\HR[F^{[1]}[v:\tlm]]^{\cl}$.
Then
\bel{DR-X}
|\dist(a,F^{[1]}[u:\tlm])-\dist(b,F^{[1]}[v:\tlm])|
\le \|a-b\|+\tlm\,\rho(u,v).
\ee
\end{lemma}
\par {\it Proof.} {\it (i)} Let $A=F^{[1]}[u:\tlm]$ and, given $z\in\Mc$, let $A_z=F(z)+\tlm\,\rho(u,z)Q_0$. Then, thanks to \rf{WF1-3L}, $A=\cap\{A_z:z\in\Mc\}$. Our goal is to prove that
\bel{O-G}
\dist(a,A)=\sup\{\dist(a,A_z):z\in\Mc\}~~~\text{provided}~~~
a\in\HR[A]^{\cl}.
\ee
\par Lemma \reff{WF1-NE} tells us that $A$ is a non-empty convex and closed subset of $\RT$. Because $A\subset A_z$ for each $z\in\Mc$, the left hand side of the above equality majorizes its right hand side.
\smsk
\par Prove the converse inequality. If $a\in A$, nothing to prove.
\par Let $a\notin A$, and let $\ve\in(0,\dist(a,A))$ be arbitrary. We know that $a\in\HR[A]^{\cl}$ so that, thanks to Lemma \reff{AB-PR}, $\Prm(a,A)$ is a singleton.
\par We let $a_{\ve}$ denote a point on the interval $(\Prm(a,A),a]$ such that $\|a_{\ve}-\Prm(a,A)\|<\ve$. Because $a_\ve\notin A$ and $A=\cap\{A_z:z\in\Mc\}$, there exists an element $\tz\in\Mc$ such that $a_\ve\notin A_{\tz}$. Note that $A\subset A_{\tz}$. Lemma \reff{AB-PR} tells us that in this case $\Prm(a,A_{\tz})$ is a  singleton such that $\Prm(a,A_{\tz})\in[\Prm(a,A),a]$.
\par Then $\Prm(a,A_{\tz})\in [\Prm(a,A),a_\ve]$; otherwise
$a_\ve\in[\Prm(a,A),\Prm(a,A_{\tz})]\subset A_{\tz}$, a contradiction.
\par This proves that $\|\Prm(a,A)-\Prm(a,A_{\tz})\|<\ve$.
Hence,
\be
\dist(a,A)&=&\|a-\Prm(a,A)\|\le \|a-\Prm(a,A_{\tz})\|+
\|\Prm(a,A_{\tz})-\Prm(a,A)\|\nn\\
&\le&
\dist(a,A_{\tz})+\|a_\ve-\Prm(a,A)\|\le \dist(a,A_{\tz})+\ve.
\nn
\ee
\par Since $\ve>0$ can be chosen as small as desired, this implies the required inequality \rf{O-G} proving part {\it (i)} of the lemma.
\msk
\par {\it (ii)} Let $A=F^{[1]}[u:\tlm]$ and $B=F^{[1]}[v:\tlm]$.
Then, thanks to \rf{DS-F1},
\be
|\dist(a,A)-\dist(a,B)|&=&
|\sup_{z\in\Mc}\,[\dist(a,F(z))-\tlm\rho(u,z)]_+
-\sup_{z\in\Mc}\,[\dist(a,F(z))-\tlm\rho(v,z)]_+|
\nn\\
&\le&
\sup_{z\in\Mc}|\,[\dist(a,F(z))-\tlm\rho(u,z)]_+
-[\dist(a,F(z))-\tlm\rho(v,z)]_+|
\nn\\
&\le&
\tlm\,\sup_{z\in\Mc}|\rho(u,z)-\rho(v,z)|
\nn
\ee
so that, thanks to the triangle inequality,
\bel{D-AB1}
|\dist(a,A)-\dist(a,B)|\le\tlm\,\rho(u,v).
\ee
\par Next,
$$
|\dist(a,A)-\dist(b,B)|\le
|\dist(a,A)-\dist(a,B)|+|\dist(a,B)-\dist(b,B)|.
$$
Because $\dist(\cdot,B)$ is a Lipschitz function, from this and \rf{D-AB1}, we have \rf{DR-X} completing the proof of the lemma.\bx
\msk
\par Let $\delta\ge 0$, and let
\bel{H12-DF}
H_1~~\text{and}~~~H_2~~~\text{be two half-planes with}~~~ \dist(H_1,H_2)\le\delta.
\ee
\par Let $\ell_i=\partial H_i$ be the boundary of the half-plane $H_i$, $i=1,2$. Let us represent the half-planes $H_i$, $i=1,2$, in the form
\bel{H12-R}
H_i=\{u\in\RT:\ip{\lh_i,u}\le \alpha_i\}~~~\text{where}~~~ \lh_i\in \SO~~~\text{and}~~~\alpha_i\in\R.
\ee
\par Thus the vector
\bel{H-OL}
\lh_i~~\text{is directed outside of}~~H_i~~\text{and}~~~\lh_i\perp \ell_i,~~~i=1,2.
\ee

\begin{proposition}\lbl{TWO-HP} Let $a_1$ and $a_2$ be two points in $\RT$ such that
\bel{A-I}
a_1\in\HR[H_1\cap(H_2+\delta Q_0)]~~~\text{and}~~~ a_2\in\HR[H_2\cap(H_1+\delta Q_0)].
\ee
Suppose that
\bel{A-II}
\Prm(a_1,H_1)\in H_2+\delta Q_0~~~\text{and}~~~
\Prm(a_2,H_2)\in H_1+\delta Q_0.
\ee
\par Then the following inequality
\bel{PRA-D}
\|\Prm(a_1,H_1)-\Prm(a_2,H_2)\|\le 2\|a_1-a_2\|+\delta
\ee
holds.
\end{proposition}
\par {\it Proof.} We will need a number of auxiliary lemmas. Let us formulate the first of them. Let
\bel{S12-D}
S_1=H_1\cap(H_2+\delta Q_0)~~~~\text{and}~~~~
S_2=H_2\cap(H_1+\delta Q_0).
\ee
We know that $\dist(H_1,H_2)\le\delta$ so that $S_1\ne\emp$ and $S_2\ne\emp$. Furthermore, thanks to \rf{A-I} and \rf{S12-D},
\bel{A12-S}
a_1\in\HR[S_1]~~~~\text{and}~~~a_2\in\HR[S_2].
\ee
\begin{lemma}\lbl{L1-P} Both $\Prm(a_1,H_1)$ and $\Prm(a_2,H_2)$ are singletons. Furthermore,
$\Prm(a_i,H_i)=\Prm(a_i,S_i)$ for every $i=1,2$, and the following inequality
$$
|\dist(a_1,H_1)-\dist(a_2,H_2)|\le
\delta+\|a_1-a_2\|
$$
holds.
\end{lemma}
\par {\it Proof.} Thanks to \rf{A12-S}, the point $a_i\in \HR[S_i]$, so that  $a_i\in\HR[H_i]$ because $S_i\subset H_i$, $i=1,2$, see \rf{S12-D}. Therefore, thanks to Lemma \reff{WMP-S}, $\Prm(a_i,H_i)$ is a singleton for every $i=1,2$.
\par Furthermore, thanks to \rf{A-II}, $\Prm(a_i,H_i)\in S_i$. But $S_i\subset H_i$ so that
$\Prm(a_i,H_i)=\Prm(a_i,S_i)$, $i=1,2$. In particular, $\dist(a_i,H_i)=\dist(a_i,S_i)$, $i=1,2$.
\par Clearly,
$$
\dhf(S_1,S_2)=\dhf(H_1\cap[H_2+\delta Q_0],
H_2\cap[H_1+\delta Q_0])\le \delta.
$$
See \rf{HD-DF}. Therefore,
$$
|\dist(a_1,H_1)-\dist(a_1,H_2)|=
|\dist(a_1,S_1)-\dist(a_1,S_2)|
\le\dhf(S_1,S_2)\le \delta.
$$
\par Note also that the function $\dist(\cdot,H_2)$ is Lipschitz. Hence,
$$
|\dist(a_1,H_1)-\dist(a_2,H_2)|\le |\dist(a_1,H_1)-\dist(a_1,H_2)|+
|\dist(a_1,H_2)-\dist(a_2,H_2)|\le
\delta+\|a_1-a_2\|
$$
proving the lemma.\bx

\begin{lemma}\lbl{L2-P} Inequality \rf{PRA-D} holds provided either $a_1\in H_1$ or $a_2\in H_2$.
\end{lemma}
\par {\it Proof.} For example, suppose that $a_2\in H_2$. Then $\Prm(a_2,H_2)=a_2$ and $\dist(a_2,H_2)=0$. Therefore,  thanks to Lemma \reff{L1-P}, $\dist(a_1,H_1)\le \delta+\|a_1-a_2\|$. Hence,
\be
\|\Prm(a_1,H_1)-\Prm(a_2,H_2)\|&=&
\|\Prm(a_1,H_1)-a_2\|\le \|\Prm(a_1,H_1)-a_1\|+\|a_1-a_2\|
\nn\\
&=&
\dist(a_1,H_1)+\|a_1-a_2\|\le \delta+2\|a_1-a_2\|
\nn
\ee
proving the lemma.\bx
\msk
\par Everywhere below, in the proof of inequality \rf{PRA-D}, we will assume that
\bel{A12-HN}
a_1\notin H_1~~~\text{and}~~~a_2\notin H_2.
\ee
\par Recall that $\ell_i$ is the boundary of the half-plane $H_i$, $i=1,2$. Let us see that the assumption $a_i\notin H_i$, $i=1,2$, implies the following property:
\bel{L-NP}
\ell_i\nparallel Ox_j~~~\text{for every}~~i,j=1,2.
\ee
\par Indeed, suppose that this statement is not true, say for $i=1$, i.e., either $\ell_1\parallel Ox_1$ or $\ell_1\parallel Ox_2$. Then $\HR[H_1]=H_1$. But
$$
a_1\in\HR[H_1\cap(H_2+\delta Q_0)]\subset \HR[H_1]
$$
so that $a_1\in H_1$ which contradicts our assumption that $a_1\notin H_1$.
\smsk
\begin{remark} {\em Our next result, Lemma \reff
{L2-P}, deals with points $a_i$ and half-planes $H_i$, $i=1,2$ such that the vectors
\bel{H-CDR}
\Prm(a_1,H_1)-a_1~~~\text{and}~~~ \Prm(a_2,H_2)-a_2
~~~\text{are co-directed.}
\ee
Recall that this property means that 
$$
\Prm(a_2,H_2)-a_2=\beta\,(\Prm(a_1,H_1)-a_1)~~~\text{for some}~~~\beta>0.
$$
\par We also recall the representation of the half-planes $H_1$ and $H_2$ in the form \rf{H12-R}, i.e., the formulae
$$
H_i=\{u\in\RT:\ip{\lh_i,u}\le \alpha_i\},~~~i=1,2,
$$
where each $\lh_i\in \SO$ is a unit vector and $\alpha_i\in\R$.
\par  Because $\lh_i\perp \ell_i\,(=\partial H_i)$, from \rf{L-NP} we have
$$
\lh_i\nparallel Ox_j~~~\text{for every}~~i,j=1,2.
$$
In particular each $\lh_i$, $i=1,2$, has non-zero coordinates.
\par Finally, let us note the following useful property of metric projections in the space $\LTI=(\RT,\|\cdot\|)$. Let $\alpha\in \R$ and let  $\lh=(h_1,h_2)\in\SO$, $h_1,h_2\ne 0$. Let
$$
H=\{u\in\RT:\ip{\lh,u}\le \alpha\},~~~\text{and let}~~~ a\notin H.
$$
Clearly, in this case $\Prm(a,H)$ is a singleton, and $\Prm(a,H)\ne a$. Then the vector
\bel{A-CH}
a-\Prm(a,H)~~\text{and the vector}~~(\sign h_1,\sign h_2)
~~\text{are co-directed.\rbx}
\ee
}
\end{remark}

\begin{lemma}\lbl{L3-P} Inequality \rf{PRA-D} holds provided condition \rf{H-CDR} is satisfied.
\end{lemma}
\par {\it Proof.} Thanks to the triangle inequality,
\bel{LK}
\|\Prm(a_1,H_1)-\Prm(a_2,H_2)\|\le \|(\Prm(a_1,H_1)-a_1)-(\Prm(a_2,H_2)-a_2)\|+\|a_1-a_2\|.
\ee
Because the vectors $\Prm(a_1,H_1)-a_1$ and $\Prm(a_2,H_2)-a_2$ are co-directed,
\be
\|(\Prm(a_1,H_1)-a_1)-(\Prm(a_2,H_2)-a_2)\|&=&
|\,\|(\Prm(a_1,H_1)-a_1)\|-\|(\Prm(a_2,H_2)-a_2)\|\,|
\nn\\
&=&
|\dist(a_1,H_1)-\dist(a_2,H_2)|.\nn
\ee
Therefore, thanks to Lemma \reff{L1-P},
$$
\|(\Prm(a_1,H_1)-a_1)-(\Prm(a_2,H_2)-a_2)\|\le \delta+\|a_1-a_2\|.
$$
Combining this inequality with \rf{LK}, we obtain the required inequality \rf{PRA-D}.\bx
\msk
\begin{lemma}\lbl{L-PH2} (i) Inequality \rf{PRA-D} holds provided
\bel{SD-12}
\dist(a_1,H_1)+\dist(a_1,H_2)\le\delta.
\ee
(ii) Inequality \rf{PRA-D} holds if $a_1\in H_2$.
\end{lemma}
\par {\it Proof.} {\it (i)} First, let us prove that
\bel{PR-H2}
\|\Prm(a_1,H_2)-\Prm(a_2,H_2)\|\le 2\|a_1-a_2\|.
\ee
\par Indeed, thanks to \rf{L-NP}, $\ell_2\nparallel Ox_1$ and $\ell_2\nparallel Ox_2$ so that $\HR[H_2]=\RT$.  Hence, $a_1,a_2\in\HR[H_2]$.
\par We know that $a_2\notin H_2$. If $a_1\in H_2$, then all conditions of Lemma \reff{L2-P} are satisfied provided $H_1=H_2$ and $\delta=0$. This lemma tells us that in this case inequality \rf{PR-H2} holds.
\par Now, suppose that $a_1\notin H_2$. Then the vectors $\Prm(a_1,H_1)-a_1$ and $\Prm(a_2,H_2)-a_2$ are {\it co-directed}. Therefore, all conditions of Lemma \reff{L3-P} are satisfied for the same case, i.e., for $H_1=H_2$ and $\delta=0$. This lemma tells us that, in these settings, inequality \rf{PR-H2} holds.
\par Thus, we have proved \rf{PR-H2} for every $a_1$ and $a_2$ satisfying \rf{A-I} and \rf{A12-HN}.
\smsk
\par From \rf{PR-H2} and the triangle inequality, we have
\be
\|\Prm(a_1,H_1)-\Prm(a_2,H_2)\|
&\le&
\|(\Prm(a_1,H_1)-a_1)-(\Prm(a_1,H_2)-a_1)\|+
\|\Prm(a_1,H_2)-\Prm(a_2,H_2)\|
\nn\\
&\le&
\|\Prm(a_1,H_1)-a_1\|+\|\Prm(a_1,H_2)-a_1\|+2\|a_1-a_2\|
\nn\\
&=&
\dist(a_1,H_1)+\dist(a_1,H_2)+2\|a_1-a_2\|.
\nn
\ee
\par Combining this inequality with \rf{SD-12}, we obtain
inequality \rf{PRA-D} proving part {\it (i)} of the lemma.
\smsk
\par {\it (ii)} Prove that if $a_1\in H_2$, then
$$
\|\Prm(a_1,H_1)-\Prm(a_1,H_2)\|\le \delta.
$$
\par Indeed, this inequality is immediate from  Lemma \reff{L2-P} applied to the case $a_1=a_2$. Now, from this and \rf{PR-H2}, we have
$$
\|\Prm(a_1,H_1)-\Prm(a_2,H_2)\|\le
\|\Prm(a_1,H_1)-\Prm(a_1,H_2)\|+
\|\Prm(a_1,H_2)-\Prm(a_2,H_2)\|\le
\delta+2\|a_1-a_2\|
$$
proving \rf{PRA-D} and the lemma.\bx
\begin{lemma}\lbl{L4-P} Inequality \rf{PRA-D} holds provided $\ell_1\parallel\ell_2$.
\end{lemma}
\par {\it Proof.} Because  $\ell_1\parallel\ell_2$ and $\ell_i\nparallel Ox_j$, $i,j=1,2$ (see \rf{L-NP}), we have
$$
\HR[H_1\cap(H_2+\delta Q_0)]=\HR[H_2\cap(H_1+\delta Q_0)]=\RT
$$
\par Let us note that, since $\ell_1\parallel \ell_2$ and $\lh_i\perp \ell_i$ (see \rf{H-OL}), the vectors $\lh_1$ and $\lh_2$ are collinear unit vectors. Therefore, either $\lh_1=\lh_2$ or $\lh_1=-\lh_2$.
\par If $\lh_1=\lh_2$, then $H_2$ is a shift of $H_1$, i.e., $H_2=H_1+p$ for some $p\in\RT$.
\par Thanks to \rf{A-CH}, in this case the vectors $\Prm(a_1,H_1)-a_1$ and $\Prm(a_2,H_2)-a_2$ are {\it co-directed}. See Fig. 9-1. Therefore, thanks to Lemma \reff{L3-P}, inequality \rf{PRA-D} holds.
\smsk
\par Let us prove \rf{PRA-D} for $\lh_2=-\lh_1$.
Part {\it (ii)} of Lemma \reff{L-PH2} tells us that
\rf{PRA-D} holds provided $a_1\in H_2$.
\par Now, let us suppose that $a_1\notin H_2$ and prove that \rf{PRA-D} holds as well. In this case $a_1\notin H_1\cup H_2$ (because  $a_1\notin H_1$, see \rf{A12-HN}) so that $H_1\cap H_2=\emp$ as it shown on Fig. 9-2.
\smsk
\begin{figure}[H]
\hspace{15mm}
\includegraphics[scale=0.22]{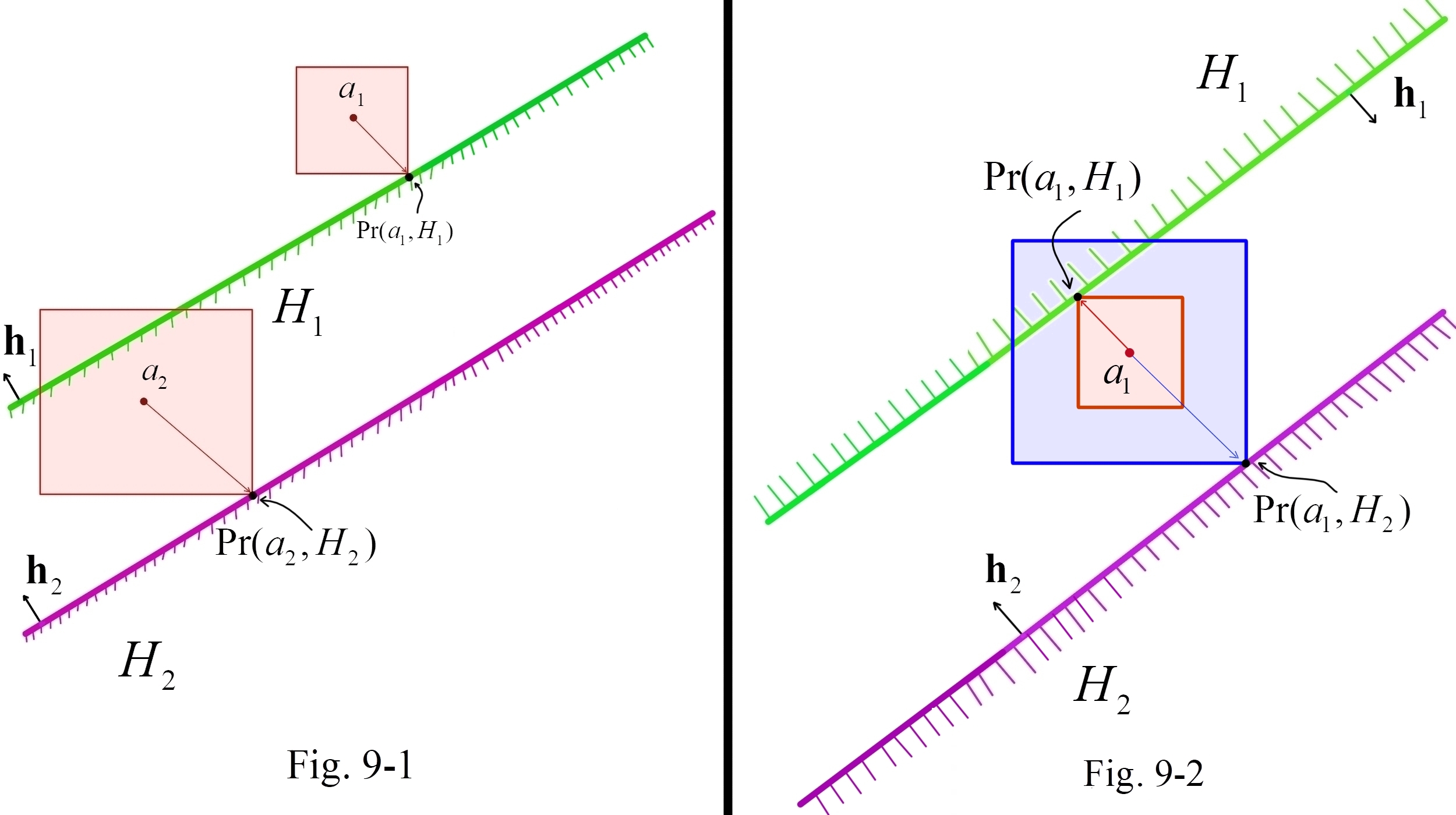}
\caption{Metric projections onto the half-planes $H_1$ and $H_2$ with the parallel boundaries.}
\end{figure}

\par Let us prove that in this case inequality \rf{SD-12} holds. Let $T$ be the closure of the set $\RT\setminus (H_1\cup H_2)$. Clearly, $T$ is the strip between the half-planes $H_1$ and $H_2$. Clearly, $\partial T=\ell\cup\ell_2$. Recall that $\ell_1\parallel \ell_2$ and $\dist(H_1,H_2)\le \delta$ so that
$$
\dist(x,H_2)\le\delta~~\text{for}~~x\in\ell_1
~~~\text{and}~~~
\dist(x,H_1)\le\delta~~\text{for}~~x\in\ell_2.
$$
\par We define a function $f$ on $T$ by letting
$f(x)=\dist(x,H_1)+\dist(x,H_2)$. Clearly, $f$ is a convex continuous function on $T$. Therefore, $\sup_T f=\sup_{\partial T} f$. But
$$
f(x)=\dist(x,\ell_2)\le\delta~~\text{on}~~\ell_1
~~~\text{and}~~~
f(x)=\dist(x,\ell_1)\le\delta~~\text{on}~~\ell_2
$$
so that $\sup_{\partial T} f\le\delta$. Hence, $\sup_{T} f\le\delta$ proving \rf{SD-12}. Therefore, thanks to part {\it (i)} of Lemma \reff{L-PH2}, inequality \rf{PRA-D} holds, and the proof of Lemma \reff{L4-P} is complete.\bx
\msk
\par Thanks to Lemma \reff{L-PH2}, part {\it (ii)}, and Lemma \reff{L4-P}, it remains to prove that inequality  \rf{PRA-D} holds provided $\ell_1\nparallel \ell_2$,
\bel{LPR-A}
a_1~~\text{and}~~a_2~~\text{satisfy \rf{A-I}},~~~
a_1\notin (H_1\cup H_2)~~~\text{and}~~~a_2\notin H_2.
\ee
\par Clearly, without loss of generality, we may assume that $\ell_1\cap\ell_2=\{0\}$. Then 
$$
H_i=\{u\in\RT:\ip{\lh_i,u}\le 0\}, ~~~~i=1,2,
$$
where $\lh_i\in\SO$, $i=1,2$, are {\it non-collinear} vectors. In these settings,
\bel{E-L}
\ell_i=\{u\in\RT: \ip{\lh_i,u}=0\}, ~~~i=1,2.
\ee
\par Let
\bel{SI-H}
\lh_i=(\cos\vf_i,\sin\vf_i)~~~\text{where the angle}~~~ \vf_i\in[0,2\pi),~~~i=1,2.
\ee
\par Because the uniform norm on the plane is invariant under reflections with respect to the coordinate axes and with respect to the bisectors of the coordinate angles, we can also assume that the angles $\vf_1$ and $\vf_2$ satisfy the following conditions:
$$
\vf_1\in(\pi/2,\pi)~~~\text{and}~~~
\vf_2\in(\vf_1,\vf_1+\pi).
$$
\par We know that $\lh_i\perp \ell_i$, $i=1,2$ so that $\lh_1\nparallel \lh_2$ (because  $\ell_1\nparallel \ell_2$). Let us also recall that $\lh_i$ is directed outside of $H_i$, $i=1,2$. We also note that, in the case under consideration, $H_1\cap H_2$ is a convex cone with the vertex at $0$. Moreover, the sets
$$
S_1=H_1\cap(H_2+\delta Q_0)~~~\text{and}~~~S_2=H_2\cap(H_1+\delta Q_0)
$$
are convex cones in $\RT$. Let
\bel{VS-12}
X_1=(s_1,s_2)~~~~\text{and}~~~~X_2=(t_1,t_2)
\ee
be the vertices of the cones $S_1$ and $S_2$ respectively.
\par Thus, $X_1$ is the point of intersection of the line $\ell_1=\partial H_1$ and the line $\tell_2=\partial(H_2+\delta Q_0)$. In turn,
$$
X_2=\ell_2\cap\tell_1~~~\text{where}~~~\ell_2=\partial H_2 ~~~\text{and}~~~\tell_1=\partial(H_1+\delta Q_0).
$$
\par Moreover, thanks to \rf{S12-D}, we have the following representations of the cones $S_1$ and $S_2$:
\bel{S-X12}
S_1=H_1\cap H_2+X_1,~~~ S_2=H_1\cap H_2+X_2.
\ee
\par Let us give explicit formulae for the points $X_1$ and $X_2$. First, we note that
$$
H_i+\delta Q_0=\{u\in\RT:\ip{\lh_i,u}\le \delta\,\|\lh_i\|_1\}, ~~~i=1,2.
$$
Here, given $u=(u_1,u_2)\in\RT$ we let  $$\|u\|_1=|u_1|+|u_2|$$ denote the $\ell^1_2$-norm in $\RT$.
Hence,
\bel{E-LW}
\tell_i=\{u\in\RT: \ip{\lh_i,u}=\delta\, \|\lh_i\|_1\}, ~~~i=1,2.
\ee
\par Let
\bel{A-DL}
A=\left(
\begin{array}{ll}
\cos \vf_1& \sin \vf_1\\
\cos \vf_2& \sin \vf_2
\end{array}
\right)
~~~\text{and let}~~~\Delta=\det A=\sin(\vf_2-\vf_1).
\ee
See \rf{SI-H}. (Clearly, $\Delta\ne 0$ because $\lh_1\nparallel\lh_2$.) We know that $$X_1=(s_1,s_2)=\ell_1\cap \tell_2$$ so that, thanks to \rf{E-L} and \rf{E-LW}, the vector $(s_1,s_2)$ is the solution of the system of linear equations
$$
A\left(
\begin{array}{l}
s_1\\s_2
\end{array}
\right)=\left(
\begin{array}{l}
\hspace*{5mm} 0\\
\delta\,\|\lh_2\|_1
\end{array}
\right).
$$
Therefore,
\bel{S1-S2}
s_1=\frac{1}{\Delta}\,\left|
\begin{array}{ll}
0& \sin \vf_1\\
\delta\,\|\lh_2\|_1& \sin \vf_2
\end{array}
\right|
=-\frac{\delta}{\Delta}\,\|\lh_2\|_1\,\sin \vf_1
~~~\text{and}~~~
s_2=\frac{\delta}{\Delta}\,\|\lh_2\|_1\,\cos \vf_1.
\ee
Thus,
\bel{X-1}
X_1=\frac{\delta}{\Delta}\,\|\lh_2\|_1\,(-\sin \vf_1,\cos \vf_1).
\ee
In the same way we prove that
\bel{X-2}
X_2=\frac{\delta}{\Delta}\,\|\lh_1\|_1\,
(\sin \vf_2,-\cos \vf_2).
\ee
\begin{lemma} Inequality \rf{PRA-D} holds provided $\vf_1\in(\pi/2,\pi)$ and $\vf_2\in(\vf_1,\pi)$.
\end{lemma}
\par {\it Proof.} We recall that in the case under consideration
$$
\lh_i=(\cos\vf_i,\sin\vf_i)\in\SO~~~\text{where}~~~ \vf_i\in(\pi/2,\pi),~~~i=1,2.
$$
Hence, $(\sign(\cos\vf_i),\sign(\sin\vf_i))=(-1,1)$.
\smsk
\par This equality and property \rf{A-CH} imply the following: for every $i=1,2$, the vector $(1,-1)$ and the vector $\Prm(a_i,H)-a_i$ are co-directed.
See Fig. 10.
\msk
\begin{figure}[H]
\hspace{-1mm}
\includegraphics[scale=0.25]{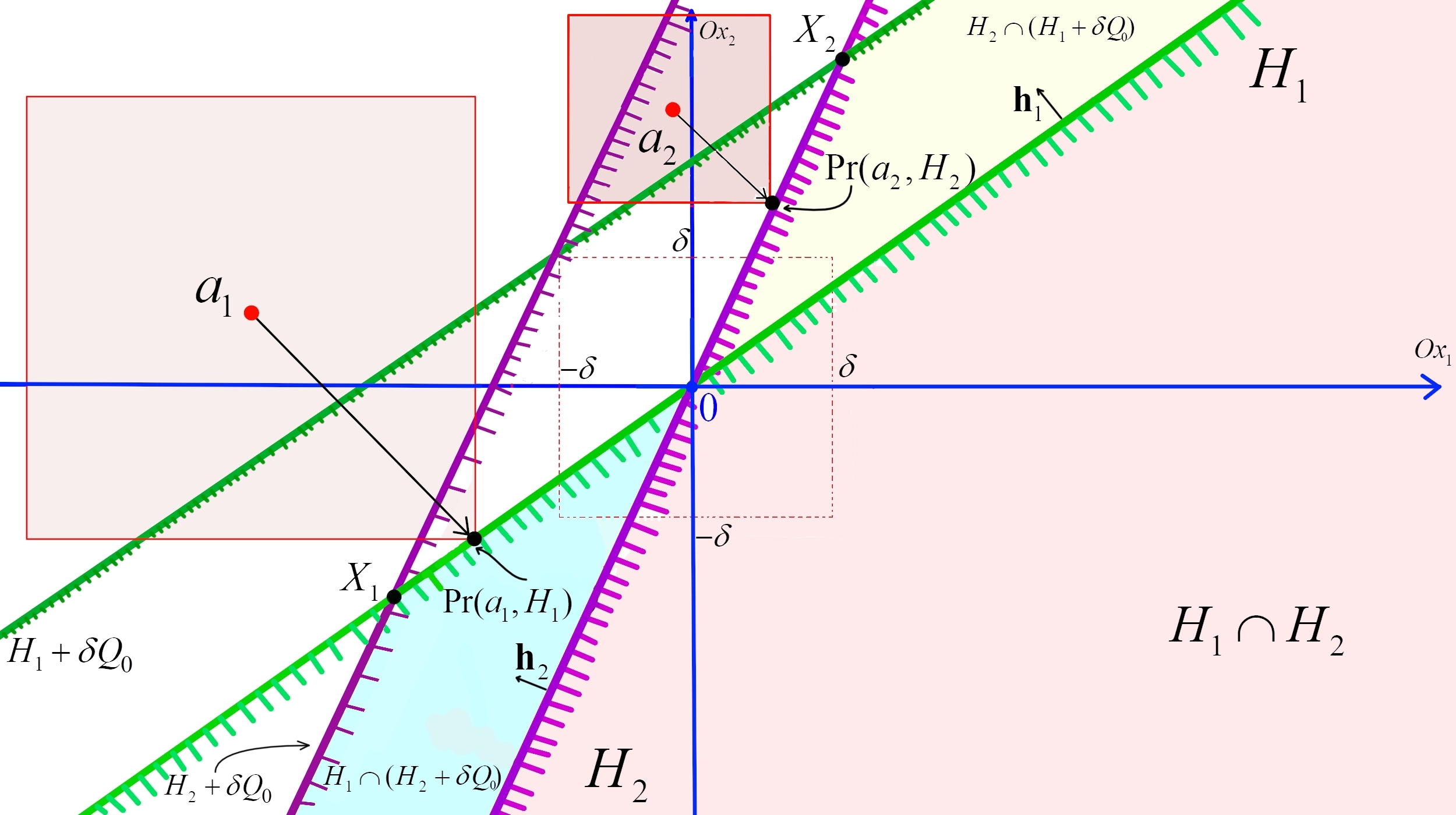}
\hspace*{-10mm}
\caption{The half-planes $H_1$ and $H_2$ with the non-parallel boundaries: the first case.}
\end{figure}
\msk
\par This proves that
$$
\Prm(a_1,H)-a_1~~~\text{and}~~~\Prm(a_2,H)-a_2
~~~\text{are co-directed vectors}.
$$
\par Therefore, thanks to Lemma \reff{L3-P}, \rf{PRA-D} holds, and the proof of the lemma is complete. \bx
\bsk
\begin{lemma}\lbl{L6-P} Inequality \rf{PRA-D} holds provided $\vf_1\in(\pi/2,\pi)$ and $\vf_2\in(\pi,\vf_1+\pi)$.
\end{lemma}
\par {\it Proof.} First, let us prove inequality \rf{PRA-D} in the case $\vf_1\in(\pi/2,\pi)$ and $\vf_2\in\left(\pi,\tfrac{3}{2}\pi\right)$. See Fig. 11.

\smsk
\begin{figure}[H]
\hspace{-7mm}
\includegraphics[scale=0.78]{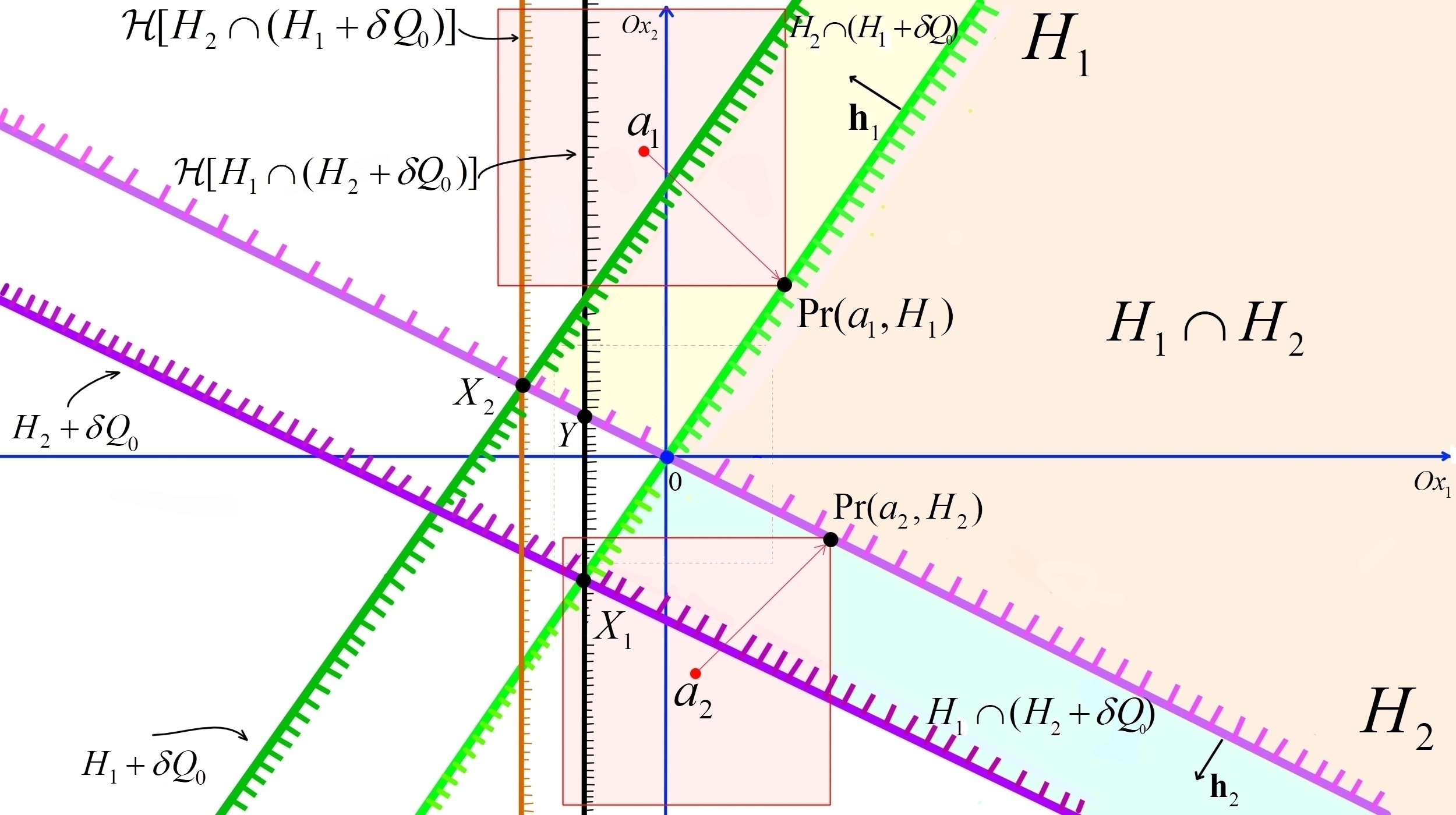}
\hspace*{10mm}
\caption{The half-planes $H_1$ and $H_2$ with the non-parallel boundaries: $\vf_2\in\left(\pi,\tfrac{3}{2}\pi\right)$.}
\end{figure}
\msk
\par In this case, the following inequalities hold:
\bel{CS-2}
\cos\vf_1<0,~\sin\vf_1>0,~
\cos\vf_2<0~~~\text{and}~~~\sin\vf_2<0.
\ee
\par Therefore, thanks to \rf{A-CH}, the vector $\Prm(a_1,H_1)-a_1$ is co-directed with the vector $(1,-1)$, and the vector $\Prm(a_2,H_2)-a_2$ is co-directed with $(1,1)$. Moreover, in this case the convex cone $H_1\cap H_2$ (with the vertex $0$) contains the positive semi axis $Ox_1^+(=\{(t,0):t\ge 0\})$. This implies the following properties of the rectangular hulls of the sets $S_1=H_1\cap (H_2+\delta Q_0)$ and  $S_2=H_2\cap (H_1+\delta Q_0)$:
$$
\HR[S_1]=\{(u_1,u_2)\in\RT:u_1\ge s_1\}~~~\text{and}~~~
\HR[S_2]=\{(u_1,u_2)\in\RT:u_1\ge t_1\}.
$$
\par We recall that $X_1=(s_1,s_2)$ and $X_2=(t_1,t_2)$. See \rf{VS-12} and Fig. 10.
\par Let $\alpha_i\in(0,\pi/2)$, $i=1,2$, be the angle between the straight line $\ell_i=\partial H_i$ and the axis $Ox_1$. On Fig. 10 we consider the case $\alpha_2\le \alpha_1$. We note that no necessity in the additional consideration of the case $\alpha_2>\alpha_1$ because it can be obtained from the case $\alpha_2<\alpha_1$ with the help of a suitable reflections with respect to the coordinate axes and the bisectors of the coordinate angles.
\par Let us prove that if $\alpha_2\le \alpha_1$, then
\bel{SGT}
t_1\le s_1\le 0.
\ee
\par First, let us note that, thanks to \rf{A-DL}, $\Delta=\sin(\vf_2-\vf_1)>0$ (because $\vf_2-\vf_1\in(0,\pi)$.) Also note that, thanks to \rf{CS-2}, $\sin \vf_1>0$. Hence, thanks to formula \rf{X-1}, $s_1\le 0$.
\par Let us see that
\bel{ST-F}
s_1-t_1=\frac{\delta}{\Delta}\sin(\vf_1+\vf_2)=
\frac{\delta}{\Delta}\sin(\alpha_1-\alpha_2).
\ee
\par Indeed, thanks to \rf{CS-2},
$$
\|\lh_1\|_1= -\cos \vf_1+\sin \vf_1~~~\text{and}~~~
\|\lh_2\|_1= -\cos \vf_2-\sin \vf_2,
$$
so that, thanks to formulae \rf{X-1} and \rf{X-2},
\be
s_1-t_1&=&(\delta/\Delta)
\{-\|\lh_2\|_1\sin\vf_1-\|\lh_1\|_1\sin\vf_2\}
\nn\\
&=&
(\delta/\Delta)\{-(-\cos \vf_2-\sin \vf_2)\sin\vf_1-
(-\cos \vf_1+\sin \vf_1)\sin\vf_2\}
\nn\\
&=&
(\delta/\Delta)
\{\cos \vf_2\sin\vf_1+\cos \vf_1\sin\vf_2\}
=(\delta/\Delta)\sin(\vf_1+\vf_2).
\nn
\ee
\par We note that $\alpha_1=\vf_1-\pi/2$ and $\alpha_2=\tfrac32\pi-\vf_2$. Hence, we have
$$
\vf_1+\vf_2=\alpha_1-\alpha_2+2\pi
$$
proving \rf{ST-F}. It remains to note that $\sin(\alpha_1-\alpha_2)>0$ (because $0<\alpha_2\le\alpha_1<\pi/2$) and $\Delta>0$ so that, thanks \rf{ST-F}, $t_1\le s_1$ proving \rf{SGT}.
\par In particular, this inequality implies the inclusion
$\HR[S_1]\subset \HR[S_2]$ as it shown on Fig. 10. Furthermore, thanks to \rf{LPR-A}
\bel{CTP}
a_1\notin H_1\cup H_2~~~\text{and}~~~a_1\in \HR[S_1]=\HR[H_1\cap(H_2+\delta Q_0)].
\ee
\par Let $Y$ be the point of intersection of the line $\ell_2$ (the boundary of $H_2$) and the line $\ellh$ passing through the point $X_1=(s_1,s_2)$ and parallel to the axis $Ox_2$. (Thus, $Y=(s_1,y_2)$ for some $y_2\in\R$). Recall that the point $X_2=(t_1,t_2)$ lies on the line $\ell_2$, and, thanks to inequality \rf{SGT}, $t_1\le s_1\le 0$. In particular, these observations shows that $Y\in[X_2,O]$ where $O=0$ is the origin.
\par Conditions \rf{CTP} shows that the point $a_1$ belongs to the triangle $\tT=\Delta(X_1,Y,O)$ with the vertices at the points $X_1$, $Y$ and $O$. Because
$Y\in[X_2,O]$, the triangle $\tT$ is a subset of the triangle $T=\Delta(X_1,X_2,O)$ with vertices at $X_1$, $X_2$ and $O=0$. Hence, $a_1\in T$.
\par Let
$$
f(x)=\dist(x,H_1)+\dist(x,H_2),~~~x\in T.
$$
\par Prove that $f(x)\le\delta$ on $T$. Indeed, $f$ is a convex continuous function on $T$ so that its maximum is attained on the set of vertices of the triangle $T$, i.e., at the points $O,X_1$ and $X$. But $f(O)=0$, $f(X_1)=\delta$ (because $X_1\in\tell_2=\partial(H_2+\delta\,Q_0)$) and $f(X_2)=\delta$ (because $X_2\in\tell_1=\partial(H_1+\delta\,Q_0)$). This proves the required inequality $f(x)\le\delta$, $x\in T$.
\par In particular,
\bel{LC-FG}
f(a_1)=\dist(a_1,H_1)+\dist(a_1,H_2)\le\delta.
\ee
\par Thus, condition \rf{SD-12} of Lemma \reff{L-PH2} holds. Thanks to this lemma, inequality \rf{PRA-D} holds proving Lemma \reff{L6-P} for the angles $\vf_1\in(\pi/2,\pi)$ and $\vf_2\in\left(\pi,\tfrac{3}{2}\pi\right)$.
\smsk
\par Let us prove Lemma \reff{L6-P} for
\bel{F12-S}
\vf_1\in(\pi/2,\pi)~~~\text{and}~~~ \vf_2\in\left(\frac32\,\pi,\vf_1+\pi\right).~~~~~~
\text{See Fig. 12.}
\ee

\begin{figure}[H]
\hspace*{7mm}
\includegraphics[scale=0.71]{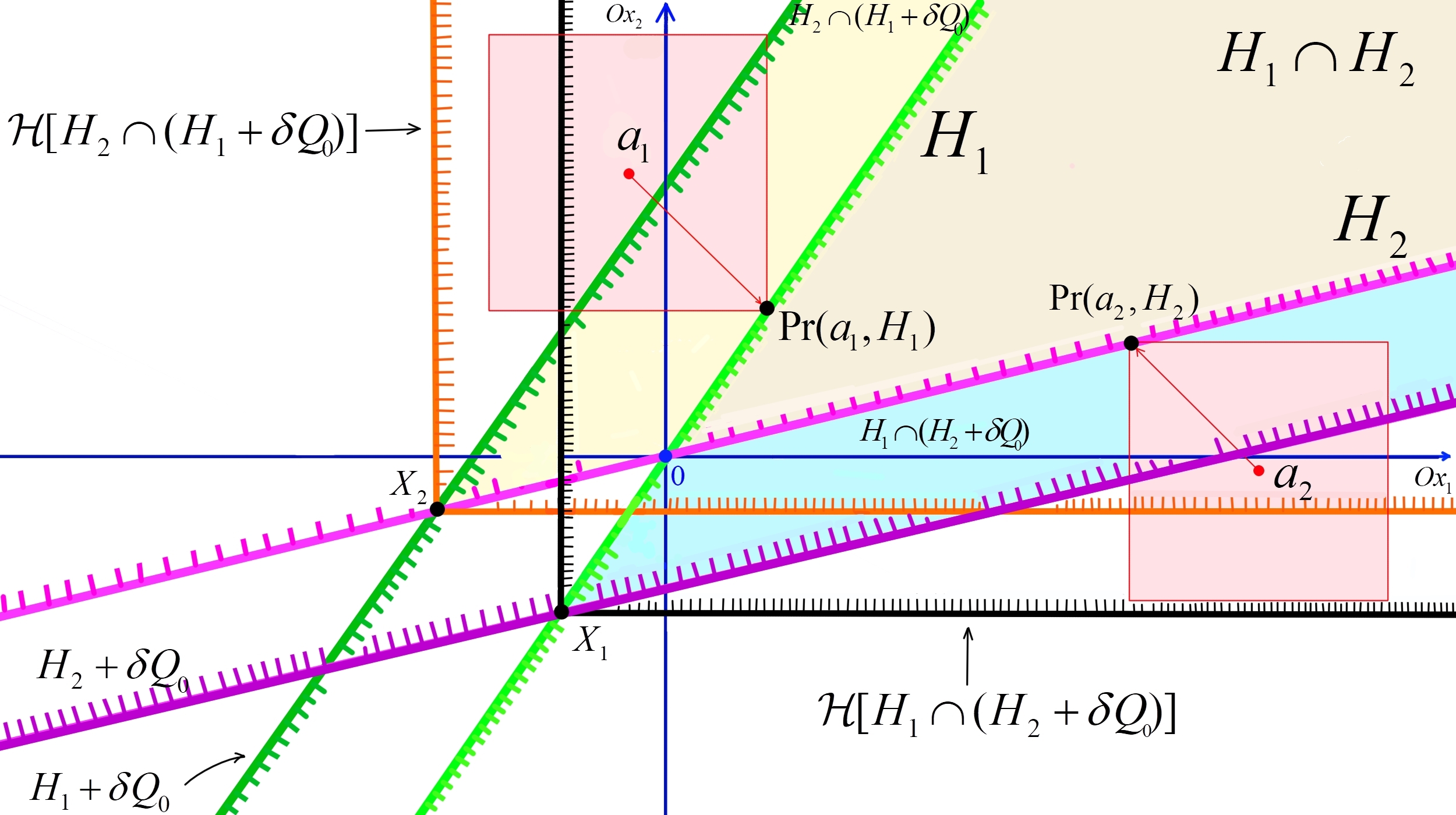}
\hspace*{10mm}
\caption{The half-planes $H_1$ and $H_2$ with the non-parallel boundaries:
$\vf_2\in\left(\frac32\,\pi,\vf_1+\pi\right)$.}
\end{figure}
\par In this case,
\bel{CS-3}
\cos\vf_1<0,~\sin\vf_1>0,~ \cos\vf_2>0~~~\text{and}~~~\sin\vf_2<0.
\ee
Therefore, thanks to \rf{A-CH}, the vector $\Prm(a_1,H_1)-a_1$ is co-directed with $(1,-1)$, and
the vector $\Prm(a_2,H_2)-a_2$ is co-directed with $(-1,1)$. Moreover,
$$
\HR[S_1]=\{(u_1,u_2)\in\RT:u_1\ge s_1,\, u_2\ge s_2\}
~~~\text{and}~~~
\HR[S_2]=\{(u_1,u_2)\in\RT:u_1\ge t_1,\, u_2\ge t_2\}.
$$

\par Let us prove that
\bel{X1-X2}
X_1-X_2=\delta(1,-1)
\ee
$X_1=(s_1,s_2)$ and $X_2=(t_1,t_2)$ are the points defined by \rf{VS-12}. For explicit formulae for $s_i,t_i$, $i=1,2$, see \rf{X-1} and \rf{X-2}.
\smsk
\par Thanks to \rf{CS-3}, we have
$$
\|\lh_1\|_1=|\cos\vf_1|+|\sin\vf_1|=
-\cos\vf_1+\sin\vf_1,~~~
\|\lh_2\|_1=|\cos\vf_2|+|\sin\vf_2|=
\cos\vf_2-\sin\vf_2.
$$
Therefore, thanks to \rf{X-1} and \rf{X-2},
$$
X_1=\frac{\delta}{\Delta}\,(\cos\vf_2-\sin\vf_2)\,(-\sin \vf_1,\cos \vf_1)
~~~\text{and}~~~
X_2=\frac{\delta}{\Delta}\,(-\cos\vf_1+\sin\vf_1)\,
(\sin \vf_2,-\cos \vf_2).
$$
Hence,
\be
X_1-X_2&=&
\frac{\delta}{\Delta}\,((\cos\vf_2-\sin\vf_2)\,(-\sin \vf_1)-(-\cos\vf_1+\sin\vf_1)\sin \vf_2,
\nn\\
&&
(\cos\vf_2-\sin\vf_2)
\cos \vf_1-(-\cos\vf_1+\sin\vf_1)(-\cos \vf_2))
\nn\\
&=&
\frac{\delta}{\Delta}\,(-\cos\vf_2\sin\vf_1
+\cos\vf_1\sin\vf_2,
-\sin\vf_2\cos\vf_1+\sin\vf_1\cos\vf_2)
\nn\\
&=&\frac{\delta}{\Delta}\sin(\vf_2-\vf_1)(1,-1).
\nn
\ee
\par Thanks to \rf{A-DL}, $\Delta=\sin(\vf_2-\vf_1)$, and the proof of \rf{X1-X2} is complete.
\smsk
\par Thanks to this equality, $t_1=s_1-\delta$ and $t_2=s_2+\delta$. Furthermore, \rf{X1-X2} and \rf{S-X12} imply the following:
$$
S_2=S_1+\delta(-1,1)~~~\text{and}~~~
\HR[S_2]=\HR[S_1]+\delta(-1,1).
$$
\par Let us prove that
\bel{ST-CM}
t_1\le s_1\le 0~~~~\text{and}~~~~s_2\le t_2\le 0.
\ee
\par In fact, we know that
$$
\vf_1\in(\pi/2,\pi)~~~~\text{and}~~~~
\vf_2\in(\frac32\,\pi,\vf_1+\pi).
$$
Hence, $0<\vf_2-\vf_1<\pi$ proving that $\Delta=\sin(\vf_2-\vf_1)>0$. We also know that $\sin \vf_1>0$ so that, thanks to \rf{S1-S2}, $s_1\le 0$. In addition, $t_1=s_1-\delta\le s_1$ proving the first inequality in \rf{ST-CM}.
\par Next, thanks to \rf{X-2},
$$
t_2=-(\delta/\Delta)\,\|\lh_1\|_1\,\cos\vf_2.
$$
But $\Delta>0$ and $\cos\vf_2>0$ so that $t_2\le 0$. Moreover, $s_2=t_2-\delta$, so that $s_2\le t_2$, and the proof of \rf{ST-CM} is complete.
\smsk
\par We recall that the point $a_1$ satisfies conditions \rf{CTP}. Therefore,
$$
a_1\in T=\Delta(X_1,X_2,O).
$$
\par See Fig. 12. From this and \rf{LC-FG} it follows that
condition \rf{SD-12} of Lemma \reff{L-PH2} is satisfied. This lemma tells us that inequality \rf{PRA-D} holds proving Lemma \reff{L6-P} for the case \rf{F12-S}.
\smsk
\par The proof of Lemma \reff{L6-P} is complete.\bx
\smsk
\par Finally, the results of Lemmas \reff{L1-P} -- \reff{L6-P} imply the required inequality \rf{PRA-D} completing the proof of Proposition \reff{TWO-HP}.\bx
\msk
\par We are in a position to complete the proof of Proposition \reff{LSEL-F}.
\smsk
\par {\it Proof of inequality \rf{LIP-FL}.}  Let us fix elements $x,y\in\Mc$. We set
$$
A_1=F^{[1]}[x:\tlm],~~~A_2=F^{[1]}[y:\tlm]
$$
and $a_1=g(x)$, $a_2=g(y)$. Recall that $g:\Mc\to\R$ is the mapping satisfying \rf{WG-HF} and \rf{G-LIP}.
\par We also recall that the mapping $F^{[1]}[\cdot:\tlm]$ is defined by \rf{WF1-3L}. Thus,
\bel{A-R12}
A_i=\cap\{A_i^{[u]}: u\in\Mc\},~~i=1,2,
\ee
where given $u\in\Mc$, we set
\bel{AU-12}
A_1^{[u]}=F(u)+\tlm\,\rho(u,x)Q_0
~~~\text{and}~~~
A_2^{[u]}=F(u)+\tlm\,\rho(u,y)Q_0.
\ee
\par Lemma \reff{WF1-NE} tells us that each $A_i$, $i=1,2$, is a non-empty closed convex subset of $\RT$.
\par Thanks to inequality \rf{G-LIP}, we have
\bel{A-12}
\|a_1-a_2\|\le\lambda\,\rho(x,y),
\ee
and, thanks to \rf{WG-HF},
\bel{HA-12}
a_i\in\HR[A_i]^{\cl},~~~i=1,2.
\ee
Furthermore, formula \rf{FX-PR} tells us that
$$
f(x)=\Prm(a_1,A_1)~~~~\text{and}~~~~f(y)=\Prm(a_2,A_2).
$$
(We also recall that, thanks to Lemma \reff{WMP-S},
the metric projection $\Prm(a_i,A_i)$ is well defined, i.e., $\Prm(a_i,A_i)$ is a singleton.)
\par In these settings, the required inequality \rf{LIP-FL}
reads as follows:
\bel{NV-12}
\|\Prm(a_1,A_1)-\Prm(a_2,A_2)\|\le (2\lambda+\tlm)\,\rho(x,y).
\ee
\par Let us note that this inequality is immediate from \rf{A-12} provided $a_i\in A_i$, $i=1,2$, (because in this case $\Prm(a_i,A_i)=a_i$).
\par Suppose that either $a_1\notin A_1$ or $a_2\notin A_2$. Without loss of generality, we may assume that $a_1\notin A_1$. Fix $\ve>0$ and prove that there exists a half-plane $H_1\in\HPL$ such that
\bel{H1-A1}
H_1\supset A_1,~~~~H_1+\tlm\,\rho(x,y)Q_0\supset A_2,
\ee
and
\bel{D-A1}
\|\Prm(a_1,A_1)-\Prm(a_1,H_1)\|<\ve.
\ee
\par We construct the half-plane $H_1$ as follows: Because $a_1\notin A_1$, we have $\Prm(a_1,A_1)\ne a_1$ so that
$(\Prm(a_1,A_1),a_1]$ is a non-empty semi-open interval in $\RT$. Let us pick a point 
$$
a^{(\ve)}\in(\Prm(a_1,A_1),a_1]
$$
such that
\bel{H-EP1}
\|a^{(\ve)}-\Prm(a_1,A_1)\|<\ve.
\ee
Because $\Prm(a_1,A_1)$ is the nearest to $a_1$ point on $A_1$, we have 
$$
(\Prm(a_1,A_1),a_1]\cap A_1=\emp.
$$
Therefore,
$$
a^{(\ve)}\notin A_1=\cap\{A_1^{[u]}: u\in\Mc\}.
$$
\par See \rf{A-R12}. This implies the existence of an element $u\in\Mc$ such that $a^{(\ve)}\notin A_1^{[u]}$. We let $B$ denote the set $A_1^{[u]}$. Thus,
\bel{AE-B}
a^{(\ve)}\notin B=A_1^{[u]}=F(u)+\tlm\,\rho(u,x) Q_0.
\ee
See \rf{AU-12}. Thanks to \rf{HA-12} and \rf{A-R12},
\bel{B-A1}
a_1\in\HR[A_1]^{\cl}~~~~\text{and}~~~~A_1\subset B.
\ee
Therefore, thanks to Lemma \reff{AB-PR}, the metric projections $\Prm(a_1,A_1)$ and $\Prm(a_1,B)$ are singletons such that
$$
\Prm(a_1,B)\in[\Prm(a_1,A_1),a_1].
$$ 
See Fig. 13.

\begin{figure}[H]
\hspace{9mm}
\includegraphics[scale=0.7]{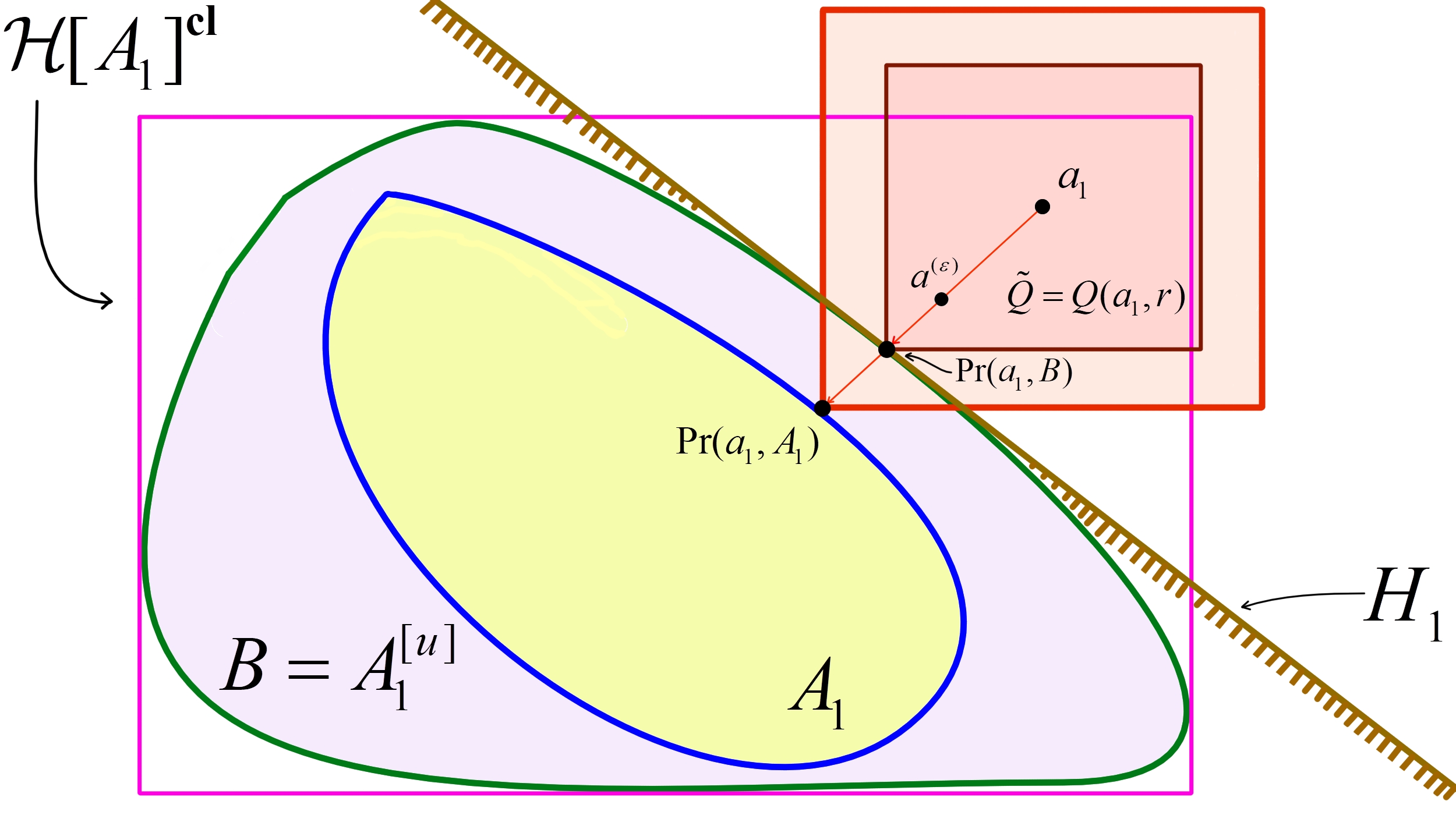}
\hspace*{10mm}
\caption{Metric projections of $a_1$ onto $A_1$ and $B$.}
\end{figure}
\par We note that $\Prm(a_1,B)\in [\Prm(a_1,A_1),a^{(\ve)}]$; indeed, otherwise $a^{(\ve)}\in [\Prm(a_1,A_1),\Prm(a_1,B)]\subset B$, a contradiction. See \rf{AE-B}. Hence, thanks to \rf{H-EP1},
\bel{PR-E1}
\|\Prm(a_1,A_1)-\Prm(a_1,B)\|\le \|a^{(\ve)}-\Prm(a_1,A_1)\|<\ve.
\ee
\par Let $\tQ=Q(a_1,r)$ where $r=\dist(a_1,B)$. Thus, $\tQ\cap B=\{\Prm(a_1,B)\}$. Therefore, thanks to the separation theorem, there exists a half-plane $H_1\in\HPL$ which contains $B$ and separates (not strictly) $\tQ$ and $B$. Thus, $B\subset H_1$  and $\tQ\cap H_1=\Prm(a_1,B)$
as it shown on Fig. 12. In particular, these properties imply the equality $\Prm(a_1,H_1)=\Prm(a_1,B)$.
\smsk
\par Let us see that inclusions \rf{H1-A1} and inequality \rf{D-A1} hold for the half-plane $H_1$.
In fact, \rf{D-A1} is immediate from \rf{PR-E1} and the last equality.
\par Prove \rf{H1-A1}. We know that $A_1\subset B$, see \rf{B-A1}, so that $A_1\subset B\subset H_1$. We also recall that $B=F(u)+\tlm\,\rho(u,x) Q_0$ (see \rf{AE-B}). Therefore,
$$
H_1+\tlm\,\rho(x,y)Q_0\supset B+\tlm\,\rho(x,y)Q_0=
F(u)+\tlm\,\rho(u,x) Q_0+\tlm\,\rho(x,y)Q_0=
F(u)+\tlm\,(\rho(u,x)+\rho(x,y))\,Q_0.
$$
Therefore, thanks to the triangle inequality, \rf{A-R12} and \rf{AU-12}, we have
$$
H_1+\tlm\,\rho(x,y)Q_0\supset
F(u)+\tlm\,\rho(u,y) \,Q_0=A_2^{[u]}\supset A_2
$$
proving \rf{H1-A1}.
\msk
\par Next, let us construct a half-plane $H_2\in\HPL$ having the following properties:
\bel{H2-A2}
H_2\supset A_2,~~~~H_2+\tlm\,\rho(x,y)Q_0\supset A_1,
\ee
and
\bel{D-A2}
\|\Prm(a_2,A_2)-\Prm(a_2,H_2)\|<\ve.
\ee
\par If $a_2\notin A_2$, we define $H_2$ in the same way as we have defined $H_1$ for $a_1$. In this case, properties  \rf{H2-A2}, \rf{D-A2} are a complete analog of properties
\rf{H1-A1} and \rf{D-A1} obtained for the point $a_1$.
\par If $a_2\in A_2$, we set
$$
H_2=H_1+\tlm\,\rho(x,y)Q_0.
$$
\par Clearly, $H_2$ is a half-plane. Let us see that inclusions \rf{H2-A2} and inequality \rf{D-A2} hold for this choice of $H_2$. Indeed, thanks to the second inclusion in \rf{H1-A1}, we have $H_2\supset A_2$. In turn, thanks to the first inclusion,
$$
H_2+\tlm\,\rho(x,y)Q_0=(H_1+\tlm\,\rho(x,y)Q_0)
+\tlm\,\rho(x,y)Q_0\supset H_1\supset A_1,
$$
proving \rf{H2-A2}. Finally, inequality \rf{D-A2} is trivial because $\Prm(a_2,A_2)=\Prm(a_2,H_2)(=a_2)$. (Recall that $a_2\in A_2\subset H_2$.)
\msk
\par Now, we set
\bel{DL}
\delta=\tlm\,\rho(x,y)+\ve.
\ee
\par Let us prove that the points $a_1,a_2$ and the half-planes $H_1$ and $H_2$ satisfy conditions \rf{H12-DF}, \rf{A-I} and \rf{A-II}.
\par Thanks to \rf{H1-A1}, $H_1\supset A_1$, and, thanks to \rf{H2-A2}, $H_2+\tlm\,\rho(x,y)Q_0\supset A_1$. Hence,
$$
H_1\cap(H_2+\tlm\,\rho(x,y)Q_0)\supset A_1.
$$
\par Note that $\tlm\,\rho(x,y)<\delta$; see \rf{DL}. Therefore,
\bel{WW}
H_1\cap(H_2+\delta\,Q_0)\supset A_1.
\ee
\par We also know that $A_1\ne\emp$ so that  $H_1\cap(H_2+\delta\,Q_0)\ne\emp$ as well. This proves that  $\dist(H_1,H_2)\le\delta$ so that the condition \rf{H12-DF} is satisfied.
\smsk
\par Let us prove that the points $a_1$ and $a_2$ satisfy condition \rf{A-I}. Indeed, inclusion \rf{WW} tells us that
$$
\HR[A_1]\subset \HR[H_1\cap(H_2+\delta\,Q_0)].
$$
Clearly, the set $\HR[H_1\cap(H_2+\delta\,Q_0)]$ is {\it closed} (as the rectangle hull of intersection of two half-planes). Hence,
$$
\HR[A_1]^{\cl}\subset \HR[H_1\cap(H_2+\delta\,Q_0)].
$$
But, thanks to \rf{HA-12}, $a_1\in \HR[A_1]^{\cl}$ proving that $a_1\in \HR[H_1\cap(H_2+\delta\,Q_0)]$.
\par In the same fashion we show that $a_2\in \HR[H_2\cap(H_1+\delta\,Q_0)]$ proving that condition \rf{A-I} holds.
\smsk
\par Let us show that condition \rf{A-II} is satisfied as well. Thanks to \rf{H2-A2},
$$
\Prm(a_1,A_1)\in A_1\subset H_2+\tlm\,\rho(x,y)\,Q_0,
$$
and, thanks to \rf{D-A1},
$$
\Prm(a_1,H_1)\in \Prm(a_1,A_1)+\ve\,Q_0\subset
A_1+\ve\,Q_0.
$$
Therefore, thanks to \rf{DL},
$$
\Prm(a_1,H_1)\in A_1+\ve\,Q_0\subset (H_2+\tlm\,\rho(x,y)\,Q_0)+\ve\,Q_0=
H_2+\delta\,Q_0.
$$
\par In the same way we show that $\Prm(a_2,H_2)\in
H_1+\delta\,Q_0$ completing the proof of \rf{A-II}.
\msk
\par Therefore, thanks to Proposition \reff{TWO-HP}, inequality \rf{PRA-D} holds. This inequality together with \rf{A-12} and \rf{DL} imply the following:
$$
\|\Prm(a_1,H_1)-\Prm(a_2,H_2)\|\le 2\|a_1-a_2\|+\delta
\le 2\lambda\,\rho(x,y)+ \tlm\,\rho(x,y)+\ve.
$$
From this, \rf{D-A1} and \rf{D-A2}, we have
\be
\|\Prm(a_1,A_1)-\Prm(a_2,A_2)\|&\le&
\|\Prm(a_1,A_1)-\Prm(a_1,H_1)\|+
\|\Prm(a_1,H_1)-\Prm(a_2,H_2)\|
\nn\\
&+&\|\Prm(a_2,H_2)-\Prm(a_2,A_2)\|
\le (2\lambda+\tlm)\,\rho(x,y)+3\ve.
\nn
\ee
\par Since $\ve>0$ is arbitrary, this implies \rf{NV-12} proving the required inequality \rf{LIP-FL} and completing the proof of Proposition \reff{LSEL-F}.\bx
\smsk
\par Finally, combining part (i) of Proposition \reff{WLS-T} with Proposition \reff{LSEL-F}, we obtain the statement of Theorem \reff{W-CR}.\bx

\SECT{5. Lipschitz selection criteria in the two dimensional case.}{5}

\addtocontents{toc}{5. Lipschitz selection criteria in the two dimensional case.\hfill \thepage\par\VST}

\indent\par {\bf 5.1 Constructive criteria for Lipschitz selections: proofs.}
\addtocontents{toc}{~~~~5.1 Constructive criteria for Lipschitz selections: proofs.\hfill\thepage\par\VST}

\indent
\par We begin with the proof of Theorem \reff{CR-LS1}. This proof is based on the following result.
\begin{proposition}\lbl{WR-LM3} Let $\Mf=\MR$ be a pseudometric space, and let $F:\Mc\to\CRT$ be a set-valued mapping. Given constants $\tlm$ and $\lambda$, $0\le\lambda\le\tlm$, let us assume that
\bel{R-LMW}
\RL[x,x':\tlm]\cap \left\{\RL[y,y':\tlm]+\lambda\,\rho(x,y)\,Q_0\right\}\ne\emp~~~
\text{for every}~~~x,x',y,y'\in\Mc.
\ee
\par Then for every $x,x',x'',y,y',y''\in\Mc$, the following property
$$
\Wc_F[x,x',x'':3\tlm]\cap\left\{\Wc_F[y,y',y'':3\tlm]
+\lambda\,\rho(x,y)Q_0\right\}\ne\emp
$$
holds.
\end{proposition}
\par {\it Proof.} Clearly, property \rf{R-LMW} guarantees that the rectangle $\Rc_F[x,y:\tlm]\ne\emp$ for every $x,y\in\Mc$.
\par Let $y,y',y''\in\Mc$. Clearly,
$\Wc_F[y,y',y'':\tlm]=\Wc_F[y,y'',y':\tlm]$,
see \rf{WC-DF}, so that, without loss of generality, we may assume that $\rho(y,y')\le\rho(y,y'')$. Thanks to the triangle inequality, we have
$$
\rho(y',y'')\le\rho(y,y')+\rho(y,y'')\le 2\rho(y,y'')~~~~
\text{so that}~~~~~\rho(y,y')+\rho(y',y'')\le 3\rho(y,y'').
$$
Hence,
\be
\Wc_F[y,y',y'':3\tlm]&=&
\HR[\{F(y')+3\tlm\,\rho(y,y')\,Q_0\}\cap
\{F(y'')+3\tlm,\rho(y,y'')\,Q_0\}]
\nn\\
&\supset&
\HR[\{F(y')+\tlm\,\rho(y,y')\,Q_0\}\cap
\{(F(y'')+\tlm\,\rho(y',y'')Q_0)+
\tlm\,\rho(y,y')Q_0\}].
\nn
\ee
\par Clearly, for every $A,B\subset\RT$, $A\cap B\ne\emp$, and every $r\ge 0$, the following inclusion
$$
A\cap B+r Q_0\subset (A+r Q_0)\cap(B+r Q_0)
$$
holds. From this, property \rf{N-HS} and inequality $\lambda\le\tlm$, we have
\be
\Wc_F[y,y',y'':3\tlm]&\supset&
\HR[F(y')\cap
\{F(y'')+\tlm\,\rho(y',y'')Q_0\}+\tlm\,\rho(y,y')Q_0]
\nn\\
&=&\HR[F(y')\cap
\{F(y'')+\tlm\,\rho(y',y'')Q_0\}]+\tlm\,\rho(y,y')Q_0
\nn\\
&\supset&
\HR[F(y')\cap
\{F(y'')+\tlm\,\rho(y',y'')Q_0\}]+\lambda\,\rho(y,y')Q_0.
\nn
\ee
This and definition \rf{RL} imply the following inclusion:
\bel{W-0}
\Wc_F[y,y',y'':3\tlm]\supset
\RL[y',y'':\tlm]+\lambda\,\rho(y,y')Q_0.
\ee
\par Now, let us consider elements $x,x',x'',y,y',y''\in\Mc$. We may assume that $\rho(x,x')\le \rho(x,x'')$ and $\rho(y,y')\le \rho(y,y'')$. We know that in this case \rf{W-0} holds. In the same way we prove that
$$
\Wc_F[x,x',x'':3\tlm]\supset
\RL[x',x'':\tlm]+\lambda\,\rho(x,x')Q_0.
$$
\par From this inclusion, \rf{W-0} and the triangle inequality, we have
\be
\Ac&=&\Wc_F[x,x',x'':3\tlm]\cap\{\Wc_F[y,y',y'':3\tlm]
+\lambda\,\rho(x,y)Q_0\}
\nn\\
&\supset&
\{\RL[x',x'':\tlm]+\lambda\,\rho(x,x')Q_0\}\cap
\{\RL[y',y'':\tlm]+\lambda\,\rho(y,y')Q_0
+\lambda\,\rho(x,y)Q_0\}
\nn\\
&\supset&
\Bc=\{\RL[x',x'':\tlm]+\lambda\,\rho(x,x')Q_0\}\cap
\{\RL[y',y'':\tlm]+\lambda\,\rho(x,y')Q_0\}.
\nn
\ee
\par Thanks to \rf{R-LMW},
$$
\RL[x',x'':\tlm]\cap \{\RL[y',y'':\tlm]+\lambda\,\rho(x',y')\,Q_0\}\ne\emp
$$
so that there exists points $p_1\in \RL[x',x'':\tlm]$ and $p_2\in\RL[y',y'':\tlm]$ such that $\|p_1-p_2\|\le \lambda\,\rho(x',y')$. Therefore, thanks to the triangle inequality,
$$
\|p_1-p_2\|\le \lambda\,\rho(x',y')\le \lambda\,\rho(x,x')+\lambda\,\rho(x,y').
$$
\par This inequality implies the existence of a point $w\in[p_1,p_2]$ such that $\|p_1-w\|\le \lambda\,\rho(x,x')$ and $\|p_2-w\|\le \lambda\,\rho(x,y')$. Therefore, $w\in\Bc\subset\Ac$, and the proof of the proposition is complete.\bx
\msk


\par {\it Proof of Theorem \reff{CR-LS1}.} The necessity part of Theorem \reff{CR-LS1} and inequality $\inf\lambda\le \,|F|_{\Mf}$ are immediate from part (i) of Proposition \reff{P-GR-RT}.
\par Let us prove the sufficiency part. Thanks to Proposition \reff{WR-LM3}, condition \rf{R-MLM} of Theorem \reff{CR-LS1} implies condition \rf{WNEW} of Theorem \reff{W-CR} with  $\tlm=3\lambda$.
\par Theorem \reff{W-CR} tells us that there exists a Lipschitz selection $f$ of $F$ with
$$
\|f\|_{\Lip(\Mc)}\le 2\lambda+\tlm=2\lambda+3\lambda=5\lambda.
$$
This proves the sufficiency and inequality $|F|_{\Mf}\le 5\inf\lambda$.
\par The proof of Theorem \reff{CR-LS1} is complete.\bx
\msk

\begin{theorem}\lbl{FP-RT} Let $(\Mc,\rho)$ be a pseudometric space, and let $F:\Mc\to\CRT$ be a set-valued mapping. Suppose that $\Mf$ and $F$ satisfy Condition \reff{CND-T}.
\par The mapping $F$ has a Lipschitz selection if and only if there exists a constant $\lambda\ge 0$ such that for every subset $\Mc'\subset\Mc$ consisting of at most $N=4$ points, the restriction $F|_{\Mc'}$ of $F$ to $\Mc'$ has a Lipschitz selection $f_{\Mc'}$ with Lipschitz  seminorm $\|f_{\Mc'}\|_{\Lip(\Mc')}\le \lambda$.
\par Furthermore,
$\inf\lambda\le \,|F|_{\Mf}\le \gamma\inf\lambda$  with $\gamma=3$.
\end{theorem}
\par {\it Proof.} The necessity part of the theorem and the inequality $\inf\lambda\le \,|F|_{\Mf}$ are obvious. Let us prove the sufficiency.
\par Suppose that there exist a constant $\lambda\ge 0$ such that for every subset $\Mc'\subset\Mc$ with $\#\Mc'\le 4$, the restriction $F|_{\Mc'}$ of $F$ to $\Mc'$ has a Lipschitz selection $f_{\Mc'}$ with Lipschitz  seminorm $\|f_{\Mc'}\|_{\Lip(\Mc')}\le \lambda$. Our goal is to prove the existence of a selection $f$ of $F$ with $\|f\|_{\Lip(\Mc)}\le 3\lambda$.
\smsk
\par Let us show that condition \rf{WNEW} of Theorem \reff{W-CR} is satisfied with $\tlm=\lambda$.
\par Let $x,x',x'',y,y',y''\in\Mc$, and let $\Mc'=\{x',x'',y',y''\}$. We know that the restriction $F|_{\Mc'}$ of $F$ to $\Mc'$ has a Lipschitz selection $f_{\Mc'}$ with $\|f_{\Mc'}\|_{\Lip(\Mc')}\le \lambda$.
Thus,
$$
f_{\Mc'}(x')\in F(x'),~ f_{\Mc'}(x'')\in F(x''),~ f_{\Mc'}(y')\in F(y')~~~\text{and}~~~f_{\Mc'}(y'')\in F(y'').
$$
\par It is well known that any Lipschitz mapping from a closed subset of a pseudometric space $\MR$ into $\LTI$ (i.e., $\RT$ equipped with  the uniform norm) can be extended to all of the space $\MR$ with preserving the Lipschitz constant. Let $\tf:\Mc\to\RT$ be such an extension of the mapping $f_{\Mc'}$ from $\Mc'$ to $\Mc$. Thus, $\tf|_{\Mc'}=f_{\Mc'}$ and $\|\tf\|_{\Lip(\Mc)}
=\|f_{\Mc'}\|_{\Lip(\Mc')}\le \lambda$ so that the following inequalities hold:
$$
\|f_{\Mc'}(x')-\tf(x)\|=\|\tf(x')-\tf(x)\|
\le\lambda\,\rho(x,x'),
~~~~~~
\|f_{\Mc'}(x'')-\tf(x)\|=\|\tf(x'')-\tf(x)\|\le \lambda\,\rho(x,x''),
$$
and
$$
\|f_{\Mc'}(y')-\tf(y)\|=\|\tf(y')-\tf(y)\|\le \lambda\,\rho(y,y'),
~~~~~~
\|f_{\Mc'}(y'')-\tf(y)\|=
\|\tf(y'')-\tf(y)\|\le\lambda\,\rho(y,y'').
$$
Hence,
$$
\tf(x)\in \{F(x')+\lambda\,\rho(x',x)\,Q_0\}
\cap \{F(x'')+\lambda\,\rho(x'',x)\,Q_0\}
$$
and
$$
\tf(y)\in \{F(y')+\lambda\,\rho(y',y)\,Q_0\}
\cap \{F(y'')+\lambda\,\rho(y'',y)\,Q_0\}
$$
so that $\tf(x)\in\Wc_F[x,x',x'':\lambda]$ and $\tf(y)\in\Wc_F[y,y',y'':\lambda]$, see \rf{WC-DF}.
\smsk
\par Furthermore, because $\|\tf\|_{\Lip(\Mc)}\le \lambda$,
we have $\|\tf(x)-\tf(y)\|\le \lambda\,\rho(x,y)$. Hence, 
$$
\tf(x)\in\Wc_F[x,x',x'':\tlm]\cap \{\Wc_F[y,y',y'':\tlm]+\lambda\,\rho(x,y)\,Q_0\}
$$
proving that condition \rf{WNEW} of Theorem \reff{W-CR} holds. Thanks to this theorem, $F$ has a Lipschitz selection $f:\Mc\to\RT$ with $\|f\|_{\Lip(\Mc)}\le 2\lambda+\tlm=3\lambda$, and the proof of Theorem \reff{FP-RT} is complete.\bx
\msk
\begin{remark} {\em We have proved Theorem \reff{FP-RT}
for the Banach space $\LTI=(\RT,\|\cdot\|)$, i.e., $\RT$
equipped with the uniform norm $\|\cdot\|$. Let us consider a slightly more general version of Theorem \reff{FP-RT} for a Banach space $X=(\RT,\|\cdot\|_X)$ supplied with a certain Banach norm $\|\cdot\|_X$.
\par Let $B_X$ be the unit ball $X$. According to a result of E. Asplund \cite{Asp-1960}, $B_X$ contains a parallelogram $P$ centered at $(0,0)$ which expanded by $\frac{3}{2}$ will cover $B_X$. Since the Banach space $X_P$ with the unit ball $P$ is linearly isometric to $\LTI$,
we conclude that an analog of Theorem \reff{FP-RT} holds for an {\it arbitrary} Banach space $X$ with the constant $\frac{3}{2}\cdot 3\lambda=4.5\lambda$ (instead of $3\lambda$ as for $\LTI$).
\rbx}
\end{remark}
\par Combining Theorem \reff{W-CR} with the part (ii) of Proposition \reff{P-GR-RT}, we obtain the following criterion for Lipschitz selections.
\begin{theorem}\lbl{CR-W5} Let $\Mf=\MR$ be a pseudometric space, and let $F:\Mc\to\CRT$ be a set-valued mapping satisfying Condition \reff{CND-T}.
\par The mapping $F$ has a Lipschitz selection if and only if there exists a constant $\lambda\ge 0$ such that for every $x,x',x'',y,y',y''\in\Mc$, we have
\bel{WT-A2}
\Wc_F[x,x',x'':\lambda]\cap \{\Wc_F[y,y',y'':\lambda]+\lambda\,\rho(x,y)\,Q_0\}\ne\emp.
\ee
\par Furthermore, the following inequality holds:
\bel{NF-WQ}
\inf\lambda\le \,|F|_{\Mf}\le 3\inf\lambda.
\ee
\end{theorem}
\par Comparing this result with Theorem \reff{CR-LS1}, we note that Theorem \reff{CR-W5} provides better upper bound for $|F|_{\Mf}$ ($3\inf\lambda$ in \rf{NF-WQ} rather than
$5\inf\lambda$ as in inequality \rf{OP-NF}). The price of such a refinement is as follows: property \rf{WT-A2} exploits rectangles $\Wc_F[\cdot,\cdot,\cdot:\lambda]$ each of these depend on {\it three} arbitrary elements of $\Mc$. In turn, condition \rf{R-MLM} of Theorem \reff{CR-LS1} formulates in terms of rectangles $\RL[\cdot,\cdot:\lambda]$ each depending only on {\it two} elements of $\Mc$. In other words, in the criterion of Theorem \reff{CR-W5}, we use more information about geometric properties of the set-valued mapping $F$, and, as a result, obtain better estimates for the Lipschitz seminorm of its Lipschitz selection.
\msk

\indent\par {\bf 5.2 Criteria for Lipschitz selections in  terms of intersections of sets.}
\addtocontents{toc}{~~~~5.2 Criteria for Lipschitz selections in  terms of intersections of sets.\hfill\thepage\par\VST}
\msk
\indent
\par In this section we prove two constructive Lipschitz selection criteria formulated in terms of intersections of certain families of rectangles in $\RT$. We begin with a Lipschitz selection criterion of such a kind for mappings taking values in the family $\RCT$ of all rectangles in $\RT$ with sides parallel to the coordinate axes.
\begin{proposition} Let $\Tc:\Mc\to\Rc(\RT)$ be a set-valued mapping. Let us assume that either $\Mc$ is finite or all rectangles $\Tc(x)$, $x\in\Mc$, are closed and at least one of them is bounded.
\par Then $\Tc$ has a Lipschitz selection if and only if there exists a constant $\lambda\ge 0$ such that the set
\bel{TC-DF-1}
\Tc^{[1]}[x:\lambda]=
\bigcap_{z\in \Mc}\,
\left[\Tc(z)+\lambda\,\rho(x,z)\,Q_0\right]
\ee
is non-empty for every $x\in\Mc$. Furthermore,
$$
|\Tc|_{\Mf}=
\inf\{\lambda: \Tc^{[1]}[x:\lambda]\ne\emp~~~\text{for all}~~~x\in\Mc\}.
$$
\end{proposition}
\par {\it Proof.} {\it (Necessity.)} Suppose that $\Tc$ has a Lipschitz selection $\tau:\Mc\to\RCT$ with $\|\tau\|_{\Lip(\Mc)}\le\lambda$. Then, given $x\in\Mc$, we have  $\|\tau(x)-\tau(z)\|\le\lambda\,\rho(x,z)$ for every $z\in\Mc$. But $\tau(x)\in \Tc(x)$ and $\tau(z)\in \Tc(z)$ (because $\tau$ is a selection of $\Tc$) so that $\tau(x)\in \Tc(z)+\lambda\,\rho(x,z) Q_0$. Therefore, thanks to \rf{TC-DF-1}, $\tau(x)\in \Tc^{[1]}[x:\lambda]$ proving the necessity.
\par {\it (Sufficiency.)} Suppose that $\Tc^{[1]}[x:\lambda]\ne\emp$ for every $x\in\Mc$. Then, thanks to definition \rf{TC-DF-1}, $\Tc(x)\cap[\Tc(z)+\lambda\,\rho(x,z)\,Q_0]\ne\emp$ for every $x,z\in\Mc$, proving the existence of points
$\tta(x)\in \Tc(x)$, $\tta(z)\in \Tc(z)$ such that $\|\tta(x)-\tta(z)\|\le\lambda\,\rho(x,z)$. This property and the assumption of the proposition's hypothesis enables us to apply Proposition \reff{FP-RC} to the set-valued mapping $\Tc$. This proposition tells us that $\Tc$ has a Lipschitz selection $\tau:\Mc\to\RT$ with  $\|\tau\|_{\Lip(\Mc)}\le \eta$ completing the proof of the proposition.\bx
\msk

\par The following theorem is the main result of this section.
\begin{theorem}\lbl{CR-LS-RI} Suppose that a pseudometric space $\Mf=\MR$ and a set-valued mapping $F:\Mc\to\CRT$ satisfy Condition \reff{CND-T}.
\par Then $F$ has a Lipschitz selection if and only if there exists a constant $\lambda\ge 0$ such that
\bel{RC-CH1}
\bigcap_{y,y'\in\Mc}
\left\{\,\RL[y,y':\lambda]+\lambda\,\rho(x,y) Q_0\right\} \ne\emp
~~~~~~\text{for every}~~x\in\Mc.
\ee
\par Moreover, in this case, $\inf\lambda\le \,|F|_{\Mf}\le 5\inf\lambda$.
\end{theorem}
\par {\it Proof.}  {\it (Necessity).} Let $f:\Mc\to\RT$ be a Lipschitz selection of $F$, and let $\lambda=\|f\|_{\Lip(\Mc)}$. Then for every $x,y,y'\in\Mc$ the following is true: $f(x)\in F(x),f(y)\in F(y)$, $f(y')\in F(y')$,
$$
\|f(x)-f(y)\|\le \lambda\,\rho(x,y)~~~\text{and}~~~
\|f(y)-f(y')\|\le \lambda\,\rho(y,y').
$$
\par Hence,
$
f(y)\in F(y)\cap [F(y')+\lambda\,\rho(y,y') Q_0]
$
and $f(x)\in f(y)+\lambda\,\rho(x,y) Q_0,$ proving that
$$
f(x)\in (F(y)\cap [F(y')+\lambda\,\rho(y,y') Q_0])
+\lambda\,\rho(x,y) Q_0~~~~\text{ for all}~~~~y,y'\in\Mc.
$$
Clearly,
\be
(F(y)\cap\{F(y')+\lambda\,\rho(y',y) Q_0\}) +\lambda\,\rho(x,y) Q_0&\subset&
\HR[F(y)\cap\{F(y')+\lambda\,\rho(y',y) Q_0\}] +\lambda\,\rho(x,y) Q_0
\nn\\
&=&
\RL[y,y':\lambda]+\lambda\,\rho(x,y) Q_0
\nn
\ee
so that $f(x)\in \RL[y,y':\lambda]+\lambda\,\rho(x,y) Q_0$ for every $y,y'\in\Mc$.
\par This property implies \rf{RC-CH1} and inequality $\inf\lambda\le \,|F|_{\Mf}$ completing the proof of the necessity.
\msk
\par {\it (Sufficiency.)} We assume that property \rf{RC-CH1} holds. Thanks to this property,
$$
\RL[x,x':\lambda]\cap \{\RL[y,y':\lambda]+\lambda\,\rho(x,y)\,Q_0\}\ne\emp
$$
for every $x,x',y,y'\in\Mc$ proving that condition \rf{R-MLM} of Theorem \reff{CR-LS1} holds. This theorem tells us that $F$ has a Lipschitz selection with Lipschitz seminorm at most $5\lambda$. This proves the inequality  $|F|_{\Mf}\le 5\inf\lambda$ and completes the proof of Theorem \reff{CR-LS-RI}.\bx

\begin{remark} {\em We note that, thanks to representation \rf{P1XP2}, property \rf{RC-CH1} is equivalent to the following one: for every $x\in\Mc$ and
and every $i=1,2$, the set
$$
\bigcap_{y,y'\in\Mc}
\Prj_i
[\left\{F(y)\cap(F(y')+\lambda\,\rho(y,y')Q_0)\right\}
+\lambda\,\rho(x,y) Q_0]\ne\emp.
$$
\par We recall that $\Prj_i$ denotes the operator of the orthogonal projection onto the axis $Ox_i$.\rbx}
\end{remark}

\SECT{6. Projection Algorithm for nearly optimal Lipschitz selections.}{6}
\addtocontents{toc}{6. Projection Algorithm for nearly optimal Lipschitz selections.\hfill \thepage\par\VST}
\indent
\par In this section we present a number of nearly optimal algorithms for Lipschitz selections based on the geometrical construction suggested in the proof of Theorem \reff{W-CR}. We refer to these algorithms as the {\it ``Projection Algorithms''}.
\smsk

\indent\par {\bf 6.1 The $\lmv$-Projection Algorithm.}
\addtocontents{toc}{~~~~6.1 The $\lmv$-Projection Algorithm.\hfill \thepage\par\VST}

\smsk
\par Let $\Mf=\MR$ be a pseudometric space, and let $F:\Mc\to\CRT$ be a set-valued mapping. Given a constant $\lambda\ge 0$ we define a set-valued mapping $F^{[1]}[\cdot:\lambda]$ on $\Mc$ by letting
\bel{F-1}
F^{[1]}[x:\lambda]=\bigcap_{y\in \Mc}\,
\left[F(y)+\lambda\,\rho(x,y)\,Q_0\right],~~~~x\in\Mc.
\ee
\par We refer to the mapping $F^{[1]}[\cdot:\lambda]$ as {\it the $\lambda$-balanced refinement of $F$}. (This important object was introduced in \cite{FS-2018} and \cite{S-2022}. We have already met various variants of it in \rf{TAU-1}, \rf{F1-PR1}, \rf{WF1-3L}, \rf{T-ONE} and \rf{TC-DF-1}.)
\begin{pra}\lbl{PA-DS} {\em Given a vector $\lmv=(\lambda_1,\lambda_2)$ with non-negative coordinates $\lambda_1,\lambda_2\ge 0$, a  pseudometric space $\Mf=\MR$, and a set-valued mapping $F:\Mc\to\CRT$, the Projection Algorithm either produces a selection $f_{\lmv;F}$ of $F$ (the outcome {\bf ``Success''}) or stops (the outcome {\bf ``No go''}). This procedure includes the following five main steps.
\msk
\par {\bf STEP 1.} At this step we construct the $\lambda_1$-balanced refinement of $F$, i.e., the set-valued mapping
\bel{F-LM1}
F^{[1]}[x:\lambda_1]=\bigcap_{y\in \Mc}\,
\left[F(y)+\lambda_1\,\rho(x,y)\,Q_0\right].
\ee
\par {\it If $F^{[1]}[x:\lambda_1]=\emp$ for some $x\in\Mc$, the algorithm produces the outcome {\bf ``No go''} and stops.}
\msk
\par {\bf STEP 2.} Let us assume that the above condition does not hold, i.e., {\it for every element $x\in\Mc$ the $\lambda_1$-balanced refinement $F^{[1]}[x:\lambda_1]$ is non-empty}.
\smsk
\par In this case, for each $x\in\Mc$, we construct {\it the rectangular hull} of $F^{[1]}[x:\lambda_1]$, the set
\bel{TC-DF}
\Tc_{F,\lambda_1}(x)=\Hc[F^{[1]}[x:\lambda_1]].
\ee
\par Thus, $\Tc_{F,\lambda_1}$ is a set-valued mapping which maps $\Mc$ into the family $\RCT$ of all rectangles in $\RT$. See Fig. 14.
\msk
\begin{figure}[h!]
\hspace{30mm}
\includegraphics[scale=0.64]{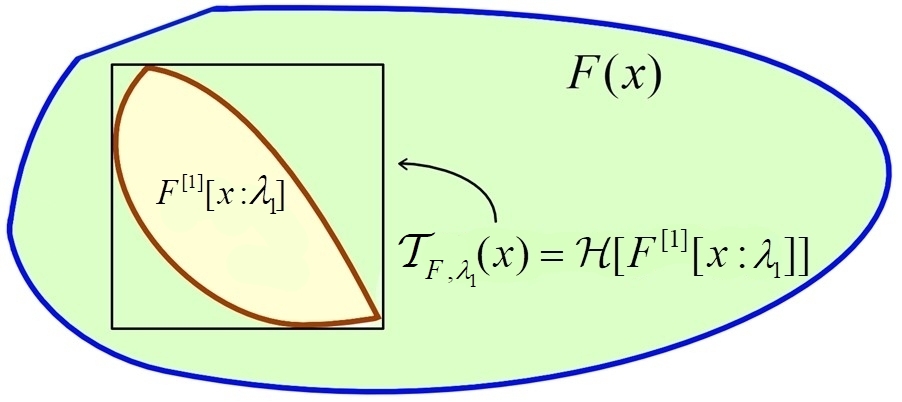}
\caption{The second step of the Projection Algorithm.}
\end{figure}
%
\par {\bf STEP 3.} For every $x\in\Mc$, we construct the $\lambda_2$-balanced refinement of the mapping $\Tc_{F,\lambda_1}$, i.e., the rectangle $\Tc^{[1]}_{F,\lambda_1}[x:\lambda_2]$ defined by
\bel{RT-1}
\Tc^{[1]}_{F,\lambda_1}[x:\lambda_2]=\bigcap_{y\in \Mc}\,
\left[\Tc_{F,\lambda_1}(y)+\lambda_2\,\rho(x,y)\,Q_0\right].
\ee
See Fig. 15.
\begin{figure}[h!]
\hspace{33mm}
\includegraphics[scale=0.65]{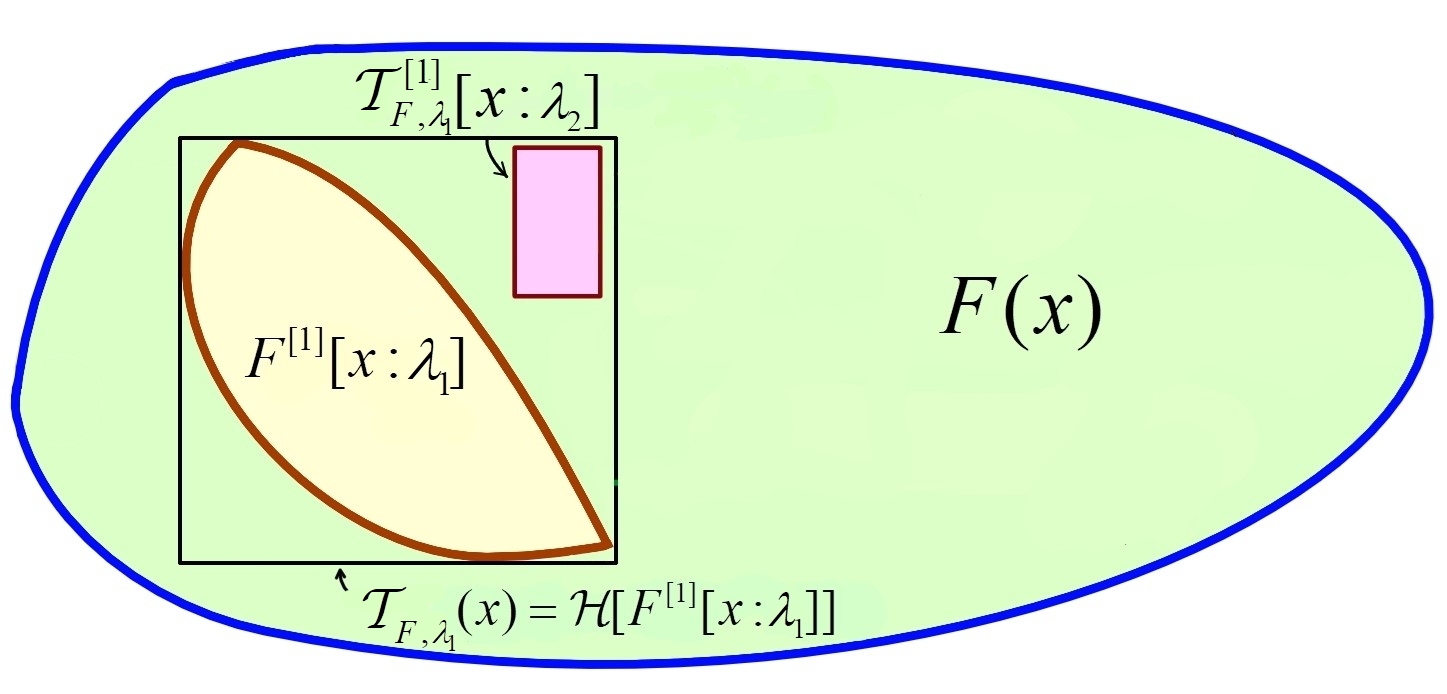}
\caption{The third step of the Projection Algorithm.}
\end{figure}
\par If
\bel{S-ST}
\Tc^{[1]}_{F,\lambda_1}[x:\lambda_2]=\emp~~~~\text{for some}~~~x\in\Mc,
\ee
the algorithm produces the outcome {\bf ``No go''} and stops.
\bsk
\par {\bf STEP 4.} At this step, we assume that {\it for each $x\in\Mc$ the rectangle $\Tc^{[1]}_{F,\lambda_1}[x:\lambda_2]\ne\emp$}.
We let $\Tc^{[1]}_{F,\lambda_1}[x:\lambda_2]^{\cl}$ denote {\it the closure} of the rectangle $\Tc^{[1]}_{F,\lambda_1}[x:\lambda_2]\ne\emp$.
\smsk
\par Let $O=(0,0)$ be the origin. We define a mapping $g_F:\Mc\to\RT$ by letting
\smsk
\bel{G-CNT}
g_F(x)=\cent\left(\Prj
\left(O,\Tc^{[1]}_{F,\lambda_1}[x:\lambda_2]^{\cl}\right)\right), ~~~~~~x\in\Mc.
\hspace{18mm}
\ee
\par See Fig. 16.
\begin{figure}[h!]
\hspace{25mm}
\includegraphics[scale=0.75]{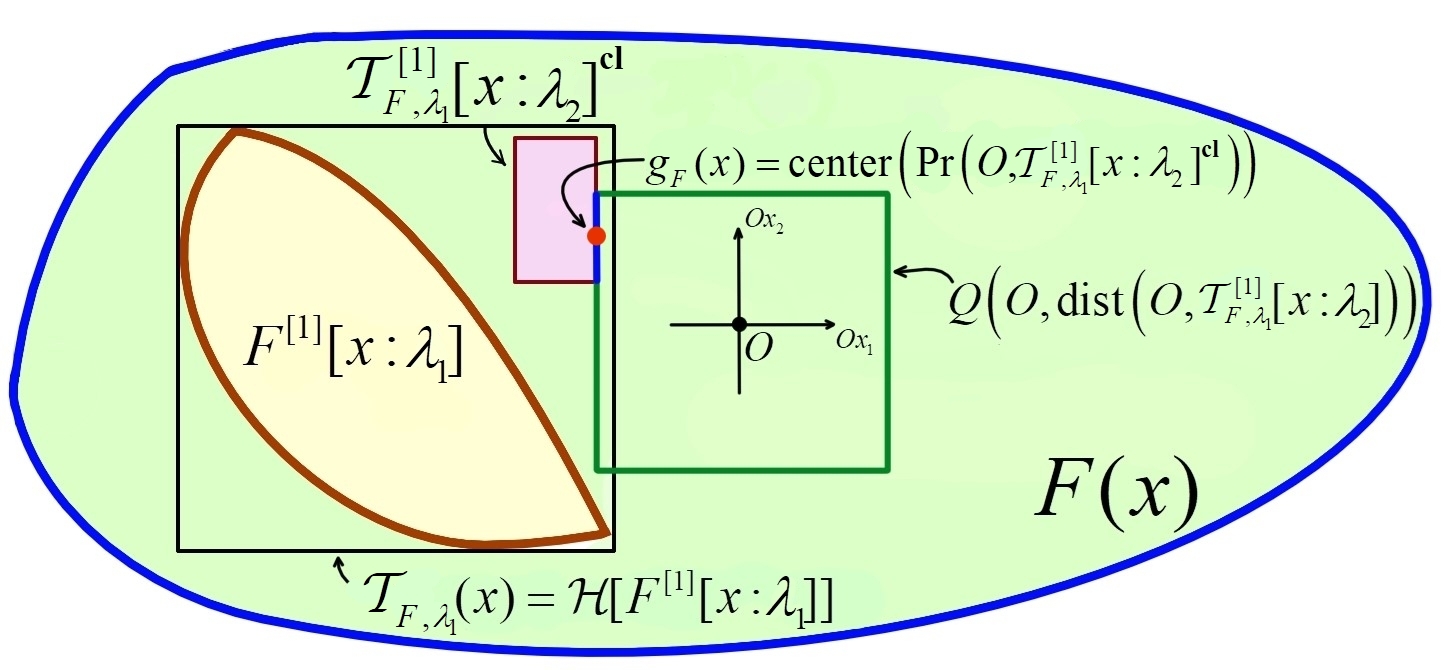}
\caption{The mapping $g_F:\Mc\to\RT$, a selection of the set-valued mapping $\Tc_{F,\lambda_1}^{\cl}$.}
\end{figure}
\par  We recall that $\Prm(\cdot,S)$ denotes the operator of metric projection onto a closed convex subset $S\subset\RT$, see \rf{MPR}, and  $\cent\,(\cdot)$ denotes the center of a centrally symmetric bounded set in $\RT$.\\

\par {\bf STEP 5.} We define the mapping $f_{\lmv;F}:\Mc\to\RT$ by letting\msk
\bel{F-LF}
f_{\lmv;F}(x)=\Prm(g(x),F^{[1]}[x:\lambda_1]), ~~~~~~x\in\Mc.
\ee
See Fig. 17.
\msk
\begin{figure}[h!]
\hspace{25mm}
\includegraphics[scale=1.18]{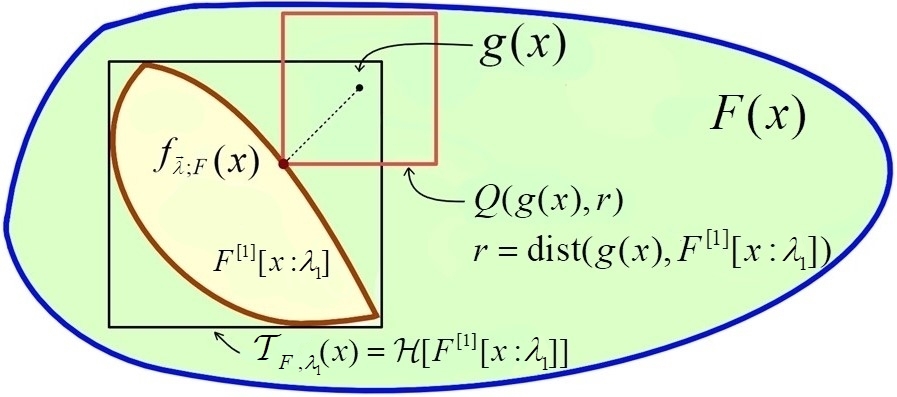}
\caption{The final step of the Projection Algorithm.}
\end{figure}
\smsk
\par At this stage, the algorithm produces the outcome
{\bf ``Success''} and stops.\bx}
\end{pra}
\par To specify the dependence on the parameters $\lambda_1$ and $\lambda_2$, we refer to the above algorithm as the $\lmv$-Projection Algorithm ($\lmv$-PA for short).

\begin{remark}\lbl{REM-SC} {\em We note that the $(\lambda_1,\lambda_2)$-Projection Algorithm
produces the outcome {\bf ``Success''} if and only if
the following conditions are satisfied:
\bel{SC-P}
F^{[1]}[x:\lambda_1]\neq\emp~~~\text{and}~~~
\Tc^{[1]}_{F,\lambda_1}[x:\lambda_2]\neq\emp~~~~
\text{for every}~~~x\in\Mc.
\ee
\par Indeed, this algorithm does not stop at {\bf STEP 1} which is equivalent to the first condition in \rf{SC-P}. Also, the algorithm does not stop at {\bf STEP 3} which, thanks to \rf{S-ST}, is equivalent to the second condition in \rf{SC-P}.\rbx}
\end{remark}
\par The next theorem describes the main properties of
the $\lmv$-PA.
\begin{theorem}\lbl{ALG-T} Let $\lambda_1,\lambda_2\ge 0$, and let $\lmv=(\lambda_1,\lambda_2)$. Let $\Mf=\MR$ be a pseudometric space, and let $F:\Mc\to\CRT$ be a set-valued mapping. Suppose that $\Mf$ and $F$ satisfy Condition \reff{CND-T}.
\smsk
\par If $\lmv$-Projection Algorithm produces the outcome {\bf ``No go''} (see {\bf STEP 1} and {\bf STEP 3}), then we gua\-ran\-tee that there does not exist a Lipschitz selection of $F$ with Lipschitz seminorm at most $\min\{\lambda_1,\lambda_2\}$.

\par  Otherwise, the $\lmv$-Projection Algorithm produces the outcome {\bf ``Success''} (see {\bf STEP 5})
and returns the mapping $f_{\lmv;F}:\Mc\to\RT$
defined by formula \rf{F-LF}. This mapping has the following properties:
\smsk
\par ($\bigstar A$) The mapping $f_{\lmv;F}$ is well defined. This means the following: (a) for every $x\in\Mc$, the set $F^{[1]}[x:\lambda_1]\ne\emp$ and the rectangle $\Tc^{[1]}_{F,\lambda_1}[x:\lambda_2]\ne\emp$ (see \rf{RT-1}), (b) the mapping $g_F$ (see \rf{G-CNT}) is well defined, and (c) the metric projection defined by the right hand side of \rf{F-LF} is a {\it singleton}.
\smsk
\par ($\bigstar B$) $f_{\lmv;F}$ is a Lipschitz selection of $F$ with the Lipschitz seminorm $$\|f_{\lambda;F}\|_{\Lip(\Mc)}\le \lambda_1+2\lambda_2.$$
\end{theorem}
\par {\it Proof.} Let us see that whenever the $\lmv$-PA produces the outcomes {\bf ``No go''} there does not exist a Lipschitz selection of $F$ with Lipschitz seminorm at most $\lambda^{\min}=\min\{\lambda_1,\lambda_2\}$.
\par Indeed, suppose that $F$ has a Lipschitz selection $f:\Mc\to\RT$ with $\|f\|_{\Lip(\Mc)}\le\lambda^{\min}$.
Then, $f(x)\in F(x)$ for every $x\in\Mc$ and
$$
\|f(x)-f(y)\|\le\lambda^{\min}\,\rho(x,y)~~~~\text{for all}
~~~~x,y\in\Mc.
$$
Hence,
$$
f(x)\in F(y)+\lambda^{\min}\,\rho(x,y)Q_0~~~\text{for every}~~~y\in\Mc.
$$
Therefore, thanks to definition \rf{F-1}, we have
$$
f(x)\in F^{[1]}[x:\lambda^{\min}]\subset F^{[1]}[x:\lambda_1]~~~~\text{for every}~~~~x\in\Mc.
$$
\par Furthermore, from this and \rf{TC-DF},
we have $f(x)\in \Tc_{F,\lambda_1}(x)$, $x\in\Mc$, proving that $f$ is a Lipschitz selection of the set-valued mapping $\Tc_{F,\lambda_1}$ with $\|f\|_{\Lip(\Mc)}\le\lambda^{\min}$.
\par Repeating the above argument for $\Tc_{F,\lambda_1}$ rather than for $F$, we conclude that
$$
f(x)\in \Tc^{[1]}_{F,\lambda_1}[x:\lambda^{\min}]
\subset \Tc^{[1]}_{F,\lambda_1}[x:\lambda_2]~~~~
\text{for every}~~~~x\in\Mc.
$$
\par See \rf{RT-1}. Thus, if there exists $x\in\Mc$ such that
$$
F^{[1]}[x:\lambda_1]=\emp~~~~~\text{(the outcome {\bf ``No go''} at {\bf STEP 1})}
$$
or
$$
\Tc^{[1]}_{F,\lambda_1}[x:\lambda_2]=\emp~~~~~\text{(the outcome {\bf ``No go''} at {\bf STEP 3})},
$$
we can guarantee that there does not exist a Lipschitz selection of $F$ with Lipschitz seminorm $\le\lambda^{\min}$.
\smsk
\par Next, let us assume that the $\lambda$-PA produces the outcome {\bf ``Success''} (see {\bf STEP 5}). Then, thanks to Remark \reff{REM-SC}, the set-valued mapping $F$ satisfies condition \rf{SC-P}.
\par Let us prove property ($\bigstar A$). Clearly, the part (a) is immediate from \rf{SC-P}. To prove part (b), we note that metric projection of the origin $O$ onto the rectangle $\Tc^{[1]}_{F,\lambda_1}[x:\lambda_2]^{\cl}$ is either a singleton or a closed line segment. Therefore, $g_F$ is well defined on $\Mc$ (because the set $\Prm(O,\Tc^{[1]}_{F,\lambda_1}[x:\lambda_2])$ is a non-empty {\it bounded} centrally symmetric set). See formula \rf{G-CNT}.
\par Let us prove part (c). Thanks to \rf{TC-DF} and \rf{RT-1},
\bel{G-IN-T}
g_F(x)\in\Tc^{[1]}_{F,\lambda_1}[x:\lambda_2]^{\cl}
\subset\Tc_{F,\lambda_1}^{\cl}(x)
=\Hc[F^{[1]}[x:\lambda_1]]^{\cl}.
\ee
(To prove the inclusion ``$\subset$'' in \rf{G-IN-T} put $y=x$ in \rf{RT-1}.) Therefore, thanks to Lemma \reff{WMP-S} and definition \rf{G-CNT}, $f_{\lmv;F}(x)$, the metric projection of $g_F(x)$ onto $F^{[1]}[x:\lambda_1]$, see \rf{F-LF}, is a singleton. The proof of property ($\bigstar A$) is complete.
\smsk
\par Let us prove property ($\bigstar B$) of the theorem.
We will follow the scheme of the proof of the sufficiency part of the key Theorem \reff{W-CR} for the special case
\bel{LM-12}
\tlm=\lambda_1~~~~\text{and}~~~~\lambda=\lambda_2.
\ee
\par Thanks to \rf{SC-P}, $\Tc^{[1]}_{F,\tlm}[x:\lambda]\ne\emp$ for every $x\in\Mc$. Part (a) of Proposition \reff{X2-C} tells us that in this case the set-valued mapping $\Tc^{[1]}_{F,\tlm}[\cdot:\lambda]$ is Lipschitz with respect to the Hausdorff distance, i.e.,
\bel{T-HD}
\dhf\left(\Tc^{[1]}_{F,\tlm}[x:\lambda],
\Tc^{[1]}_{F,\tlm}[y:\lambda]\right)
\le \lambda\,\rho(x,y)~~~\text{for all}~~~x,y\in\Mc.
\ee
\par It is clear that for every two rectangles $S_1,S_2\in\RCT$ we have $\dhf(S_1,S_2)=\dhf(S_1^{\cl},S_2^{\cl})$, so that, thanks to \rf{T-HD},
\bel{TCL-HD}
\dhf\left(\Tc^{[1]}_{F,\tlm}[x:\lambda]^{\cl},
\Tc^{[1]}_{F,\tlm}[y:\lambda]^{\cl}\right)
\le \lambda\,\rho(x,y)~~~\text{for every}~~~x,y\in\Mc.
\ee
\par Let $\delta(x)=\dist(0,\Tc^{[1]}_{F,\tlm}[x:\lambda])$. Clearly, $\delta(x)=\dist(0,\Tc^{[1]}_{F,\tlm}[x:\lambda]^{\cl})$. Then,
$$
|\delta(x)-\delta(y)|=
|\dist(0,\Tc^{[1]}_{F,\tlm}[x:\lambda])-
\dist(0,\Tc^{[1]}_{F,\tlm}[y:\lambda])|
\le \dhf\left(\Tc^{[1]}_{F,\tlm}[x:\lambda],
\Tc^{[1]}_{F,\tlm}[y:\lambda]\right)
$$
proving that
\bel{L-DLT}
|\delta(x)-\delta(y)|\le \lambda\,\rho(x,y)~~~
\text{for every}~~~x,y\in\Mc.
\ee
\par We note that
$$
\Prm(O,\Tc^{[1]}_{F,\tlm}[x:\lambda]^{\cl})=
Q(O,\delta(x))\cap \Tc^{[1]}_{F,\tlm}[x:\lambda]^{\cl}.
$$
From this and Lemma \reff{H-RB}, we have
\be
\Prm(O,\Tc^{[1]}_{F,\tlm}[x:\lambda]^{\cl})+\lambda\,\rho(x,y)Q_0
&=&
\{Q(O,\delta(x))+\lambda\,\rho(x,y)Q_0\}\cap \{\Tc^{[1]}_{F,\tlm}[x:\lambda]^{\cl}+\lambda\,\rho(x,y)Q_0\}
\nn\\
&=&
Q(O,\delta(x)+\lambda\,\rho(x,y))\cap \{\Tc^{[1]}_{F,\tlm}[x:\lambda]^{\cl}+\lambda\,\rho(x,y)Q_0\}.
\nn
\ee
\par Note that, thanks to \rf{L-DLT}, $\delta(y)\le\delta(x)+\lambda\,\rho(x,y)$, and, thanks to \rf{TCL-HD},
$$
\Tc^{[1]}_{F,\tlm}[y:\lambda]^{\cl}\subset\
\Tc^{[1]}_{F,\tlm}[x:\lambda]^{\cl}+\lambda\,\rho(x,y)Q_0.
$$
Hence,
$$
\Prm(O,\Tc^{[1]}_{F,\tlm}[x:\lambda]^{\cl})
+\lambda\,\rho(x,y)Q_0
\supset
Q(O,\delta(y))\cap \Tc^{[1]}_{F,\tlm}[y:\lambda]^{\cl}=
\Prm(O,\Tc^{[1]}_{F,\tlm}[y:\lambda]^{\cl}).
$$
\par By interchanging the roles of $x$ and $y$, we obtain also
$$
\Prm(O,\Tc^{[1]}_{F,\tlm}[y:\lambda]^{\cl})
+\lambda\,\rho(x,y)Q_0
\supset\Prm(O,\Tc^{[1]}_{F,\tlm}[y:\lambda]^{\cl}).
$$
These two inclusions imply the inequality
\bel{TH-PR}
\dhf\left(\Prm(O,\Tc^{[1]}_{F,\tlm}[x:\lambda]^{\cl}),
\Prm(O,\Tc^{[1]}_{F,\tlm}[y:\lambda]^{\cl})\right)
\le \lambda\,\rho(x,y),~~~~x,y\in\Mc.
\ee
\par As we have noted above, the set $\Prm(O,\Tc^{[1]}_{F,\tlm}[x:\lambda])$ is either a singleton or a closed line segment. This line segment lies on one of the sides of the square $Q(O,\delta(x))$, proving that the metric projection $\Prm(O,\Tc^{[1]}_{F,\tlm}[x:\lambda]^{\cl})$ is a {\it bounded rectangle}, i.e., an element of the family $\RCT$. Therefore, thanks to part (ii) of Claim \reff{TWO} and \rf{TH-PR}, the mapping $g_F$ defined by formula \rf{G-CNT} is Lipschitz with Lipschitz seminorm at most $\lambda$.
\smsk
\par Let us also not that, thanks to \rf{G-IN-T} and \rf{LM-12}, for every $x\in\Mc$, we have
$$
g_F(x)\in\Tc_{F,\tlm}^{\cl}(x)=\Hc[F^{[1]}[x:\tlm]]^{\cl}.
$$
These properties of $g_F$ show that the mapping $g=g_F$ satisfies conditions \rf{WG-HF} and \rf{G-LIP}.
\smsk
\par Furthermore, comparing \rf{FX-PR} with \rf{F-LF}, we conclude that $f=f_{\lmv;F}$ where $f$ is the mapping defined by the formula \rf{FX-PR}. Proposition \reff{LSEL-F} tells us that this mapping is a selection of $F$ with Lipschitz seminorm
$$
\|f\|_{\Lip(\Mc)}=\|f_{\lmv;F}\|_{\Lip(\Mc)}\le \tlm +2\lambda=\lambda_1+2\lambda_2.
$$
\par The proof of Theorem \reff{ALG-T} is complete.\bx

\begin{remark} {\em Theorem \reff{ALG-T} holds for various versions of the $\lmv$-Projection Algorithm related to the specific choice of the mapping $g_F$ at {\bf STEP 4} of the algorithm. The proof of Theorem \reff{ALG-T} shows that $g_F$ should be a Lipschitz selection of the set-valued mapping
$$
\Tc_{F,\lambda_1}^{\cl}(x)=
\HR[F^{[1]}[\cdot:\lambda_1]]^{\cl}~~~~
\text{with}~~~~\|g_F\|_{\Lip(\Mc)}\le\lambda_2,
$$
i.e., $g_F$ have to satisfy condition \rf{WG-HF} and inequality \rf{G-LIP} with constants $\tlm=\lambda_1$ and $\lambda=\lambda_2$.
\par Let us note some of these versions.
\smsk
\par {\it(i)} We can define $g_F(x)$ by formula \rf{G-CNT} with replacing the origin $O$ by an {\it arbitrary} point in $\RT$;
\msk
\par  {\it(ii)} Suppose that
\bel{TC-BO}
\text{for every}~~x\in\Mc~~\text{the set}~~ F^{[1]}[x:\lambda_1]~~\text{is bounded}.
\ee
Then the rectangle $\Tc^{[1]}_{F,\lambda_1}[x:\lambda_2]$, (i.e., the  rectangular hull of $F^{[1]}[x:\lambda_1]$) is also bounded for all $x\in\Mc$. In this case, we can define $g_F$ by the formula
\bel{GF-CNT}
g_F(x)=\cent
\left(\Tc^{[1]}_{F,\lambda_1}[x:\lambda_2]\right),
~~~~~x\in\Mc.
\ee
See Fig. 18.
\msk
%
%
\begin{figure}[h!]
\hspace{30mm}
\includegraphics[scale=0.65]{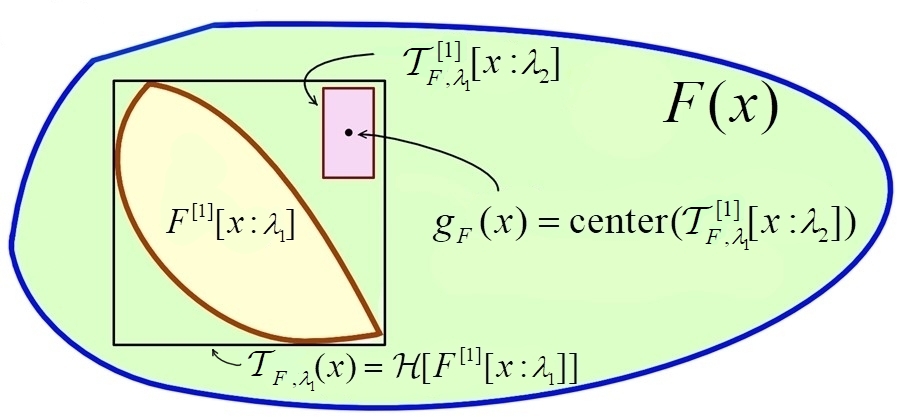}
\caption{{\bf STEP 4} of the Projection Algorithm: the case of bounded sets $F^{[1]}[x:\lambda_1]$.}
\end{figure}
\msk
\par Then, thanks to part (iii) of Proposition \reff{WLS-T}, $g_F$ is Lipschitz with $\|g_F\|_{\Lip(\Mc)}\le\lambda_2$. Thus, for this choice of $g_F$ both property \rf{WG-HF} and inequality \rf{G-LIP} hold.
\par We note that, thanks to Lemma \reff{WF1-NE}, property \rf{TC-BO} holds provided the set $\Mc$ is {\it infinite},
so that in this case we can define $g_F$ by formula \rf{GF-CNT}.
\par Of course, if $\Mc$ is finite, we can not guarantee that property \rf{TC-BO} holds. However, in this case, the following property of the rectangles $\{\Tc^{[1]}_{F,\lambda_1}[x:\lambda_2]:x\in\Mc\}$ maybe useful: if one of the rectangles of this family is a {\it bounded set}, then {\it all rectangles} from this family are bounded as well, i.e., \rf{TC-BO} holds. This property is immediate  from the fact that at {\bf STEP 4} the mapping $\Tc^{[1]}_{F,\lambda_1}[\cdot:\lambda_2]$ is Lipschitz with respect to the Hausdorff distance. See inequality \rf{T-HD}. (Recall that $\lambda=\lambda_2$.)
\smsk
\par {\it(iii)} Suppose that there exists $r>0$ such that, for every $x\in\Mc$, the intersection of the square $Q(O,r)$ and the rectangle
$\Tc^{[1]}_{F,\lambda_1}[x:\lambda_2]$ is non-empty. (For instance, such $r$ exists provided the set $\Mc$ is finite.) In this case, we can define $g_F$ by the formula
$$
g_F(x)=\cent
\left(\Tc^{[1]}_{F,\lambda_1}[x:\lambda_2]^{\cl}\cap Q(O,r)\right),~~~~~x\in\Mc.
$$
\par Clearly,
$$
g_F(x)\in\Tc^{[1]}_{F,\lambda_1}[x:\lambda_2]^{\cl} \subset\Tc_{F,\lambda_1}^{\cl}(x)~~~\text{on}~~\Mc
$$
so that property \rf{WG-HF} holds. The proof of \rf{G-LIP} in this case follows the same scheme as the proof of this inequality for $g_F$ defined by formula \rf{G-CNT}.\rbx}
\end{remark}

\indent\par {\bf 6.2 Projection Algorithms and a solution to the second main problem.}
\addtocontents{toc}{~~~~6.2 Projection Algorithms and a solution to the second main problem.\hfill \thepage\par\VST}

\indent
\msk
\par We are in a position to present a {\it solution to Problem \reff{PR2}}, the second main problem of the paper.
\par Let $\Mf=\MR$ be a pseudometric space, and let $F:\Mc\to\CRT$ be a set-valued mapping.

\begin{definition}\lbl{DF-LR} {\em Let
\bel{LG-SET}
\Lc_{\Rc}(F)=\{\lambda\ge 0:
\RL[x,x':\lambda]\cap \{\RL[y,y':\lambda]+\lambda\,\rho(x,y)\,Q_0\}\ne\emp
~\,\text{for all}~x,x',y,y'\in\Mc\}.
\ee
\par We set
\bel{LM-RT1}
\Lambda_{\Rc}(F)=\inf \Lc_{\Rc}(F)
\ee
provided $\Lc_{\Rc}(F)\ne\emp$, and $\Lambda_{\Rc}(F)=+\infty$ if $\Lc_{\Rc,T}(F)=\emp$.}\rbx
\end{definition}
\par Thus, $\Lambda_{\Rc}(F)$ is the infimum of constants $\lambda\ge 0$ such that condition \rf{R-MLM} of Theorem \reff{CR-LS1} is satisfied. Hence, thanks to the necessity part of Theorem \reff{CR-LS1} and \rf{OP-NF}, we have
\bel{LM-BL}
\Lambda_{\Rc}(F)\le \,|F|_{\Mf}.
\ee
\par Let us prove the following useful property of the constant $\Lambda_{\Rc}(F)$.

\begin{lemma}\lbl{L-ATT} Let $\Mf=\MR$ be a pseudometric space. Suppose that either $F:\Mc\to\KRT$ or $F:\Mc\to\HPL$, see \rf{KRT-DF} and \rf{HPL-DF}. If $\Lambda_{\Rc}(F)<\infty$, then the infimum in the right hand side of definition \rf{LM-RT1} is attained. In other words, in these settings, condition \rf{R-MLM} of Theorem \reff{CR-LS1} is satisfied with $\lambda=\Lambda_{\Rc}(F)$.
\end{lemma}
\par {\it Proof.} Let $T=(x,x',y,y')\subset\Mc$ be an (ordered) quadruple of elements of $\Mc$, and let
$$
\Lc_{\Rc,T}(F)=\{\lambda\ge 0:
\RL[x,x':\lambda]\cap \{\RL[y,y':\lambda]+\lambda\,\rho(x,y)\,Q_0\}\ne\emp\}.
$$
Let
\bel{LMT-R}
\Lambda_{\Rc,T}(F)=\inf \Lc_{\Rc,T}(F)
\ee
provided $\Lc_{\Rc,T}(F)\ne\emp$, and $\Lambda_{\Rc,T}(F)=+\infty$ if $\Lc_{\Rc,T}(F)=\emp$.
\par Clearly, thanks to Definition \reff{DF-LR} and \rf{LMT-R},
\bel{SUP-L}
\Lambda_{\Rc}(F)=\sup \Lambda_{\Rc,T}(F)
\ee
where the supremum is taken over all (ordered) quadruples $T=(x,x',y,y')\subset\Mc$.
\par Let us prove that the infimum in the right hand side of \rf{LMT-R} is attained, i.e., for each (ordered) quadruple $T=(x,x',y,y')$ of elements of $\Mc$, we have
\bel{LM-ATR}
\RL[x,x':\lambda]\cap \{\RL[y,y':\lambda]+\lambda\,\rho(x,y)\,Q_0\}\ne\emp
~~~~\text{for every}~~~\lambda\ge\Lambda_{\Rc,T}(F).
\ee
\par Let us note that if this property is true, then the proof of the lemma is immediate. Indeed, thanks to \rf{SUP-L}, $\Lambda_{\Rc,T}(F)\le\Lambda_{\Rc}(F)$ for every $T$, so that, thanks to \rf{LM-ATR}, for every $x,x',y,y'\in\Mc$, we have
$$
\RL[x,x':\lambda]\cap \{\RL[y,y':\lambda]+\lambda\,\rho(x,y)\,Q_0\}\ne\emp
~~~~\text{with}~~~\lambda=\Lambda_{\Rc}(F)
$$
proving the lemma.
\smsk
\par We turn to the proof of property \rf{LM-ATR}. Thanks to part (i) of Remark \reff{RMK-A}, the quantity
\bel{LMT-WQ}
\Lambda_{\Rc,T}(F)=\inf \lambda
\ee
where the infimum is taken over all constants $\lambda\ge 0$ such that for  every $i=1,2$, there exist points
\bel{ST-0}
A(i)=(a_1(i),a_2(i))\in F(x),~~A'(i)\in F(x'),~~ B(i)=(b_1(i),b_2(i))\in F(y),~~B'(i)\in F(y')
\ee
satisfying the following inequalities:
\bel{ST-1}
\|A(i)-A'(i)\|\le \lambda\rho(x,x'),~~~
\|B(i)-B'(i)\|\le \lambda\rho(y,y')~~\text{and}~~|a_i(i)-b_i(i)|\le \lambda\rho(x,y).
\ee
\par Furthermore, property \rf{LM-ATR} holds provided the infimum in the right hand side of \rf{LMT-WQ} is attained.
\par Let us prove this property. We know that $\Lambda_{\Rc}(F)$ is finite, i.e., $\Lambda_{\Rc}(F)\le\lambda_0$ for some $\lambda_0\ge 0$. Therefore, thanks to \rf{SUP-L}, $\Lambda_{\Rc,T}(F)\le\lambda_0$ as well. Thus, without loss of generality, we can assume that in \rf{LMT-WQ} we have $0\le\lambda\le\lambda_0$. Furthermore, this observation and Definition \reff{DF-LR} imply the existence of a constant $\lambda$ and points $A(i)$, $A'(i)$, $B(i)$, $B'(i)$ satisfying this inequality and constraints \rf{ST-0} and \rf{ST-1}.
\par Our aim is to show that under all these conditions, the infimum in \rf{LMT-WQ} is attained. To this end, let us consider a tuple
$$
\Tc=(\lambda,A(1),A'(1),B(1),B'(1),A(2),A'(2),B(2),B'(2))
$$
where $A(i)=(a_1(i),a_2(i))$, $A'(i)=(a'_1(i),a'_2(i))$, $B(i)=(b_1(i),b_2(i))$, $B'(i)=(b'_1(i),b'_2(i))$, $i=1,2$, are points in $\RT$. Because $\Tc$ depends on $17$ parameters, we will identify this tuple with a point in $\R^{17}$. Then the constrains \rf{ST-0}, \rf{ST-1} and the inequality $0\le\lambda\le\lambda_0$ determine a certain {\it non-empty convex closed subset} $E\subset\R^{17}$. We have to prove that the minimum of the function $G(T)=\lambda$, $T\in E$, is attained on the constraint set $E$.
\par First, let us show this provided $F:\Mc\to\KRT$. In this case $F(x)$, $F(x')$, $F(y)$ and $F(y')$ are {\it compact subsets} of $\RT$ so that the set $E$ is a {\it compact subset of $\R^{17}$}. Therefore, the continuous function $G$ takes the minimum value on the set $E$ proving \rf{LM-ATR} in the case under consideration.
\par Let us prove \rf{LM-ATR} provided $F:\Mc\to\HPL$.
In this case the sets $F(x)$, $F(x')$, $F(y)$ are half-planes. We also recall that $\|\cdot\|$ is the {\it uniform norm} in $\RT$. These observation tell us the set $E$ is determined by {\it the finite number of linear constrains}, i.e., $E$ is a {\it non-empty convex polytope} (not necessarily bounded) in $\R^{17}$.
\par We also note that the objective function $G(T)=\lambda$ is a {\it linear functional bounded from below} on $E$ (because $\lambda\ge 0$). Thus, we have a problem of linear programming on a non-empty convex polytope $E$ in $\R^{17}$ (i.e., with at least one feasible solution in the terminology of the linear programming theory), and with the objective function is bounded from below on the set $E$, the set of all feasible solutions. In this case, there exists an optimal solution to this problem, see, e.g., \cite{MG-2007}, Theorem 4.2.3.
\par Thus, in the both cases, the infimum  in the right hand side of \rf{LMT-WQ} is attained, and the proof of the lemma is complete.\bx
\smsk
\par Applying to $F$ the $\lmv$-Projection Algorithm with an appropriate parameter $\lmv$, we obtain the following solution to Problem \reff{PR2} for the class $\Tf$ of the set-valued mappings satisfying Condition \reff{CND-T}.

\begin{theorem}\lbl{PR2-SL} Let $\Mf=\MR$ be a pseudometric space, and let $F:\Mc\to\CRT$ be a set-valued mapping. Suppose that $\Mf$ and $F$ satisfy Condition \reff{CND-T}. Then the following statements hold:
\smsk
\par (i) There exists a Lipschitz selection of $F$ if and only if $\Lambda_{\Rc}(F)<\infty$.
\smsk
\par (ii) Suppose that $0<\Lambda_{\Rc}(F)<\infty$. Then for every constant $\gamma>1$, the $\lmv$-Projection Algorithm with
$$
\lmv=(\,3\gamma\Lambda_{\Rc}(F),\,\gamma\Lambda_{\Rc}(F)\,)
$$
produces the outcome {\bf ``Success''} and returns the Lipschitz selection $f_{\lmv;F}$ of $F$ with Lipschitz seminorm
$$
\|f_{\lmv;F}\|_{\Lip(\Mc)}\le 5\gamma\,\FM.
$$
\par (iii) If $\Lambda_{\Rc}(F)=0$, then for every $\ve>0$, the $\lmv$-Projection Algorithm with $\lmv=(3\ve,\ve)$ produces the outcome {\bf ``Success''} and returns the Lipschitz selection $f_{\lmv;F}$ of $F$ with $\|f_{\lmv;F}\|_{\Lip(\Mc)}\le 5\ve$.
\smsk
\par (iv) Suppose that either $F:\Mc\to\KRT$ or $F:\Mc\to\HPL$. If $\Lambda_{\Rc}(F)$ is finite then for every $\gamma\ge 1$, the $\lmv$-PA with
$\lmv=(\,3\gamma\Lambda_{\Rc}(F),\,\gamma\Lambda_{\Rc}(F)\,)
$ produces the outcome {\bf ``Success''} and returns the Lipschitz selection $f_{\lmv;F}$ of $F$ with Lipschitz seminorm $\|f_{\lmv;F}\|_{\Lip(\Mc)}\le 5\gamma\,\FM$.
\end{theorem}
\par {\it Proof.} Part (i) of the theorem is immediate from Theorem \reff{CR-LS1}.
\smsk
\par Let us prove parts (ii) and (iii). We set $\lmh=\gamma\Lambda_{\Rc}(F)$ if $\Lambda_{\Rc}(F)\in(0,\infty)$, and $\lmh=\ve$ if $\Lambda_{\Rc}(F)=0$. Because $\gamma>1$ and $\ve>0$, we have $\Lambda_{\Rc}(F)<\lmh<\infty$. Therefore, thanks to Definition \reff{DF-LR}, condition \rf{R-MLM} of Theorem \reff{CR-LS1} is satisfied with $\lambda=\lmh$. Thus
\bel{C-LHH}
\RL[x,x':\lmh]\cap \{\RL[y,y':\lmh]+\lmh\,\rho(x,y)\,Q_0\}\ne\emp
~~~~~\text{for every}~~x,x',y,y'\in\Mc
\ee
proving that $F$ satisfies the hypothesis of Proposition \reff{WR-LM3} with $\tlm=\lambda=\lmh$. Thanks to this proposition, for every $x,x',x'',y,y',y''\in\Mc$, we have:
$$
\Wc_F[x,x',x'':3\tlm]\cap\left\{\Wc_F[y,y',y'':3\tlm]
+\lambda\,\rho(x,y)Q_0\right\}\ne\emp.
$$
\par This property tells us that $F$ satisfies the hypothesis of Theorem \reff{W-CR} with
\bel{LM-TW}
\tlm=3\lmh~~~\text{and}~~~\lambda=\lmh.
\ee
\par Now, let us apply to $F$ the $\lmv=(\lambda_1,\lambda_2)$-Projection Algorithm with $\lambda_1=3\lmh$ and $\lambda_2=\lmh$, and prove that this algorithm produces the outcome {\bf ``Success''}. In fact, we know that $F$ satisfies the hypothesis of Theorem \reff{W-CR}. Lemma \reff{WF1-NE} tells us that in this case the set $F^{[1]}[x:\tlm]\ne\emp$ for every $x\in\Mc$. Furthermore, thanks to \rf{TC-NEMP}, the rectangle $\Tc_{F,\tlm}^{[1]}[x:\lambda]\ne\emp$ on $\Mc$. Thus, condition \rf{SC-P} holds proving that $\lmv$-PA produces the outcome {\bf ``Success''}. See Remark \reff{REM-SC}.
\smsk
\par Parts ($\bigstar A$) and ($\bigstar B$) of Theorem \reff{ALG-T} tell us that in this case, the $\lmv$-Projection Algorithm returns the mapping $f_{\lmv;F}:\Mc\to\RT$ which is a Lipschitz selection of $F$ with the Lipschitz seminorm
\bel{NF-OK}
\|f_{\lambda;F}\|_{\Lip(\Mc)}\le \lambda_1+2\lambda_2
=3\lmh+2\lmh=5\lmh.
\ee
\par If $\Lambda_{\Rc}(F)=0$, then $\lmh=\ve$ so that $\|f_{\lambda;F}\|_{\Lip(\Mc)}\le 5\ve$ proving part (iii)
of the theorem. If $\Lambda_{\Rc}(F)>0$, then
$\|f_{\lambda;F}\|_{\Lip(\Mc)}\le 5\gamma\Lambda_{\Rc}(F)$.
From this inequality and \rf{LM-BL}, we have  $\|f_{\lambda;F}\|_{\Lip(\Mc)}\le 5\gamma\,\FM$, and the proof of part (ii) of the theorem is complete.
\smsk
\par Let us prove part (iv) of the theorem. Thanks to Lemma \reff{L-ATT}, for arbitrary $\gamma\ge 1$, condition \rf{C-LHH} holds with $\lmh=\gamma\Lambda_{\Rc}(F)$. In other words, the set-valued mapping $F$ satisfies the hypothesis of Proposition \reff{WR-LM3} with $\tlm=\lambda=\lmh$. This enables us to repeat the proof
of part (ii) of the theorem given above (for this choice of the constant $\lmh$). This proof leads us to the statement of part (iv) completing the proof of the theorem.\bx
\msk
\par Theorem \reff{PR2-SL} leads us to a solution to Problem \reff{PR2} for pseudometric spaces $\Mf$ and set-valued mappings $F$ satisfying Condition \reff{CND-T}.
For simplicity, let us demonstrate this solution for the case of a {\it finite set $\Mc$} and a set-valued mapping $F:\Mc\to\Tf$ where $\Tf$ is either the family $\KRT$ or the family $\HPL$. (See \rf{KRT-DF} and \rf{HPL-DF}.)
\begin{alg}\lbl{ALG-PR2} (A constructive algorithm for a nearly optimal Lipschitz selection.) {\em Given a finite pseudometric space $\Mf=\MR$, and a set-valued mapping $F:\Mc\to\Tf$, the algorithm produces a nearly optimal Lipschitz selection $f$ of $F$ (the outcome {\bf ``Success''}) or stops (the outcome {\bf ``No go''}). This procedure includes the following two main steps.
\msk
\par {\bf STEP 1.} At this step we compute the constant
$\Lambda_{\Rc}(F)$. If it turns out that $\Lambda_{\Rc}(F)=+\infty$, then a {\it Lipschitz selection of $F$ does not exist}. In this case, the algorithm produces the outcome {\bf ``No go''} and stops.
\par If we determine that $\Lambda_{\Rc}(F)<+\infty$, we calculate this constant up to some constant $\gamma\ge 1$. Thus, at this step the algorithm returns the number $\gamma\,\Lambda_{\Rc}(F)$.
\smsk
\par {\bf STEP 2.} At this step we apply to $F$ the $\lmv=(\lambda_1,\lambda_2)$-Projection Algorithm with

\bel{PRM-A1}
\lambda_1=3\gamma\,\Lambda_{\Rc}(F)~~~~\text{and}
~~~~\lambda_2=\gamma\,\Lambda_{\Rc}(F).
\ee
\par The Projection Algorithm produces the outcome {\bf ``Success''} and returns the mapping $f_{\lmv;F}$ which is a Lipschitz selection of $F$ with Lipschitz seminorm
\bel{A-FN}
\|f_{\lambda;F}\|_{\Lip(\Mc)}\le 5\gamma\,\Lambda_{\Rc}(F).
\ee
\par At this stage, Algorithm \reff{ALG-PR2} produces the outcome {\bf ``Success''} and stops.\bx}
\end{alg}
\par Note that Algorithm \reff{ALG-PR2} is immediate from part (iv) of Theorem \reff{PR2-SL}.

\begin{remark}\lbl{PR2-RM} {\em (i) Clearly, if $\Mf=\MR$ is a {\it metric space} and $\Mc$ is {\it finite} then {\it any selection $f$} of the set-valued mapping $F$ from  Algorithm \reff{ALG-PR2} is Lipschitz, i.e., $\|f\|_{\Lip(\Mc)}<\infty$. Hence, $\FM<\infty$, so that,
thanks to inequality \rf{LM-BL}, $\Lambda_{\Rc}(F)<\infty$ as well. Thus, if $\rho$ is a metric and $\Mc$ is finite, Algorithm \reff{ALG-PR2} always produces the outcome {\bf ``Success''}.
\par However, if $\rho$ is a pseudometric, i.e., $\rho(x,y)=0$ may hold with certain $x\ne y$, in general,
the quantity $\Lambda_{\Rc}(F)$ may take the value $+\infty$ (even for a finite set $\Mc$). For instance, if there exist elements $x,x'\in\Mc$ such that $\rho(x,x')=0$ and $F(x)\cap F(x')=\emp$, then, thanks to
\rf{RL}, $\RL[x,x':\lambda]=\emp$ for every $\lambda\ge 0$.
Therefore, inequality \rf{R-MLM} does not hold for every
$y,y'\in\Mc$ and any $\lambda\ge 0$ proving that $\Lambda_{\Rc}(F)=+\infty$. See Definition \reff{DF-LR}. It is also clear that in this case $F$ does not have a Lipschitz selection (with respect to the pseudometric $\rho$).
\smsk
\par (ii) In general, the {\it precise computation} of the quantity $\Lambda_{\Rc}(F)$ maybe a very difficult technical problem. For this reason, we introduce in  Algorithm \reff{ALG-PR2} the parameter $\gamma\ge 1$
which enables us to calculate $\Lambda_{\Rc}(F)$ up to this parameter. Inequality \rf{A-FN} tells us that in this case
the algorithm construct a selection of $F$ whose Lipschitz
seminorm is bounded by a constant linearly depending on $\gamma$.\rbx}
\end{remark}

\indent\par {\bf 6.3 The constant $\Lambda_{\Rc}(F)$ and other related constants.}
\addtocontents{toc}{~~~~6.3 The constant $\Lambda_{\Rc}(F)$ and other related constants.\hfill \thepage\par\VST}

\indent
\smsk
\par In this section we give several remarks related to efficient algorithms for computing the constant $\Lambda_{\Rc}(F)$. See Definition \reff{DF-LR}.
\par Theorem \reff{CR-LS-RI} motivates us to introduce the following constant.

\begin{definition}\lbl{DF-LR-INT} {\em Let $\Mf=\MR$ be a pseudometric space, and let $F:\Mc\to\CRT$ be a set-valued mapping. Let
\bel{LG-SET-INT}
\Lc_{\Rc}^{(int)}(F)=\{\lambda\ge 0:
\bigcap_{y,y'\in\Mc}
\left\{\,\RL[y,y':\lambda]+\lambda\,\rho(x,y) Q_0\right\} \ne\emp
~\,\text{for all}~x\in\Mc\}.
\ee
\par Let
\bel{LM-RT1-INT}
\Lambda_{\Rc}^{(int)}(F)=\inf \Lc_{\Rc}^{(int)}(F)
\ee
provided $\Lc_{\Rc}^{(int)}(F)\ne\emp$, and let $\Lambda_{\Rc}^{(int)}(F)=+\infty$ if $\Lc_{\Rc}^{(int)}(F)=\emp$.}\rbx
\end{definition}
\par Thus, $\Lambda_{\Rc}^{(int)}(F)$ is the infimum of constants $\lambda\ge 0$ such that condition \rf{RC-CH1} of Theorem \reff{CR-LS-RI} holds.
\begin{lemma}\lbl{TWO-C2} Let $\Mf=\MR$ be a pseudometric space, and let $F:\Mc\to\CRT$ be a set-valued mapping. Suppose that either $\Mc$ is finite or at least one of the sets $F(x)$, $x\in\Mc$, is compact.
\par Then $\Lc_{\Rc}^{(int)}(F)=\Lc_{\Rc}(F)$ and $\Lambda_{\Rc}^{(int)}(F)=\Lambda_{\Rc}(F)$.
\end{lemma}
\par {\it Proof.} Clearly, thanks to \rf{LM-RT1} and \rf{LM-RT1-INT}, it suffices to prove the equality $\Lc_{\Rc}^{(int)}(F)=\Lc_{\Rc}(F)$.
\par If $\lambda\in\Lc_{\Rc}^{(int)}(F)$, then, thanks to \rf{LG-SET-INT}, for every $x\in\Mc$, we have
\bel{DF-LMK}
\bigcap_{y,y'\in\Mc}
\left\{\,\RL[y,y':\lambda]+\lambda\,\rho(x,y) Q_0\right\} \ne\emp.
\ee
\par Clearly, thanks to this property, for every $x,x',y,y'\in\Mc$, we have
\bel{LG-H3}
\RL[x,x':\lambda]\cap \{\RL[y,y':\lambda]+\lambda\,\rho(x,y)\,Q_0\}\ne\emp
\ee
proving that $\lambda\in\Lc_{\Rc}(F)$, see \rf{LG-SET}.
\smsk
\par Conversely, suppose that $\lambda\in\Lc_{\Rc}(F)$, i.e., \rf{LG-H3} holds for all $x,x',y,y'\in\Mc$. Prove that in this case, for every $x\in\Mc$, property  \rf{DF-LMK} holds as well.
\par We set
$$
\Vc=\{\RL[y,y':\lambda]+\lambda\,\rho(x,y)Q_0:y,y'\in\Mc\}.
$$
\par Let us prove that any two members of the family $\Vc$ have a common point. Indeed, thanks to \rf{LG-H3}, given
$y,y',z,z'\in\Mc$, there exist points $a\in\RL[y,y':\lambda]$ and $b\in\RL[z,z':\lambda]$ such that $\|a-b\|\le \lambda\,\rho(y,z)$. Therefore, thanks to the triangle inequality,
$\|a-b\|\le \lambda\,(\rho(y,x)+\rho(x,z))$ so that there exists a point $w\in[a,b]$ such that
$\|a-w\|\le \lambda\,\rho(y,x)$ and $\|b-w\|\le \lambda\,\rho(z,x)$. Hence,
$$
w\in \{\RL[y,y':\lambda]+\lambda\,\rho(x,y)Q_0\}\cap
\{\RL[z,z':\lambda]+\lambda\,\rho(x,z)Q_0\}
$$
proving the required property of the family $\Vc$.
\par From this property and the lemma's hypothesis, it follows that the family $\Vc$ satisfies the hypothesis of Lemma \reff{H-R} (i.e., Helly's theorem for rectangles). Thanks to this lemma, the family $\Vc$ has non-empty intersection proving the required property \rf{DF-LMK}.
\par The proof of the lemma is complete.\bx
\smsk
\par Below we will see how the equality $\Lambda_{\Rc}^{(int)}(F)=\Lambda_{\Rc}(F)$, i.e., the representation of $\Lambda_{\Rc}(F)$ in the form $\Lambda_{\Rc}(F)=\inf \Lc_{\Rc}^{(int)}(F)$, see \rf{LM-RT1-INT}, will help us calculate this constant in an efficient way.
\par But now we introduce one more useful constant directly related to the Finiteness Principle in $\RT$, see Theorem \reff{FP-RT}.

\begin{definition}\lbl{FP-4} {\em Let $\Mf=\MR$ be a pseudometric space, and let $F:\Mc\to\CRT$ be a set-valued mapping. We let $\Lc^{(\Fc\Pc)}(F)$ denote the family of all constants $\lambda\ge 0$ such that for every subset $\Mc'\subset\Mc$ consisting of at most {\it four} points, the restriction $F|_{\Mc'}$ of $F$ to $\Mc'$ has a Lipschitz selection $f_{\Mc'}$ with Lipschitz  seminorm $\|f_{\Mc'}\|_{\Lip(\Mc')}\le \lambda$.
\par Let
\bel{LFP-W}
\Lambda^{(\Fc\Pc)}(F)=\inf \Lc^{(\Fc\Pc)}(F)
\ee
provided $\Lc^{(\Fc\Pc)}(F)\ne\emp$, and let $\Lambda^{(\Fc\Pc)}(F)=+\infty$ if $\Lc^{(\Fc\Pc)}(F)=\emp$.}\rbx
\end{definition}
\par Thus, $\Lambda^{(\Fc\Pc)}(F)$ is the smallest $\lambda\ge 0$ such that the hypothesis of Theorem \reff{FP-RT} holds. Clearly,
\bel{LM-FP-N}
\Lambda^{(\Fc\Pc)}(F)=\sup\, \{|F|_{\Mc'}|_{\Mf}:
M'\subset M,\,\#M'\le 4\}.\hspace*{15mm}\text{See \rf{FM}.}
\ee
\par Furthermore, thanks to Definitions \reff{DF-LR} and \reff{FP-4} and formula \rf{LM-FP-N}, the following inequality
\bel{LM-DT}
\Lambda_{\Rc}(F)\le \Lambda^{(\Fc\Pc)}(F)\le |F|_{\Mf}
\ee
holds provided $F:\Mc\to\CRT$ is an arbitrary set-valued mapping. Moreover, if $\Mc$ and $F$ satisfy Condition \reff{CND-T}, then, thanks to Theorem \reff{CR-LS1} and Theorem \reff{FP-RT}, we have
\bel{NF-MLM2}
|F|_{\Mf}\le\min\,\left\{5\Lambda_{\Rc}(F),
3\Lambda^{(\Fc\Pc)}(F)\right\}.
\ee

\begin{lemma}\lbl{L-FP4} Let $\Mf=\MR$ be a pseudometric space, and let $F:\Mc\to\Tf$ where $\Tf$ is either $\KRT$ or $\HPL$. If $\Lambda_{\Rc}(F)<\infty$, then $\Lambda^{(\Fc\Pc)}(F)\in\Lc^{(\Fc\Pc)}(F)$. Cf. \rf{LFP-W}.
\end{lemma}
\par {\it Proof.} We follow the scheme of the proof of Lemma \reff{L-ATT}, and use equality \rf{LM-FP-N}. We leave the details to the interested reader.\bx

\begin{theorem}\lbl{FP-AL} Let $\Mf=\MR$ be a pseudometric space, and let $F:\Mc\to\CRT$ be a set-valued mapping. Suppose that either $F:\Mc\to\KRT$ or $F:\Mc\to\HPL$.
\par If $\Lambda^{(\Fc\Pc)}(F)<\infty$, then for every $\gamma\ge 1$, the $\lmv$-Projection Algorithm with
$$
\lmv=(\,\gamma\Lambda^{(\Fc\Pc)}(F),
\,\gamma\Lambda^{(\Fc\Pc)}(F)\,)
$$
produces the outcome {\bf ``Success''} and returns the Lipschitz selection $f_{\lmv;F}$ of $F$ with Lipschitz seminorm $\|f_{\lmv;F}\|_{\Lip(\Mc)}\le 3\gamma\,\FM$.
\end{theorem}
\par {\it Proof.} We follow the approach suggested in the proof of part (iv) of Theorem \reff{PR2-SL}. Thanks to Lemma \reff{L-FP4}, $\gamma\Lambda^{(\Fc\Pc)}(F)\in\Lc^{(\Fc\Pc)}(F)$ for every $\gamma\ge 1$. In other words, for every subset $\Mc'\subset\Mc$ consisting of at most four points, the restriction $F|_{\Mc'}$ of $F$ to $\Mc'$ has a Lipschitz selection $f_{\Mc'}$ with Lipschitz  seminorm $\|f_{\Mc'}\|_{\Lip(\Mc')}\le \lambda$ where $\lambda=\gamma\Lambda^{(\Fc\Pc)}(F)$. As we have shown in the sufficiency part of the proof of Theorem \reff{FP-RT},
in this case condition \rf{WNEW} of Theorem \reff{W-CR} is satisfied with $\tlm=\lambda$. Combining this condition with the hypothesis of Theorem \reff{FP-AL}, we conclude that $\Mf$ and $F$ satisfy the hypothesis of Theorem \reff{W-CR} with $\tlm=\lambda=\gamma\Lambda^{(\Fc\Pc)}(F)$ and $\lambda=\gamma\Lambda^{(\Fc\Pc)}(F)$.
\smsk
\par We then repeat the proof of part (ii) of Theorem \reff{PR2-SL} from definition \rf{LM-TW} up to and including inequality \rf{NF-OK}, setting $\tlm=\lambda$ and $\lambda_1=\lambda_2=\lambda$ in this proof. As a result, we show that in this case the $\lmv$-Projection Algorithm with $\lmv=(\lambda,\lambda)$ produces the outcome {\bf ``Success''} and returns the mapping $f_{\lmv;F}:\Mc\to\RT$ which is a Lipschitz selection of $F$ with the Lipschitz seminorm
$$
\|f_{\lambda;F}\|_{\Lip(\Mc)}\le \lambda_1+2\lambda_2
=3\lambda=3\gamma\Lambda^{(\Fc\Pc)}(F).
$$
\par But $\Lambda^{(\Fc\Pc)}(F)\le |F|_{\Mf}$, see \rf{LM-DT}, and the proof of the theorem is complete.\bx

\begin{alg}\lbl{AL-FPR} {\em Theorem \reff{FP-AL} leads us to a new algorithm for a nearly optimal Lipschitz selection of a set-valued mapping. We obtain this algorithm by a slight modification of Algorithm \reff{ALG-PR2}. More specifically, at {\bf STEP 1} of this algorithm we replace the parameter $\Lambda_{\Rc}(F)$ with the parameter $\Lambda^{(\Fc\Pc)}(F)$. At {\bf STEP 2} of Algorithm \reff{ALG-PR2} we set
\bel{PRM-A2}
\lambda_1=\gamma\,\Lambda^{(\Fc\Pc)}(F)~~~~\text{and}
~~~~\lambda_2=\gamma\,\Lambda^{(\Fc\Pc)}(F).
\ee
\par Then, thanks to Theorem \reff{FP-AL}, inequality \rf{A-FN} transforms into the following one:
$$
\|f_{\lambda;F}\|_{\Lip(\Mc)}\le 3\gamma\,\Lambda_{\Rc}(F).
$$
\par At this stage, Algorithm \reff{AL-FPR} produces the outcome {\bf ``Success''} and stops.\bx}
\end{alg}
\smsk
\par Let us make two remarks related to some preliminary estimates of the computational efficiency of Algorithm \reff{ALG-PR2} and Algorithm \reff{AL-FPR}.
\begin{remark}\lbl{TWO-ALG} {\em  We note that at the second step of Algorithms \reff{ALG-PR2} and \reff{AL-FPR} the $\lmv$-Projection Algorithm is applied with some parameters $\lmv=(\lambda_1,\lambda_2)$, see \rf{PRM-A1} and  \rf{PRM-A2}. In the forthcoming paper \cite{S-2023} we will show that if $\Mc$ is an $N$-element pseudometric space and $F$ satisfies certain rather mild geometric conditions, the running time of the $\lmv$-Projection Algorithm is $O(N^3)$. In particular, this is true provided
each $F(x)$ is a half-plane, i.e., $F:\Mc\to\HPL$.
\smsk
\par The most difficult part of Algorithms \reff{ALG-PR2} and \reff{AL-FPR} is {\bf STEP 1}, i.e., calculation (up to some constant $\gamma\ge 1$) the values of the constants $\Lambda_{\Rc}(F)$ and $\Lambda^{(\Fc\Pc)}(F)$ respectively. Although the definitions of these constants are given in explicit geometric terms, their calculation may require a huge number of computer operations. (Even if every set $F(x)$, $x\in\Mc$, is a disc in $\RT$, it is not clear how to calculate $\Lambda_{\Rc}(F)$ and $\Lambda^{(\Fc\Pc)}(F)$ (even up to an absolute  constant)).
\smsk
\par Let us consider the simplest non-trivial case of {\it a set-valued mapping from $\Mc$ (with $\#\Mc=N$) to the family $\HPL$ of all half-planes in $\RT$}.
\par It is not difficult to compute $\Lambda_{\Rc}(F)$ with
$O(N^4)$ running time. Indeed, to calculate the constant $\Lambda_{\Rc}(F)$ we have to consider all possible (ordered) quadruples $T=(x,x',y,y')\subset\Mc$ and solve the corresponding linear programming problem determined by \rf{LMT-WQ}, \rf{ST-0} and \rf{ST-1}. For each $T$, solving this problem will take $O(1)$ running time. Since the number of such quadruples $T$ is $O(N^4)$, the total running time is also $O(N^4)$. The same estimate $O(N^4)$ is obtained for the running time of computing the constant
$\Lambda^{(\Fc\Pc)}(F)$.
\par Thus, for both Algorithm \reff{ALG-PR2} and Algorithm \reff{AL-FPR} this approach provides the running time}
$$
O(N^4)\,~\text{(at {\bf STEP 1})}+O(N^3)\,
~\text{(at {\bf STEP 2})}=O(N^4). \text{\rbx}
$$
\end{remark}

\begin{remark}\lbl{N3-E} {\em For the same case of a set-valued mapping $F:\Mc\to\HPL$ defined on a finite set $\Mc$ with $\#\Mc=N$, there exists another algorithm for calculating the constant $\Lambda_{\Rc}(F)$ with $O(N^3)$ running time. Let us briefly explain the main idea of this algorithm.
\par Definition \reff{DF-LR-INT} and Lemma \reff{TWO-C2} tell us that, in the case under the consideration,  $\Lambda_{\Rc}^{(int)}(F)=\Lambda_{\Rc}(F)$. Thus, $\Lambda_{\Rc}^{(int)}(F)=\inf \lambda$ where $\lambda$ runs over all non-negative numbers such that for every $x\in\Mc$ there exists a point $u\in\RT$ such that
\bel{U-ST}
u\in\RL[y,y':\lambda]+\lambda\,\rho(x,y) Q_0
\ee
for all $y,y'\in\Mc$.
\par Fix $x\in\Mc$ and consider the family $\Vc_x$ of all points $v=(u,\lambda)\in\R^3$ with $\lambda\ge 0$ such that property \rf{U-ST} holds. The reader can easily see that
{\it $\Vc_x$ is intersection of at most four half-spaces in $\R^3$}. Each of these half-spaces depends only on parameters determined the half-planes $F(y)$ and $F(y')$.
Hence, the total number of linear constraints determining the set $\Vc_x\subset\R^3$ is bounded by $O(N^2)$.
\par Let $\Lambda_{\Rc,x}(F)=\{\inf \lambda: (u,\lambda)\in\Vc_x\}$. Then,
\bel{L-MX}
\Lambda_{\Rc}(F)=\max_{x\in\Mc} \Lambda_{\Rc,x}(F).
\ee
\par Let us see that there exists an algorithm which for every $x\in\Mc$, computes the quantity $\Lambda_{\Rc,x}(F)$ with the running time $O(N^2)$. Indeed, we know that the problem
\bel{VX}
\text{\it Find}~~~\inf\lambda~~~
\text{\it under the condition}~~~(u,\lambda)\in \Vc_x,
\ee
is a {\it linear program in three variables with $O(N^2)$ linear constraints}. The following theorem is the  classical results on low-dimensional linear programming by N. Megiddo \cite{M-1983} and M. E. Dyer \cite{Dy-1984}.
\begin{theorem}\lbl{MD-L} A linear program in three variables and $m$ constraints can be solved in $O(m)$ time.
\end{theorem}
\par This theorem implies the following useful corollary.
\begin{corollary}\lbl{CR-HD} (i) There exists an algorithm which given a convex polygon $G\subset\RT$ determined by $N$ linear constraints, constructs its rectangular hull $\Hc[G]$ in $O(N)$ running time.
\par (ii) There exists an algorithm which given a point $A\in\RT$ and a convex polygon $G\subset\RT$ determined by $N$ linear constraints, calculates the distance from $A$ to $G$ using at most $O(N)$ computer operations.
\end{corollary}
\par {\it Proof.} These results are well known in the theory of geometric algorithms. Each of them can be readily reduced to a certain linear program in three variables and $CN$ constraints where $C$ is an absolute constant. Applying Theorem \reff{MD-L} to the corresponding linear program, we obtain the statements of the corollary. We leave the details to the interested reader.\bx
\smsk
\par Thus, for every $x\in\Mc$ the problem \rf{VX} can be solved in $O(N^2)$ time so that, thanks to \rf{L-MX}, the constant $\Lambda_{\Rc}(F)$ can be calculated in $N\cdot O(N^2)=O(N^3)$ time.
\smsk
\par Combining this algorithm for $\Lambda_{\Rc}(F)$ with Algorithm \reff{ALG-PR2}, we obtain  a constructive algorithm which, given a set-valued mapping $F:\Mc\to\HPL$ assigns a Lipschitz selection $f$ of $F$ with $\|f\|_{\Lip(\Mc)}\le 5\,|F|_{\Mf}$. The running time of this algorithm is}
$$
O(N^3)\,~\text{({\bf STEP 1})}+O(N^3)\,
~\text{({\bf STEP 2})}=O(N^3).~~~~~~\text{\rbx}
$$
\end{remark}
\smsk

\indent\par {\bf 6.4 Lipschitz selections of polygon-set valued mappings.}
\addtocontents{toc}{~~~~6.4 Lipschitz selections of polygon-set valued mappings.\hfill \thepage\par\VST}

\indent
\smsk
\par In this section, we provide some remarks related to a recent paper by C. Fefferman and B. Peguero\-les \cite{FP-2019}. Let $D$ and $L$ be positive integers and let ${\mathscr P}_L(\R^D)$ be the family of all compact convex polytope in $\RD$ defined by at most $L$ linear constraints. In \cite{FP-2019} C. Fefferman and B. Pegueroles solved a selection problem for set-valued mappings $F$ from a set $E\subset\RN$ into ${\mathscr P}_L(\R^D)$ {\it with a slight enlarging of the targets} $F(x)$, $x\in E$. (See also \cite[Chapter 7.2]{FI-2020}.) Let us recall a particular case of this result related to Lipschitz selections. Let $E$ be a finite subset of $\RN$ with $\#E=N$, and let $F:E\to{\mathscr P}_L(\R^D)$ be a set-valued mapping. Given $\tau>0$ and $x\in E$, we let $(1+\tau)\lozenge F(x)$ denote the convex set obtained by dilating $F(x)$ about its center of mass by a factor $(1+\tau)$.
\begin{theorem}\lbl{FP-2019} Let $M,\tau>0$ be given. The algorithm described in \cite{FP-2019} produces one of the following two outcomes, using at most $C_1\,N \log N$ computer operations and $C_1\,N$ units of computer memory.
\smsk
\par {\bf Outcome 1 (``No Go''):} We guarantee that there does not exist $f\in\Lip(E,\RD)$ with the Lipschitz seminorm $\le M$ such that $f(x)\in F(x)$ for all $x\in E$. \smsk
\par {\bf Outcome 2 (``Success''):} The algorithm produces a mapping $f:E\to\RD$ with the Lipschitz seminorm at most
$C_2\,M$ satisfying $f(x)\in(1 + \tau)\lozenge F(x)$ for all $x\in E$.
\par Here, $C_1>0$ depends only on $\tau$, $L$, $n$ and $D$, and $C_2>0$ depends only on $n$ and $D$.
\end{theorem}
\par This theorem motivates us to consider Problem \reff{PR2} for the family $\Tf=\PL$ of all non-empty {\it convex polygons in $\RT$} each defined by at most $L$ linear constraints. Let us note that this problem is a particular case of the same problem but for the family $\Tf=\HPL$ of all half-planes in $\RT$. Indeed, let $\Mf=\MR$ be a pseudometric space, and let $F:\Mc\to\PL$ be a set-valued mapping. We know that each polygon $F(x)$, $x\in\Mc$, can be represented as an intersection of at most $L$ half-planes. We denote this family of half-planes by $\Hcr_F(x)$; thus
\bel{FH-CR}
F(x)=\cap \{H:H\in\Hcr_F(x)\}~~~\text{for every}~~~ x\in\Mc.
\ee
\par We introduce a new pseudometric space $\tMf=(\tMc,\trh)$ whose elements are all couples $u=(x,H)$ where $x\in\Mc$ and $H\in \Hcr_F(x)$. We define on $\tMc$ The pseudometric $\trh$ on $\tMc$ is defined by
$$
\trh(u,u')=\rho(x,x')~~~\text{for every} ~~~u=(x,H),~u'=(x',H')\in\tMc.
$$
\par Finally, we define a new set-valued mapping $\tF:\tMc\to\HPL$ by letting
\bel{TF-4}
\tF((x,H))=H~~~\text{provided}~~~x\in \Mc~~\text{and}~~ H\in\Hcr_F(x).
\ee
\par In particular, $\trh((x,H),(x,H'))=0$ for every $H,H'\in\Hcr_F(x)$. This observation implies the following simple claim.
\begin{claim}\lbl{CL-PL} (i) If a mapping $f:\Mc\to\RT$ is a Lipschitz selection of $F$ then the mapping $\tf:\tMc\to\RT$ defined by
\bel{FW-FR}
\tf((x,H))=f(x),~~~~(x,H)\in\tMc,
\ee
is a Lipschitz (with respect to $\trh$) selection of $\tF$, and the following equality
\bel{FW-NM}
\|\tf\|_{\Lip(\tMc,\trh)}=\|f\|_{\Lip(\Mc,\rho)}
\ee
holds.
\par (ii) Conversely, if a mapping $\tf:\tMc\to\RT$ is a Lipschitz (with respect to $\trh$) selection of $\tF$, then there exists a unique mapping $f:\Mc\to\RT$ satisfying   \rf{FW-FR}. Furthermore, $f$ is a Lipschitz selection of $F$, and equality \rf{FW-NM} holds.
\end{claim}
\par {\it Proof.} Part (i) of the claim is obvious. Let us prove part (ii). Because $\tf$ is Lipschitz (with respect to $\trh$), for every $x\in\Mc$ and every $H,H'\in\Hcr_F(x)$, we have
\bel{DF-FM3}
\|\tf((x,H))-\tf((x,H'))\|\le \|\tf\|_{\Lip(\tMc,\trh)}\,\rho(x,x)=0.
\ee
\par This enables us to define a mapping $f:\Mc\to\RT$ by letting
\bel{FW-FR-RV}
f(x)=\tf((x,H))~~\text{where}~~H\in\Hcr_F(x)~~\text{is arbitrary}.
\ee
Thanks to \rf{DF-FM3}, $f$ is well defined (i.e., $\tf((x,H))$ depends only on $x$ and independent of $H$ provided $H\in\Hcr_F(x)$). Furthermore, $f$ is a selection of $F$ on $\Mc$. Indeed, because $\tf$ which is a Lipschitz selection of $\tF$, the point $\tf((x,H))\in \tF((x,H))$ for every $H\in\Hcr_F(x)$. Hence, thanks to \rf{FW-FR-RV},
$f(x)\in \tF((x,H))$ for all $H\in\Hcr_F(x)$ proving that
$$
f(x)\in\cap \{H:H\in\Hcr_F(x)\}=F(x).~~~~\text{See \rf{FH-CR}.}
$$
\par Finally, thanks to \rf{FW-FR-RV} and part (i) of the claim, the equality \rf{FW-NM}holds.\bx
\msk
\par Applying the Projection Algorithm \reff{PA-DS} to the pseudometric space $\tMf=(\tMc,\trh)$ and the set-valued mapping $\tF:\tMc\to\HPL$, we obtain the following theorem.
\begin{theorem}\lbl{PA-PG} Let a positive integer $L$ and a constant $M>0$ be given. Let $\Mf=\MR$ be a finite pseudometric space with $\#\Mc=N$, and let $F:\Mc\to\PL$ be a set-valued mapping.
\par The $(M,M)$-Projection Algorithm produces one of the following two outcomes, using $O(L^3 N^3)$ computer operations.
\smsk
\par {\sc Outcome 1} {\bf (``No Go''):} We guarantee that there does not exist $f\in\Lip(\Mc)$ with the Lipschitz seminorm $\|f\|_{\Lip(\Mc)}\le M$ such that $f(x)\in F(x)$ for all $x\in\Mc$.
\smsk
\par {\sc Outcome 2} {\bf (``Success''):} The algorithm produces a mapping $f:\Mc\to\RT$ with the Lipschitz seminorm $\|f\|_{\Lip(\Mc)}\le 3M$ satisfying $f(x)\in F(x)$ for all $x\in\Mc$.
\end{theorem}
\par {\it Proof.} We set $\lmv=(\lambda_1,\lambda_2)$ where $\lambda_1=\lambda_2=M$. Then we apply Theorem \reff{ALG-T} with this $\lmv$ to $\tMf=(\tMc,\trh)$ and $\tF:\tMc\to\HPL$. This theorem tells us that the $\lmv$-Projection Algorithm produces the following two outcomes:
\par {\sc Outcome 1} {\bf ``No go''}. In this case we guarantee that there does not exist a Lipschitz selection of $\tF$ with Lipschitz seminorm at most $\min\{\lambda_1,\lambda_2\}=M$;
\par {\sc Outcome 2} {\bf ``Success''}. In this case, the $\lmv$-PA returns a mapping $\tf:\tMc\to\RT$ which is a Lipschitz selection of $\tF$ with
$$
\|\tf\|_{\Lip(\tMc,\trh)}\le \lambda_1+2\lambda_2=3M.
$$
\par Let us see that {\sc Outcome 1} and {\sc Outcome 2} provide the statements given in the formulation of Theorem
\reff{FP-2019}. Indeed, suppose that in case of {\sc Outcome 1} there exists a Lipschitz selection $f$ of $F$ with $\|f\|_{\Lip(\Mc)}\le M$. Then the mapping
$\tf:\tMc\to\RT$ defined by formula \rf{FW-FR} is a Lipschitz selection of $\tF$ on $\tMc$ with
$$
\|\tf\|_{\Lip(\tMc,\trh)}=\|f\|_{\Lip(\Mc,\rho)}\le M,
$$
a contradiction.
\par We turn to {\sc Outcome 2}. We know that the mapping $\tf:\tMc\to\RT$ is Lipschitz with respect to the pseudometric $\trh$. Therefore, thanks to part (ii) of Claim \reff{CL-PL}, there exists a unique mapping $f:\Mc\to\RT$ satisfying \rf{FW-FR} such that $f$ is a Lipschitz selection of $F$. Furthermore, in this case  equality \rf{FW-NM} holds proving that
$$
\|f\|_{\Lip(\Mc)}=\|\tf\|_{\Lip(\tMc,\trh)}\le 3M.
$$
\par Thus, both the statement of {\sc Outcome 1} and the statement of {\sc Outcome 2} hold.
\par Finally, let us estimate the running time of the $\lmv$-Projection Algorithm with $\lmv=(M,M)$ which produces {\sc Outcome 1} and {\sc Outcome 2}. As we have noted at the beginning of Remark \reff{TWO-ALG}, the running time of the $\lmv$-Projection Algorithm is $O(N^3)$ provided $F:\Mc\to\HPL$ and $\#\Mc=N$. (As we have noted in this remark, this result will be proved in \cite{S-2023}.) We also know that the pseudometric space $\tMf=(\tMc,\trh)$ contains at most $L\cdot(\#\Mc)$ elements. Therefore, the $(M,M)$-Projection Algorithm produces {\sc Outcome 1} and {\sc Outcome 2}, using at most $O(L^3 N^3)$ computer operations.
\par The proof of the theorem is complete.\bx
\msk
\par Comparing the results of Theorem \reff{FP-2019} and Theorem \reff{PA-PG}, we ask the following question: Let $\Mc=E$ where $E$ is a finite subset of $\RN$ with $\#E=N$,
and let $\rho$ be the Euclidean metric in $\RN$. {\it Does there exist an algorithm which produces {\sc Outcome 1} and {\sc Outcome 2} in Theorem \reff{PA-PG} and uses $O(N\log N)$ computer operations as in Theorem \reff{FP-2019}?} We will concern this problem in paper \cite{S-2023}.
\smsk
\par The last result of this section is the following theorem.
\begin{theorem} Let $L$ be a positive integer. Let $\Mf=\MR$ be a finite \underline{metric} space with $\#\Mc=N$.
\par There exists an algorithm with $O(L^3 N^3)$ running time, which, given $F:\Mc\to\PL$, produces its Lipschitz selection with Lipschitz seminorm at most $5|F|_{\Mf}$.
\end{theorem}
\par {\it Proof.} Let $\tF:\tMc\to\HPL$ be the set-valued mapping defined by \rf{TF-4}. Let $f:\Mc\to\RT$ be a selection of $F$. Because $\Mf=\MR$ is a finite {\it metric} space, $f$ is Lipschitz. Part (i) of Claim \reff{CL-PL} tells us that the mapping
$\tf:\tMc\to\RT$ defined by formula \rf{FW-FR}
is a Lipschitz (with respect to $\trh$) selection of $\tF$.
Hence, $|\tF|_{\tMf}<\infty$ so that, thanks to \rf{LM-BL},
$\Lambda_{\Rc}(\tF)<\infty$.
\smsk
\par Part (iv) of Theorem \reff{PR2-SL} tells us that in this case the $\lmv$-Projection Algorithm with
$$
\lmv=(\,3\Lambda_{\Rc}(\tF),\Lambda_{\Rc}(\tF)\,)
$$
produces the outcome {\bf ``Success''} and returns the Lipschitz selection $\tf=f_{\lmv;\tF}$ of $\tF$ with Lipschitz seminorm $\|\tf\|_{\Lip(\tMc,\trh)}\le 5|\tF|_{\tMf}$.
\par In Remark \reff{N3-E} we have shown the main ideas of an algorithm which produces the mapping $\tf$ in $O((\#\tMc)^3)=O(L^3 N^3)$ running time. (A detailed description and justification of this algorithm based on the results of the works \cite{M-1983} and \cite{Dy-1984} will be done in \cite{S-2023}.)
\par We note that, thanks to part (ii) of Claim \reff{CL-PL}, there exists a (unique) mapping $f:\Mc\to\RT$ satisfying \rf{FW-FR}. This mapping is a Lipschitz selection of $F$ with  $\|f\|_{\Lip(\Mc)}=\|\tf\|_{\Lip(\tMc,\trh)}$.
\par Hence, $\|f\|_{\Lip(\Mc)}\le 5|\tF|_{\tMf}$. But, thanks to part (i) of Claim \reff{CL-PL}, $|\tF|_{\tMf}\le |F|_{\Mf}$ proving the required inequality $\|f\|_{\Lip(\Mc)}\le 5|F|_{\Mf}$. It is also clear that,
thanks to formula \rf{FW-FR}, we can construct the mapping $f$ using the same number of computer operation as when constructing the mapping $\tf$, i.e., at most
$O(L^3 N^3)$.
\par The proof of the theorem is complete.\bx

\SECT{7. Lipschitz selections and iterations of balanced refinements.}{7}
\addtocontents{toc}{7. Lipschitz selections and iterations of balanced refinements.\hfill \thepage\par\VST}

\indent\par {\bf 7.1 The Stabilization Principle for balanced refinements of set-valued mappings.}
\addtocontents{toc}{~~~~7.1 The Stabilization Principle for balanced refinements of set-valued mappings.\hfill \thepage\par\VST}

\indent
\smsk
\par In Section 6, we studied a number of efficient algorithms which provide a solution to Problem \reff{PR2}. These algorithms are based on the Projection Algorithm introduced in Section 6.1.
\par In the next section, we present another approach to Problem \reff{PR2} based on the so-called {\it Iterative Algorithm for Lipschitz selections}. This algorithm relies on the results of a recent author's paper \cite{S-2022}. More specifically, the Iterative Algorithm is a new constructive and nearly optimal algorithm for Lipschitz selections based on an interesting property of successive balanced refinements of set-valued mappings. We refer to this property as {\it the Stabilization Principle for balanced refinements}.
\smsk
\par Let us formulate the Stabilization Principle for the special case of the space $X=\LTI$. Given constants $\lambda_1,\lambda_2,\lambda_3\ge 0,$ and a set-valued mapping $F:\Mc\to\CRT$, we introduce the following mappings:
\bel{F1-X}
F^{[1]}[x:\lambda_1]=
\bigcap_{z\in\Mc}
\big\{F(z)+\lambda_1\,\rho(x,z)\,Q_0\big\},~~~~x\in\Mc,
\ee
\bel{F2-X}
F^{[2]}[x:\lambda_1,\lambda_2]=
\bigcap_{z\in\Mc}
\big\{F^{[1]}[z:\lambda_1]+
\lambda_2\,\rho(x,z)\,Q_0\big\},~~~~x\in\Mc.
\ee
and
\bel{F3-X}
F^{[3]}[x:\lambda_1,\lambda_2,\lambda_3]=
\bigcap_{z\in\Mc}
\big\{F^{[2]}[z:\lambda_1]+
\lambda_3\,\rho(x,z)\,Q_0\big\},~~~~x\in\Mc.
\ee
\par Thus, the mapping $F^{[1]}[\cdot:\lambda_1]$ is the $\lambda_1$-balanced refinement of $F$. See \rf{F-1}. We refer to the mapping $F^{[2]}[\cdot:\lambda_1,\lambda_2]:\Mc\to \CRT\cup\{\emp\}$ as {\it the second order $(\lambda_1,\lambda_2)$-balanced refinement of $F$}. Accordingly, we refer to the mapping $F^{[3]}[\cdot:\lambda_1,\lambda_2,\lambda_3]:\Mc\to \CRT\cup\{\emp\}$ as {\it the third order $(\lambda_1,\lambda_2,\lambda_3)$-balanced refinement of $F$}.
\par Clearly,
\bel{IN-S}
F^{[3]}[x:\lambda_1,\lambda_2,\lambda_3]
\subset F^{[2]}[x:\lambda_1,\lambda_2]
\subset F^{[1]}[x:\lambda_1]
\ee
for every $\lambda_1,\lambda_2,\lambda_3\ge 0$ and every $x\in\Mc$.
\begin{theorem}\lbl{SPR-3}(The Stabilization Principle for $\LTI$) Let $\Mf=\MR$ be a pseudometric space, and let $\lambda\ge 0$. Let $F:\Mc\to\CRT$ be a set-valued mapping such that for every subset $\Mc'\subset\Mc$ with $\#\Mc'\le 4$, the restriction $F|_{\Mc'}$ of $F$ to $\Mc'$ has a Lipschitz selection with Lipschitz seminorm at most $\lambda$. Suppose that either $\Mc$ is finite or $F:\Mc\to\KRT$.
\par Then
\bel{F2-NE}
F^{[2]}[x:\lambda,3\lambda]\ne\emp~~~~\text{for every}~~~~x\in\Mc.
\ee
Furthermore,
\bel{SP-3}
F^{[3]}[x:\lambda,3\lambda,15\lambda]=
F^{[2]}[x:\lambda,3\lambda]~~~~\text{for every}~~~~ x\in\Mc.
\ee
\end{theorem}
\par In particular, Theorem \reff{SPR-3} implies the following property: the sequence of the successive refinements defined by
\bel{SF-RF}
F^{[k+1]}[x:\lambda_1,...,\lambda_{k+1}]=
\bigcap_{z\in\Mc}
\big\{F^{[k]}[z:\lambda_1,...,\lambda_k]+
\lambda_{k+1}\rho(x,z)B_X\big\},~~~~x\in\Mc,
\ee
stabilizes at the third step of this procedure provided $F$ satisfies the hypothesis of this theorem and $\lambda_1=\lambda$, $\lambda_2=3\lambda$ and
\bel{LM-K1}
\lambda_k=\lambda_3=15\lambda~~~\text{for all}~~~k\ge 3.
\ee
In other words, $F^{[k]}=F^{[2]}$ on $\Mc$ for every $k=3,4,...\,$. Indeed, if $F^{[k]}=F^{[2]}$ for some $k\ge 3$, then, thanks to \rf{SF-RF}, \rf{LM-K1} and \rf{SP-3}, for every $x\in\Mc$, we have
$$
F^{[k+1]}[x:\lambda_1,...,\lambda_{k+1}]=
\bigcap_{z\in\Mc}
\big\{F^{[2]}[z:\lambda_1,\lambda_2]+
\lambda_{3}\rho(x,z)B_X\big\}= F^{[3]}[x:\lambda_1,\lambda_2,\lambda_3]=
F^{[2]}[x:\lambda_1,\lambda_2].
$$
\par Also, let us note that, thanks to \rf{F1-X}, \rf{F2-X}, \rf{F3-X} and \rf{IN-S}, equality
\rf{SP-3} is equivalent to the following inequality:
\bel{SP-3D}
\dhf(F^{[2]}[x:\lambda,3\lambda],
F^{[2]}[y:\lambda,3\lambda])\le 15\lambda\,\rho(x,y)
~~~~\text{for every}~~~~ x,y\in\Mc.
\ee
(Recall that the sign $\dhf$ denotes the Hausdorff distance between sets. See \rf{HD-DF}.)
\begin{remark} {\em Theorem \reff{SPR-3} is a slight generalization of the main result of the work \cite{S-2022}, Theorem 1.9, for the case of the space $X=\LTI$. More specifically, in this theorem we give a proof of the Stabilization Principle only for set-valued mappings from $\Mc$ into $\KRT$.
\par However, obvious changes to this proof  show that this principle also holds for an arbitrary {\it finite} pseudometric space $\Mf=\MR$ and arbitrary set-valued mapping $F:\Mc\to\CRT$. These changes are mainly related to the classical Helly intersection theorem in $\RT$. It known that this theorem is true both for arbitrary {\it families of convex compacts} (this version was used in \cite{S-2022}) and arbitrary {\it finite collections of convex sets} in $\RT$.
\par This enables us to add to the hypotheses of Lemma 3.4 in \cite{S-2022} the case of a {\it finite} collection of closed sets. Using this version of this lemma in the proof given in \cite{S-2022}, we obtain the required generalization of the Stabilization Principle to the case of a finite set $\Mc$ and a mapping $F$ from $\Mc$ into $\CRT$.\rbx}
\end{remark}

\bsk
\indent\par {\bf 7.2 The Iterative Algorithm for set-valued mappings.}
\addtocontents{toc}{~~~~7.2 The Iterative Algorithm for set-valued mappings.\hfill \thepage\par\VST}
\indent
\bsk
\par Let $\Mf=\MR$ be a {\it finite} pseudometric space. In this section we describe the Iterative Algorithm for Lipschitz selections. This geometrical algorithm relies on Theorem \reff{SPR-3}.
\par An important parameter of the Iterative Algorithm is the constant $\Lambda^{(\Fc\Pc)}(F)$ introduced in Section 6.3, see Definition \reff{FP-4}.
\msk
\par Let $\Tf$ be either the family $\KRT$ or the family $\HPL$. (See \rf{KRT-DF} and \rf{HPL-DF}.)
\begin{alg}\lbl{ALG-IA} (The Iterative Algorithm for nearly optimal Lipschitz selections.)
\par {\em Given a finite pseudometric space $\Mf=\MR$ and a set-valued mapping $F:\Mc\to\Tf$, the Iterative Algorithm produces a nearly optimal Lipschitz selection of $F$ (the outcome {\bf ``Success''}) or stops (the outcome {\bf ``No go''}).
\msk
\par This procedure includes the following three steps.
\msk
\par {\bf STEP 1.} At this step we compute the constant
$\Lambda^{(\Fc\Pc)}(F)$. If it turns out that
$$
\Lambda^{(\Fc\Pc)}(F)=+\infty,
$$
the algorithm produces the outcome {\bf ``No go''} and stops. In this case, we guarantee that $F$ does not have a Lipschitz selection.
\msk
\par {\bf STEP 2.} At this and the next steps, we assume that the constant  $\Lambda^{(\Fc\Pc)}(F)<\infty$.
\par\noindent We fix a constant
$$
\lambda\in[\Lambda^{(\Fc\Pc)}(F),+\infty),
$$
and, for all $x\in\Mc$, construct the $\lambda$-balanced refinement of $F$,
\smsk
\bel{LB-1}
F^{[1]}[x:\lambda]=
\bigcap_{z\in\Mc}
\big\{F(z)+\lambda\,\rho(x,z)Q_0\big\}.
\ee
Then, for every $x\in\Mc$, we construct the second order $(\lambda,3\lambda)$-balanced refinement of $F$, i.e., the set
\smsk
\bel{LB-2}
F^{[2]}[x:\lambda,3\lambda]=
\bigcap_{z\in\Mc}
\big\{F^{[1]}[z:\lambda]+3\lambda\rho(x,z)Q_0\big\}.
\ee
\par We define a set $\Fcr[x:\lambda]$ by letting
\bel{FR-2}
\Fcr[x:\lambda]=F^{[2]}[x:\lambda,3\lambda]\cap Q(0,2r_x) \ee
where
\bel{FR-21}
r_x=\dist\left(0,F^{[2]}[x:\lambda,3\lambda]\right).
\ee
\smsk
\par {\bf STEP 3.} Finally, we define a mapping $f^{[\lambda;F]}:\Mc\to\RT$ by the formula
\bel{F-I2}
f^{[\lambda;F]}(x)=\cent(\,\Pi_{\lambda,F}(x)), ~~~~~~x\in\Mc,
\ee
where
\bel{PI-H}
\Pi_{\lambda,F}(x)=\Hc[\Fcr[x:\lambda]].
\ee
\par Thus, $\Pi_{\lambda,F}(x)$ is the rectangular hull of the set $\Fcr[x:\lambda]$. See \rf{HRS}.
We also recall that $\cent(S)$ denotes the center of a centrally symmetric set $S\subset\RT$.
\smsk
\par See Fig. 19.
\newpage
\begin{figure}[h!]
\hspace{22mm}
\includegraphics[scale=0.34]{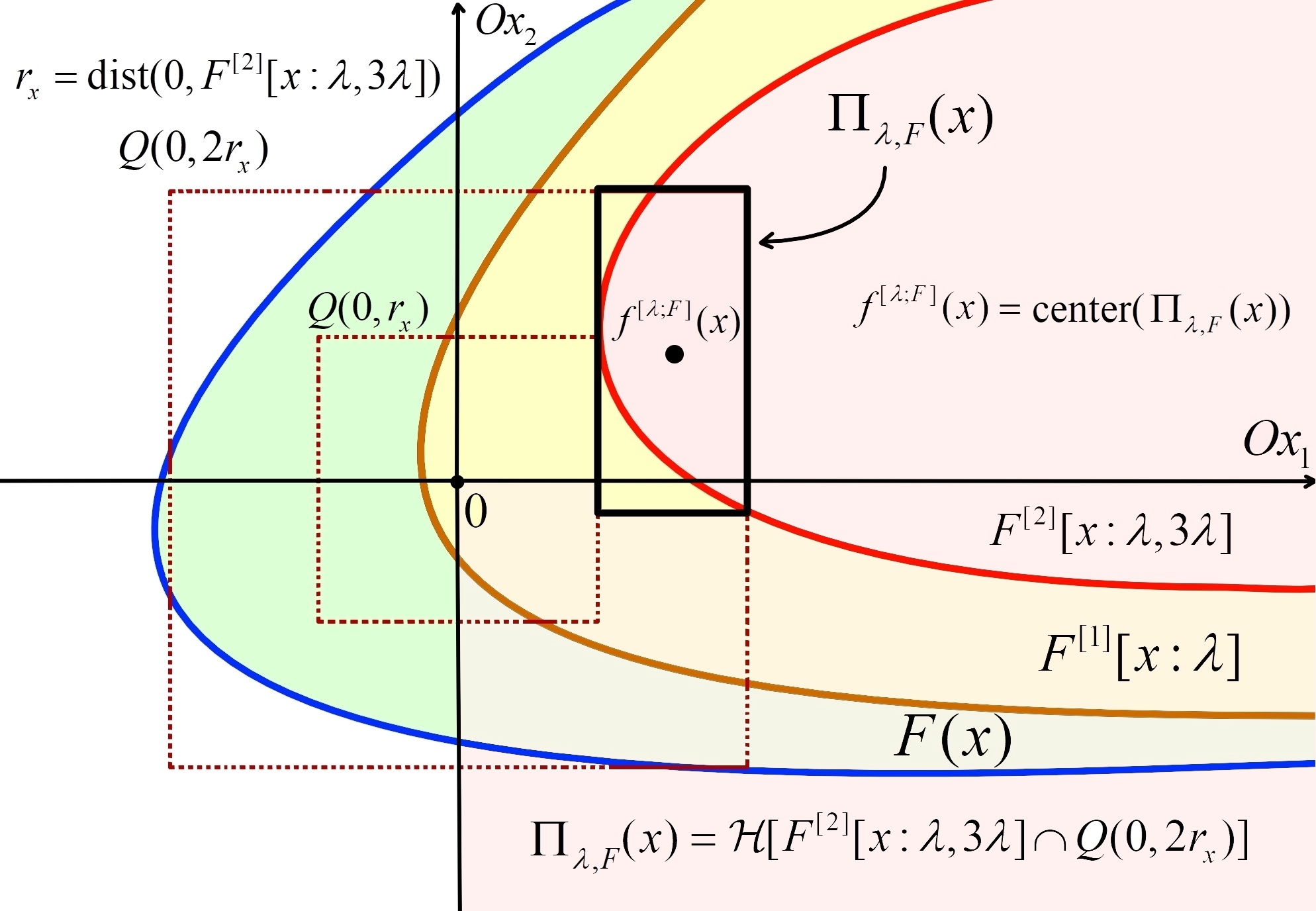}
\caption{The $\LF$-Iterative Algorithm.}
\end{figure}

\par At this stage, Algorithm \reff{ALG-IA} produces the outcome {\bf ``Success''} and stops. \bx}
\end{alg}
\smsk
\par We refer to Algorithm \reff{ALG-IA} as the $(\lambda;F)$-Iterative Algorithm.

\begin{remark} {\em Let if $\Mf=\MR$ be a {\it finite metric space} and let $F:\Mc\to\CRT$ be a set-valued mapping. Clearly, in this case {\it any selection $f$} of $F$ is Lipschitz. In particular, in this case, the constant  $\Lambda^{(\Fc\Pc)}(F)<+\infty$.
\par Thus, if $\rho$ is a metric on a finite set $\Mc$, and $\lambda\in[\Lambda^{(\Fc\Pc)}(F),+\infty)$, the $(\lambda;F)$-Iterative Algorithm \reff{ALG-IA} always produces the outcome {\bf ``Success''}. Cf. Remark \reff{PR2-RM}, part (i).\rbx}
\end{remark}

\begin{theorem}\lbl{A-L-IA} Let $\Mf=\MR$ be a finite pseudometric space. Let $\Tf$ be either the family $\KRT$ or the family $\HPL$, and let $F:\Mc\to\Tf$ be a set-valued mapping.
\par If Algorithm \reff{ALG-IA} produces the outcome {\bf ``No go''} (i.e., if $\Lambda^{(\Fc\Pc)}(F)=+\infty$, see {\bf STEP 1}), the set-valued mapping $F$ does not have a Lipschitz selection.
\smsk
\par  Otherwise, given a constant $\lambda\in[\Lambda^{(\Fc\Pc)}(F),+\infty)$, the $(\lambda;F)$-Iterative Algorithm produces the outcome {\bf ``Success''} and returns the mapping $f^{[\lambda;F]}:\Mc\to\RT$ defined by formula \rf{F-I2}. This mapping has the following properties:
\smsk
\par ($\bigstar \Acr$) The mapping $f^{[\lambda;F]}$ is well defined. This means the following: (i) For every $x\in\Mc$, the sets $F^{[1]}[x:\lambda]$ and
$F^{[2]}[x:\lambda,3\lambda]$ are non-empty; (ii) The rectangle $\Pi_{\lambda,F}(x)=\Hc[\Fcr[x:\lambda]]$ (see \rf{FR-2} and \rf{FR-21}) is a non-empty bounded subset of $\RT$.
\smsk
\par ($\bigstar \Bcr$) $f^{[\lambda;F]}$ is a Lipschitz selection of $F$ with the Lipschitz seminorm $\|f^{[\lambda;F]}\|_{\Lip(\Mc)}\le \gamma\lambda$ with $\gamma=420$.
\end{theorem}
\par {\it Proof.} Suppose that $\Lambda^{(\Fc\Pc)}(F)=+\infty$, and at the same time there exists a Lipschitz selection $f$ of $F$. Then, for every subset $\Mc'\subset\Mc$ with $\#\Mc'\le 4$, the mapping $f|_{\Mc'}$ is a Lipschitz selection of the restriction $F|_{\Mc'}$ of $F$ to $\Mc'$ with $\|f|_{\Mc'}\|_{\Lip(\Mc')}\le \lambda$ where $\lambda=\|f\|_{\Lip(\Mc)}$. Therefore, thanks to Definition \rf{FP-4}, $\lambda\in \Lc^{(\Fc\Pc)}(F)$ so that, thanks to \rf{LFP-W}, $\Lambda^{(\Fc\Pc)}(F)<\infty$, a contradiction.
\par This proves that if Algorithm \reff{ALG-IA} produces the outcome {\bf ``No go''}, the set-valued mapping $F$ does not have a Lipschitz selection.
\smsk
\par Now, let us assume that $\Lambda^{(\Fc\Pc)}(F)<\infty$
and $\lambda\in[\Lambda^{(\Fc\Pc)}(F),+\infty)$. Prove that in this case the mapping $f^{[\lambda;F]}:\Mc\to\RT$ defined by \rf{F-I2} has properties ($\bigstar \Acr$) and ($\bigstar \Bcr$).
\par Prove part (i) of the statement ($\bigstar \Acr$). Thanks to Lemma \reff {L-FP4}, the constant $\Lambda^{(\Fc\Pc)}(F)\in\Lc^{(\Fc\Pc)}(F)$, see Definition \reff{FP-4}. Therefore, $\lambda\in\Lc^{(\Fc\Pc)}(F)$ as well, so that, thanks to Definition \reff{FP-4}, for every subset $\Mc'\subset\Mc$ with $\#\Mc'\le 4$, the restriction $F|_{\Mc'}$ of $F$ to $\Mc'$ has a Lipschitz selection $f_{\Mc'}$ with $\|f_{\Mc'}\|_{\Lip(\Mc')}\le \lambda$.
\par We also recall that $\Mc$ is finite. Therefore, the pseudometric space $\Mf=\MR$ and the constant $\lambda$ satisfy the hypothesis of Theorem \reff{SPR-3}. This theorem tells us that property \rf{F2-NE} holds, i.e.,
the set
\bel{F2-NEM}
F^{[2]}[x:\lambda,3\lambda]\ne\emp~~~~\text{for every}~~~~ x\in\Mc.
\ee
\par But, thanks to \rf{IN-S}, $F^{[2]}[x:\lambda,3\lambda]\subset F^{[1]}[x:\lambda]$ so that $F^{[1]}[x:\lambda]\ne\emp$ as well. Thus, part (i) of ($\bigstar \Acr$) holds.
\par Let us note that part (ii) of ($\bigstar \Acr$) is immediate from part (i). Indeed, thanks to \rf{F2-NEM}, the set $\Fcr[x:\lambda]$ defined by \rf{FR-2} and \rf{FR-21} is a non-empty. Futhermore, $\Fcr[x:\lambda]$ is {\it bounded}. Therefore, its rectangular hull, the rectangular $\Pi_{\lambda,F}(x)$, is also non-empty and bounded. In particular, the center of $\Pi_{\lambda,F}(x)$ is well defined so that the mapping $f^{[\lambda;F]}$ is well defined on $\Mc$. See \rf{F-I2}.
\smsk
\par Let us prove property ($\bigstar \Bcr$). Let
\bel{GCR-F2}
\Gcr(x)=F^{[2]}[x:\lambda,3\lambda],~~~x\in\Mc.
\ee
Then, thanks to \rf{SP-3D},
\bel{FXY-DH1}
\dhf(\Gcr(x),\Gcr(y))\le 15\lambda\,\rho(x,y)
~~~~\text{for every}~~~~ x,y\in\Mc.
\ee
\par In these settings, definitions \rf{FR-2} and \rf{FR-21} read as follows:
\bel{FR-2N}
\Fcr[x:\lambda]=\Gcr(x)\cap Q(0,2r_x) ~~~~~~\text{where}~~
r_x=\dist\left(0,\Gcr(x)\right).
\ee
\par Let us estimate the Hausdorff distance between the sets $\Fcr[x:\lambda]$ and $\Fcr[y:\lambda]$. We will do this with the help of the following result:
\par {\it Let $C_1,C_2\subset\RT$ be convex sets. Given $i=1,2$, let $a_i\in\RT$, $r_i\ge 0$, and let $Q(a_i,r_i)$
be a square in $\RT$, see \rf{SQ-DF}. Suppose that $C_i\cap Q(a_i,r_i)\ne\emp$, $i=1,2$. Then}
\bel{HD-O2}
\dhf(C_1\cap Q(a_1,2r_1),C_2\cap Q(2a_2,r_2))
\le 14\,(\dhf(C_1,C_2)+\|a_1-a_2\|+|r_1-r_2|).
\ee
See \cite[Theorem 4]{PY-1995}, and also \cite[Lemma 3.9]{S-2004}.
\par We note, that, thanks to \rf{FR-2N},
$$
\Gcr(x)\cap Q(0,r_x)\ne\emp~~~~\text{and}~~~~\Gcr(y)\cap Q(0,r_y)\ne\emp.
$$
Therefore, thanks to \rf{HD-O2},
$$
\dhf(\Fcr[x:\lambda],\Fcr[y:\lambda])=
\dhf(\Gcr(x)\cap Q(0,2r_x),\Gcr(y)\cap Q(0,2r_y))
\le 14\,(\dhf(\Gcr(x),\Gcr(y))+|r_x-r_y|).
$$
\par Furthermore, from \rf{FXY-DH1}, we have
$$
|r_x-r_y|=|\dist(0,\Gcr(x))-\dist(0,\Gcr(y))|\le \dhf(\Gcr(x),\Gcr(y))
\le 15\lambda\,\rho(x,y).
$$
\par This inequality and \rf{FXY-DH1} imply the following:
\bel{FC-HD1}
\dhf(\Fcr[x:\lambda],\Fcr[y:\lambda])\le 14\,(15\lambda\,\rho(x,y)+15\lambda\,\rho(x,y))=
\gamma\lambda\,\rho(x,y)
\ee
with $\gamma=420$.
\smsk
\par We are in a position to prove that $f^{[\lambda;F]}$ is a Lipschitz selection of $F$ with Lipschitz seminorm at most $\gamma\lambda$. This property easily follows from \rf{FC-HD1} and the following simple claim proven in \cite[Remark 7.1]{FS-2017}: {\it For every compact convex set $S\subset\RT$, the center of its rectangular hull belongs to $S$, i.e., $\cent(\Hc[S])\in S$. Furthermore, for every compact convex sets $S_1,S_2\subset\RT$, we have}
\bel{CL-2}
\|\cent(\Hc[S_1])-\cent(\Hc[S_2])\|\le \dhf(S_1,S_2).
\ee
\par See also inequality \rf{CAB-H}.
\smsk
\par Now, thanks to this claim and definition \rf{F-I2}, for every $x\in\Mc$, we have
$$
f^{[\lambda;F]}(x)=\cent(\Pi_{\lambda,F}(x))=
\cent(\Hc[\Fcr[x:\lambda]])\in \Fcr[x:\lambda]\subset
F^{[2]}[x:\lambda,3\lambda]\subset F(x)
$$
proving that $f^{[\lambda;F]}$ is a {\it selection} of $F$.
\par Finally, combining \rf{F-I2} with inequality \rf{CL-2} (where we set $S_1=\Fcr[x:\lambda]$ and $S_2=\Fcr[y:\lambda]$) and inequality \rf{FC-HD1}, we conclude that
$\|f^{[\lambda;F]}\|_{\Lip(\Mc)}\le\gamma\lambda$ with $\gamma=420$.
\par The proof of the theorem is complete.\bx

\begin{remark} {\em Formula \rf{F-I2} can be simplified essentially provided
\bel{F2-BN}
\text{for some $\brx\in\Mc$ the set}~~~\Gcr(\brx)=F^{[2]}[\brx:\lambda,3\lambda]~~~
\text{is bounded.}
\ee
See \rf{GCR-F2}. In this case, thanks to \rf{FXY-DH1},
{\it every set $\Gcr(x)=F^{[2]}[x:\lambda,3\lambda]$, $x\in\Mc$, is bounded.}

\par This enables us to modify {\bf STEP 3} of Algorithm \reff{ALG-IA} as follows:
\msk
\par {\bf STEP 3{$^\prime$}.} We define the mapping $f^{[\lambda;F]}:\Mc\to\RT$ by the formula
$$
f^{[\lambda;F]}(x)=\cent(\,\Vc_{\lambda,F}(x)), ~~~~~~x\in\Mc,
$$
where $\Vc_{\lambda,F}(x)=\Hc[\Gcr(x)]
=\Hc[F^{[2]}[x:\lambda,3\lambda]]$. See Fig. 20.}
\end{remark}

\begin{figure}[h!]
\hspace{20mm}
\includegraphics[scale=0.31]{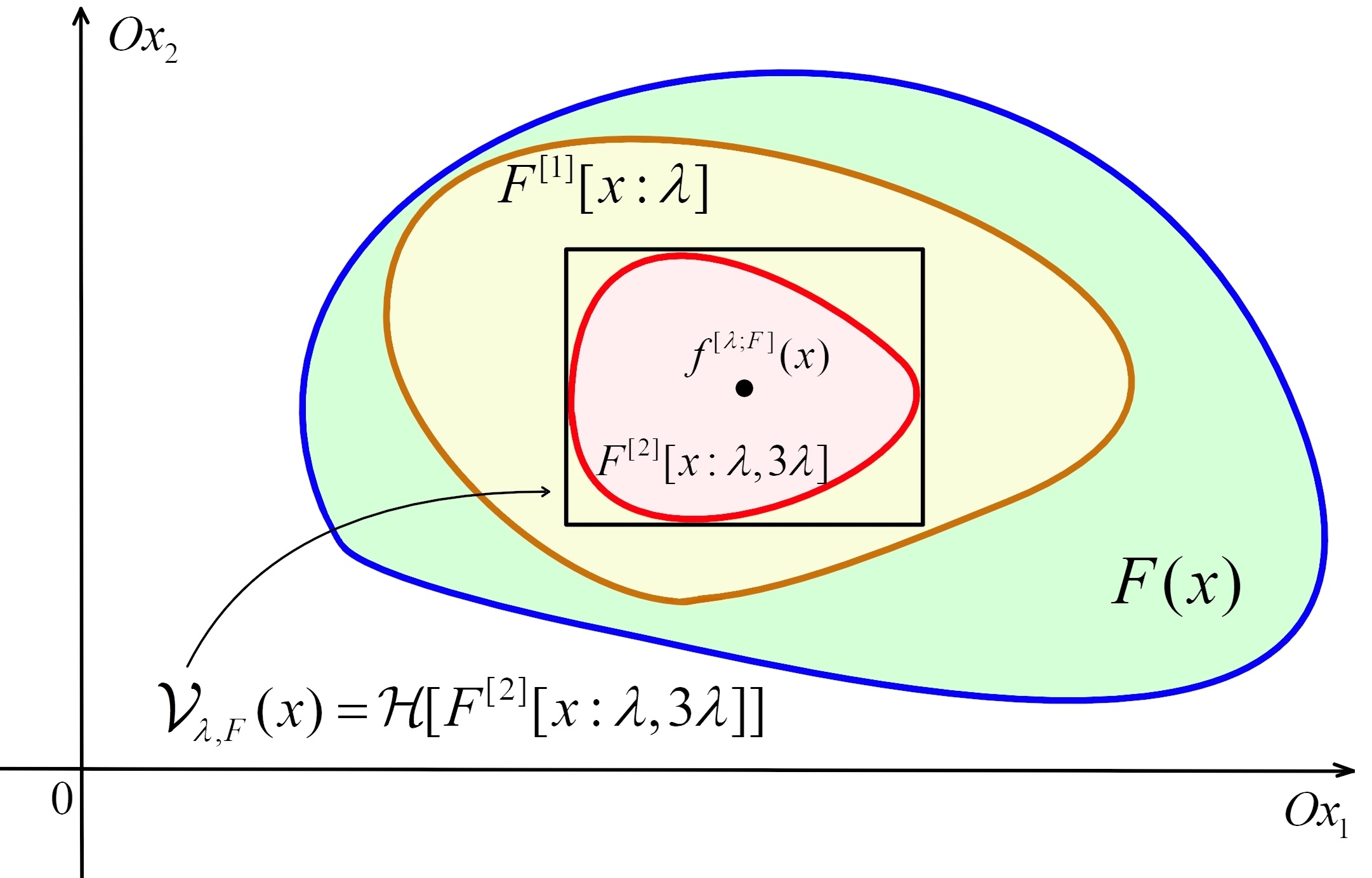}
\caption{The $\LF$-Iterative Algorithm for bounded sets $F^{[2]}[x:\lambda,3\lambda]$.}
\end{figure}

\par Because $\Gcr(x)$ is bounded for every $x\in\Mc$, the set $\Vc_{\lambda,F}(x)=\Hc[\Gcr(x)]$ is a {\it bounded} rectangle so that the mapping $f^{[\lambda;F]}$ is well defined on $\Mc$.
\par Following the scheme of the proof of Theorem \reff{A-L-IA}, we show that $f^{[\lambda;F]}$ is a selection of $F$. Furthermore, thanks to \rf{FXY-DH1} and \rf{CL-2}, we have
$$
\|f^{[\lambda;F]}(x)-f^{[\lambda;F]}(y)\|
=\|\cent(\Hc[\Gcr(x)])-\cent(\Hc[\Gcr(y)])\|\le \dhf(\Gcr(x),\Gcr(y))\le 15\lambda\rho(x,y)
$$
proving that
\bel{15-LM}
\|f^{[\lambda;F]}\|_{\Lip(\Mc)}\le 15 \lambda.
\ee
\par We note that condition \rf{F2-BN} holds for every set-valued mapping $F$ from $\Mc$ into the family $\KRT$. Therefore, in this important case, we can apply this simplified version of Algorithm \reff{ALG-IA} which produces a Lipschitz selection of $F$ with Lipschitz seminorm satisfying inequality \rf{15-LM}.\rbx

\begin{remark} {\em We recall Definitions \reff{DF-LR} and \reff{FP-4} of the constants $\Lambda_{\Rc}(F)$ and $\Lambda^{(\Fc\Pc)}(F)$ respectively. As we have noted in Remark \reff{TWO-ALG}, in general, given a set-valued mapping $F:\Mc\to\CRT$, the precise calculation of the constants $\Lambda_{\Rc}(F)$ and $\Lambda^{(\Fc\Pc)}(F)$ is a rather difficult technical problem. Nevertheless, comparing these constants, we note that the constant $\Lambda_{\Rc}(F)$ is defined in a more
constructive way than $\Lambda^{(\Fc\Pc)}(F)$. Moreover, Remark \reff{N3-E} tells us there exists an algorithm which given $F:\Mc\to\HPL$ calculates $\Lambda_{\Rc}(F)$ using at most $O(N^3)$ computer operations. (Here $N=\#\Mc$.)
\par Thus, in applications, it would be preferable to express Algorithm \reff{ALG-IA} in terms of the constant
$\Lambda_{\Rc}(F)$ rather than $\Lambda^{(\Fc\Pc)}(F)$.
\par An appropriate choice of the constant $\lambda$ at {\bf STEP 2} of this algorithm enables us to do this. Indeed, suppose we have calculated $\Lambda_{\Rc}(F)$ up to a constant $\eta\ge 1$. From inequalities \rf{LM-DT} and \rf{NF-MLM2}, we have $\Lambda^{(\Fc\Pc)}(F)\le 5\Lambda_{\Rc}(F)$. From this inequality it follows that the parameter
\bel{LRC-A}
\lambda=5\eta\Lambda_{\Rc}(F)
\ee
belongs to the interval $[\Lambda^{(\Fc\Pc)}(F),+\infty)$.
\par Therefore, thanks to Theorem \reff{A-L-IA}, if $\Lambda_{\Rc}(F)<\infty$, the $(\lambda;F)$-Iterative Algorithm \reff{ALG-IA} produces the outcome {\bf ``Success''} and returns a Lipschitz selection of $F$ with Lipschitz seminorm at most $\gamma\lambda=5\gamma\eta\Lambda_{\Rc}(F)$. (Recall that $\gamma=420$).\rbx}
\end{remark}


\par Let $\Mf=\MR$ be a finite pseudometric space, $N=\#\Mc$, and let $F:\Mc\to\HPL$ be a set-valued mapping. Let us see, how with the help of the $(\lambda;F)$-Iterative Algorithm \reff{ALG-IA} with $\lambda$ defined by \rf{LRC-A} one can produce a nearly optimal Lipschitz selection of $F$ using at most $O(N^3)$ computer operations.
\par First of all, Remark \reff{N3-E} tells us that there exists an algorithm for calculating the constant $\Lambda_{\Rc}(F)$ with $O(N^3)$ running time. Therefore, we can calculate the constant $\lambda$ from \rf{LRC-A} using at most $O(N^3)$ computer operations.
\par Let us estimate the running time of the $(\lambda;F)$-Iterative Algorithm \reff{ALG-IA}.

\begin{proposition}\lbl{RT-LIA} Let $F:\Mc\to\HPL$ be a set valued mapping defined on a finite pseudometric space $\Mf=\MR$. Suppose that the constant $\Lambda^{(\Fc\Pc)}(F)<\infty$.
\par If $\lambda\in [\Lambda^{(\Fc\Pc)}(F),+\infty)$, then
the $(\lambda;F)$-Iterative Algorithm \reff{ALG-IA} produces the outcome {\bf ``Success''} and returns the mapping $f^{[\lambda;F]}:\Mc\to\RT$ (see \rf{F-I2}) with $\|f^{[\lambda;F]}\|_{\Lip(\Mc)}\le \gamma\lambda$ using at most $O(N^2)$ computer operations. (Here $\gamma$ is an absolute constant.)
\end{proposition}
\par The proof of this proposition relies on a formula for the set $F^{[2]}[x:\lambda_1,\lambda_2]$, see \rf{F2-X}, which we present in Lemma \reff{RPR-F2} below.
\par Before stating this lemma, let us recall the definitions of some objects from the Projection Algorithm \reff{PA-DS}. These are the set $F^{[1]}[x:\lambda_1]$, i.e., the $\lambda_1$-balanced refinement of $F$ defined by \rf{F-LM1}, its rectangular hull, the set $\Tc_{F,\lambda_1}(x)=\Hc[F^{[1]}[x:\lambda_1]]$, see \rf{TC-DF}, and, finally, the $\lambda_2$-balanced refinement of $\Tc_{F,\lambda_1}$, i.e., the rectangle $\Tc^{[1]}_{F,\lambda_1}[x:\lambda_2]$ defined by
formula \rf{RT-1}.
\begin{lemma}\lbl{RPR-F2} Let $\Mf=\MR$ be a finite pseudometric space and let $F:\Mc\to\CRT$. Given constants $\lambda_1,\lambda_2\ge 0$, $\lambda_1\le \lambda_2$, let us assume that, for each $x\in\Mc$, the sets  $F^{[1]}[x:\lambda_1]$ and $F^{[2]}[x:\lambda_1,\lambda_2]$ are non-empty. Then,
\bel{F2-R3}
F^{[2]}[x:\lambda_1,\lambda_2]=F^{[1]}[x:\lambda_1]
\cap \Tc^{[1]}_{F,\lambda_1}[x:\lambda_2]~~~~\text{for every}~~~x\in\Mc.
\ee
\end{lemma}
\par {\it Proof.} We note that
$F^{[2]}[x:\lambda_1,\lambda_2]\subset
F^{[1]}[x:\lambda_1]$, see \rf{IN-S}. Furthermore,
$$
\Tc_{F,\lambda_1}(z)=\Hc[F^{[1]}[z:\lambda_1]]\supset F^{[1]}[z:\lambda_1],~~~z\in\Mc,
$$
so that
$$
\Tc^{[1]}_{F,\lambda_1}[x:\lambda_2]=
\bigcap_{z\in\Mc}
\big\{\Tc_{F,\lambda_1}(z)+\lambda_2\rho(x,z)Q_0\big\}
\supset
\bigcap_{z\in\Mc}
\big\{F^{[1]}[z:\lambda_1]+\lambda_2\rho(x,z)Q_0\big\}.
$$
From this and definition \rf{LB-2}, we have
\bel{KH-CQ}
\Tc^{[1]}_{F,\lambda_1}[x:\lambda_2]
\supset
F^{[2]}[x:\lambda_1,\lambda_2]
\ee
proving that the right hand side of \rf{F2-R3} contains its left hand side.
\par Let us prove the converse statement. We will need two auxiliary results. Here is the first of them: {\it Let $C_1,C_2\in\CRT$, and let $r\ge 0$. Suppose that $C_1\cap C_2\ne\emp$. Then}
\bel{C12-N}
C_1\cap C_2+rQ_0=(C_1+rQ_0)\cap(C_2+rQ_0)
\cap\Hc[C_1\cap C_2+rQ_0].
\ee
\par The second auxiliary result states the following: {\it Let $\Kc\subset\CRT$ be a finite family of closed convex sets in $\RT$ with non-empty intersection, and let $r\ge 0$. Then}
\bel{K-NL}
\left(\,\bigcap_{K\in\Kc} K\right)+rQ_0=
\bigcap_{K,K'\in\Kc} \left\{(K\cap K')+rQ_0\right\}.
\ee
\par These result were proved in \cite{S-2022} under the conditions $C_1,C_2\in\KRT$ and $\Kc\subset\KRT$ respectively. (See \cite[Lemma 3.4, Lemma 5.4]{S-2022}.) Obvious changes in the proofs of these results show that these property are satisfied for every $C_1,C_2\in\CRT$ and any finite family $\Kc\subset\CRT$.
\par Property \rf{K-NL} tells us that, for every $z\in\Mc$, we have
\be
F^{[1]}[z:\lambda_1]+\lambda_2\rho(x,z)Q_0&=&
\left(\bigcap_{z\in\Mc}
\big\{F(u)+\lambda_1\rho(z,u)Q_0\big\}\right)
+\lambda_2\rho(x,z)Q_0
\nn\\
&=&
\bigcap_{u,u'\in\Mc}
\left[\{F(u)+\lambda_1\rho(z,u)Q_0\}\cap
\{F(u')+\lambda_1\rho(z,u')Q_0\}
+\lambda_2\rho(x,z)Q_0\right]\nn\\
&=&
\bigcap_{u,u'\in\Mc} A(x:z,u,u')\nn
\ee
where
$$
A(x:z,u,u')=
\{F(u)+\lambda_1\rho(z,u)Q_0\}\cap
\{F(u')+\lambda_1\rho(z,u')Q_0\}
+\lambda_2\rho(x,z)Q_0.
$$
\par Thanks to property \rf{C12-N}, for every $z,u,u'\in\Mc$, we have
$$
A(x:z,u,u')=
\{F(u)+(\lambda_1\rho(z,u)+\lambda_2\rho(x,z))Q_0\}\cap
\{F(u')+(\lambda_1\rho(z,u')+\lambda_2\rho(x,z))Q_0\}
\cap \Hc[A(x:z,u,u')].
$$
\par Recall that $0\le\lambda_1\le\lambda_2$. From this, the triangle inequality and \rf{LB-1}, we have
\be
S_1&=&
F(u)+(\lambda_1\rho(z,u)+\lambda_2\rho(x,z))Q_0\nn\\
&\supset&
F(u)+\lambda_1(\rho(z,u)+\rho(x,z))Q_0
\supset
F(u)+\lambda_1\rho(x,u)Q_0\supset F^{[1]}[x:\lambda_1].\nn
\ee
In the same way we prove that
$$
S_2=
F(u')+(\lambda_1\rho(z,u')+\lambda_2\rho(x,z))Q_0\supset F^{[1]}[x:\lambda_1].
$$
\par Let us show that the set $\Hc[A(x:z,u,u')]$ contains $\Tc^{[1]}_{F,\lambda_1}[x:\lambda_2]$. Indeed, thanks to \rf{LB-1},
$$
G=\{F(u)+\lambda_1\rho(z,u)Q_0\}\cap
\{F(u')+\lambda_1\rho(z,u')Q_0\}
\supset F^{[1]}[z:\lambda_1],
$$
so that $\Hc[G]\supset\Hc[F^{[1]}[z:\lambda_1]]=
\Tc_{F,\lambda_1}(z)$. From this and \rf{N-HS}, we have
\be
S_3&=&
\Hc[A(x:z,u,u')]=\Hc[G+\lambda_2\rho(x,z)Q_0]=
\Hc[G]+\lambda_2\rho(x,z)Q_0\nn\\
&\supset&
\Tc_{F,\lambda_1}(z)+\lambda_2\rho(x,z)Q_0\supset \Tc^{[1]}_{F,\lambda_1}[x:\lambda_2]\nn
\ee
proving the required inclusion $S_3\supset \Tc^{[1]}_{F,\lambda_1}[x:\lambda_2]$.
\par These properties of the sets $S_1,S_2$ and $S_3$ imply the following:
$$
A(x:z,u,u')=S_1\cap S_2\cap S_3\supset F^{[1]}[x:\lambda_1]
\cap \Tc^{[1]}_{F,\lambda_1}[x:\lambda_2]~~~\text{for every}~~~z,u,u'\in\Mc.
$$
\par Thus, for every $z\in\Mc$, we have
$$
F^{[1]}[z:\lambda_1]+\lambda_2\rho(x,z)Q_0=
\bigcap_{u,u'\in\Mc} A(x:z,u,u')\supset
F^{[1]}[x:\lambda_1]\cap \Tc^{[1]}_{F,\lambda_1}[x:\lambda_2].
$$
Thanks to this property and definition \rf{F2-X}, we get
$$
F^{[2]}[x:\lambda_1,\lambda_2]=
\bigcap_{z\in\Mc}
\big\{F^{[1]}[z:\lambda_1]+
\lambda_2\,\rho(x,z)\,Q_0\big\}\supset
F^{[1]}[x:\lambda_1]\cap \Tc^{[1]}_{F,\lambda_1}[x:\lambda_2]
$$
proving that the left hand side of \rf{F2-R3} contains its right hand side.
\par The proof of the lemma is complete.\bx
\msk

\begin{remark} {\em Let us note the following interesting representation of the mapping $F^{[2]}[\cdot:\lambda_1,\lambda_2]$: in the settings of Lemma \reff{RPR-F2}, for every $x\in\Mc$, the following is true
\bel{N-F2-R3}
F^{[2]}[x:\lambda_1,\lambda_2]=F^{[1]}[x:\lambda_1]
\cap \Tc^{[1]}_{F,\lambda_1}[x:\lambda_2]=
F^{[1]}[x:\lambda_1]
\cap \Hc[F^{[2]}[x:\lambda_1,\lambda_2]].
\ee
\par This is immediate from Lemma \reff{RPR-F2} and inclusion \rf{KH-CQ}. Indeed, because $\Tc^{[1]}_{F,\lambda_1}[x:\lambda_2]$ is a {\it rectangle}, from \rf{KH-CQ} we have
$\Tc^{[1]}_{F,\lambda_1}[x:\lambda_2]
\supset \Hc[F^{[2]}[x:\lambda_1,\lambda_2]]$.
This and \rf{F2-R3} imply the second equality in \rf{N-F2-R3}.
\smsk
\par Note that, in general,
$$
\Tc^{[1]}_{F,\lambda_1}[x:\lambda_2]
\supsetneqq\Hc[F^{[2]}[x:\lambda_1,\lambda_2]]
~~~\text{on}~~~\Mc.
$$
However, thanks to \rf{N-F2-R3}, intersection of these rectangles with $F^{[1]}[x:\lambda_1]$ produces the same set, i.e., the set $F^{[2]}[x:\lambda_1,\lambda_2]$.
\rbx}
\end{remark}

\par {\it Proof of Proposition \reff{RT-LIA}.} We know that
$\Lambda^{(\Fc\Pc)}(F)<\infty$ and $\lambda\in [\Lambda^{(\Fc\Pc)}(F),+\infty)$. Therefore, thanks to Theorem \reff{A-L-IA}, the $(\lambda;F)$-Iterative Algorithm \reff{ALG-IA} produces the outcome {\bf ``Success''} and returns the mapping $f^{[\lambda;F]}:\Mc\to\RT$ defined by formula \rf{F-I2}. Its Lipschitz seminorm is bounded by $\gamma\lambda$ where $\gamma$ is an absolute constant.
\par Thanks to part ($\bigstar \Acr$) of Theorem \reff{A-L-IA}, in this case, the sets $F^{[1]}[x:\lambda]$ and $F^{[2]}[x:\lambda,3\lambda]$ are non-empty for every $x\in\Mc$. Furthermore, Lemma \reff{RPR-F2} tells us that in this case, for every $x\in\Mc$, the following equality
\bel{P-WN}
F^{[2]}[x:\lambda,3\lambda]=F^{[1]}[x:\lambda]
\cap \Tc^{[1]}_{F,\lambda}[x:3\lambda]
\ee
holds. Recall that, given $x\in\Mc$,
\bel{T1-F}
\Tc^{[1]}_{F,\lambda}[x:3\lambda]=
\bigcap_{z\in\Mc}
\big\{\Tc_{F,\lambda}(z)+3\lambda\,\rho(x,z)\,Q_0\big\}
\ee
and $\Tc_{F,\lambda}(x)=\Hc[F^{[1]}[x:\lambda]]$.
\smsk
\par Let us see how the mapping $f^{[\lambda;F]}$ can be constructed in $O(N^2)$ running time.
\par At {\it the first step} of this procedure, for every $x\in\Mc$, we construct the rectangle $\Tc_{F,\lambda}(x)$. Because $F:\Mc\to\HPL$ and $\#\Mc=N$, the convex polygon $F^{[1]}[x:\lambda]$
is determined by $O(N)$ linear constrains. See \rf{LB-1}.
\par Part (i) of Corollary \reff{CR-HD} tells us that there exists an algorithm which, for every $x\in\Mc$, constructs the set $\Tc_{F,\lambda}(x)$, the rectangular hull of $F^{[1]}[x:\lambda]$, in $O(N)$ running time.
\par Thus, at this step we are able to construct all $N$ rectangles $\Tc_{F,\lambda}(x)$ in $O(N^2)$ running time.
\smsk
\par At {\it the second step} of the procedure, following \rf{T1-F}, we construct each rectangle $\Tc^{[1]}_{F,\lambda}[x:3\lambda]$, $x\in\Mc$, in $O(N)$ running time. Thus, the constructing all $N$ rectangles $\Tc^{[1]}_{F,\lambda}[x:3\lambda]$ will require at most $O(N^2)$ computer operations.
\smsk
\par We turn to {\it the third step} of the procedure.
Thanks to representation \rf{P-WN}, each set  $F^{[2]}[x:\lambda,3\lambda]$ is a convex polygon determined by at most $N+4$ linear constraints. Part (ii) of Corollary \reff{CR-HD} tells us that there exists an algorithm which, for every $x\in\Mc$, calculates the distance from $O$ to $F^{[2]}[x:\lambda,3\lambda]$ (i.e., the quantity $r_x$, see \rf{FR-21}) using at most $O(N)$ computer operations. Thus, calculating all $N$ numbers $r_x$, $x\in\Mc$, can be proceed in $O(N^2)$ running time.
\smsk
\par At {\it the fourth step} of the procedure, we treat the sets $\Fcr[x:\lambda]=F^{[2]}[x:\lambda,3\lambda]\cap Q(0,2r_x)$, $x\in\Mc$, see \rf{FR-2}. Because each  $F^{[2]}[x:\lambda,3\lambda]$ is determined by at most $N+4$ linear constraints, the set $\Fcr[x:\lambda]$ is a convex polygon determined by at most $N+8$ linear constraints. Therefore, thanks to part (i) of Corollary \reff{CR-HD}, there exists an algorithm which, for every $x\in\Mc$, constructs the set $\Pi_{\lambda,F}(x)$, the rectangular hull of $\Fcr[x:\lambda]$ (see \rf{PI-H}), in $O(N)$ running time.
\par Thus, at this step we construct all $N$ rectangles
$\Pi_{\lambda,F}(x)$ using at most $O(N^2)$ computer operations.
\smsk
\par Finally, at {\it the fifth step} of the procedure, for every $x\in\Mc$, we construct the point $f^{[\lambda;F]}(x)$ as {\it the center of the rectangle} $\Pi_{\lambda,F}(x)$. See formula \rf{F-I2}. Because the number of all such rectangles is $N$, this procedure will
take $4N$ computer operations.
\smsk
\par Thus, we construct the mapping $f^{[\lambda;F]}:\Mc\to\RT$ in five steps, and each of these steps requires at most $O(N^2)$ computer operations.
This shows that $f^{[\lambda;F]}$ can be constructed in $O(N^2)$ running time completing the proof of Proposition \reff{RT-LIA}.\bx
\newpage
\fontsize{12}{13.5}\selectfont

\end{document}